\documentclass[12pt]{amsart}
\usepackage{graphicx,amsmath,amscd,amsfonts,amsthm,amssymb,fullpage}
\usepackage{times}
\usepackage{enumerate}
\usepackage[all,cmtip]{xy}
\usepackage{diagxy}

\newtheorem{thm}{Theorem}[section]
\newtheorem{cor}[thm]{Corollary}
\newtheorem{lem}[thm]{Lemma}

\newtheorem{prop}[thm]{Proposition}

\theoremstyle{definition}
\newtheorem{defin}[thm]{Definition}
\newtheorem{rem}[thm]{Remark}

\numberwithin{equation}{section}

\newcommand{\R}{\mathbb R}
\newcommand{\C}{\mathbb C}

\newcommand{\Z}{\mathbb Z}

\newcommand{\p}{\mathfrak p}
\newcommand{\g}{\mathfrak g}

\newcommand{\gk}{\mathfrak k}
\newcommand{\ga}{\mathfrak a}

\newcommand{\gn}{\mathfrak n}
\newcommand{\gu}{\mathfrak u}
\newcommand{\gh}{\mathfrak h}
\newcommand{\gl}{\mathfrak l}

\newcommand{\gm}{\mathfrak m}

\newcommand{\bm}{{\bf m}}

\newcommand{\cO}{\mathcal O}

\newcommand{\cA}{\mathcal A}

\newcommand{\cF}{\mathcal F}

\newcommand{\cS}{\mathcal S}

\newcommand{\cM}{\mathcal M}
\newcommand{\cN}{\mathcal N}

\newcommand{\cC}{\mathcal C}

\newcommand{\cB}{\mathcal B}

\newcommand{\cJ}{\mathcal J}

\newcommand{\tD}{\widetilde{\Delta}}
\newcommand{\tb}{\widetilde{\beta}}
\newcommand{\tE}{\widetilde{\Sigma}}

\newcommand{\q}{\mathfrak q}

\newcommand{\Ll}{\mathcal L}

\newcommand{\Ad}{\textup{Ad}}

\newcommand{\Vol}{\textup{Vol}}

\newcommand{\be}{\begin{equation}}
\newcommand{\ee}{\end{equation}}
\newcommand{\bes}{\begin{equation*}}
\newcommand{\ees}{\end{equation*}}

\newcommand{\ba}{\begin{eqnarray}}
\newcommand{\ea}{\end{eqnarray}}
\newcommand{\bas}{\begin{eqnarray*}}
\newcommand{\eas}{\end{eqnarray*}}

\begin{document}

\title{$L^p$ Norms of Higher Rank Eigenfunctions and Bounds for Spherical Functions}

\author{Simon Marshall}
\address{S. Marshall: Department of Mathematics,
University of Wisconsin -- Madison,
480 Lincoln Drive,
Madison,
WI 53706, USA;}
\email{marshall@math.wisc.edu}

\begin{abstract}
We prove almost sharp upper bounds for the $L^p$ norms of eigenfunctions of the full ring of invariant differential operators on a compact locally symmetric space, as well as their restrictions to maximal flat subspaces.  Our proof combines techniques from semiclassical analysis with harmonic theory on reductive groups, and makes use of new asymptotic bounds for spherical functions that are of independent interest.
\end{abstract}

\maketitle


\section{Introduction}

If $M$ is a compact Riemannian manifold of dimension $n$ and $\psi$ is a Laplace eigenfunction on $M$ satisfying $\Delta \psi = \lambda^2 \psi$, it is a well studied problem to investigate the asymptotic behaviour of the $L^p$ norms of $\psi$ as $\lambda \rightarrow \infty$.  The fundamental upper bound for these norms was established by Sogge \cite{So} (see also Avacumovi\'c \cite{Av} and Levitan \cite{Le} in the case $p = \infty$), who proves that

\begin{equation}
\label{sogge}
\| \psi \|_p \ll \lambda^{\delta(n,p)} \| \psi \|_2
\end{equation}
where $\delta(n,p)$ is the piecewise linear function of $1/p$ given by

\begin{equation}
\label{delta0}
\delta(n,p) = \bigg\{ \begin{array}{ll} n( \tfrac{1}{2} - \tfrac{1}{p} ) -1/2, & 0 \le \tfrac{1}{p} \le \tfrac{n-1}{2(n+1)}, \\
\tfrac{n-1}{2}( \tfrac{1}{2} - \tfrac{1}{p} ), & \tfrac{n-1}{2(n+1)} \le \tfrac{1}{p} \le \tfrac{1}{2}.
\end{array}
\end{equation}
Moreover, these bounds were shown by Sogge \cite{So} to be sharp when $M$ is the round $n$-sphere $S^n$.\\

It is sometimes possible to improve the upper bound in (\ref{sogge}) by assuming that $M$ has additional symmetry, or that $\psi$ is an eigenfunction of extra differential operators that commute with $\Delta$.  In the extreme case of the flat torus $T^n$, for instance, if one assumes that $\psi$ is an eigenfunction of all the translations $\{ i \partial / \partial x_j \}$ then $\psi$ is a complex exponential, and so we have $\| \psi \|_p \le C \| \psi \|_2$ for all $p$ and some $C$ depending only on $T^n$.  A more interesting example of this phenomenon is given by Sarnak in his letter to Morawetz \cite{Sa}.  He proves that if $X$ is a compact locally symmetric space of dimension $n$ and rank $r$, and $\psi$ is an eigenfunction of the full ring of differential operators on $X$ with Laplace eigenvalue $\lambda^2$, then

\begin{equation}
\label{Sarnak}
\| \psi \|_\infty \ll \lambda^{(n-r)/2} \| \psi \|_2.
\end{equation}
(Notations are standard and given in $\mathsection$\ref{Lpnotation}.)  Note that (\ref{Sarnak}) represents an improvement in the exponent of (\ref{sogge}) from $(n-1)/2$ to $(n-r)/2$.  This upper bound is also sharp in the case when $X$ is of compact type, and Sarnak states that it should be considered as the `local bound' for the sup norm of a higher rank eigenfunction.\\

The goal of this paper is to derive the correct local bound for all $L^p$ norms of an eigenfunction in higher rank, by combining real interpolation with an analysis of spherical functions.  Our main result in this direction is stated below, which in the compact case differs from the sharp bound only by a factor of $( \log t )^{1/2}$ at the kink point.

\begin{thm}
\label{main}
Let $X$ be a compact locally symmetric space of dimension $n$ and rank $r$ that is a quotient of the globally symmetric space $S = G / K$, and assume that $S$ is irreducible and not Euclidean.  Let $\ga_0$ be a real Cartan subalgebra of $G$, and let $\ga_0^*$ and $\ga^*$ be its real and complex dual respectively.  If $f \in C^\infty(X)$ is an eigenfunction of the ring of invariant differential operators, we say $f$ has spectral parameter $\nu \in \ga^*$ if it has the same eigenvalues as the function $\exp( (\rho + i\nu)(A(x)))$ on $S$.

Let $B^* \subset \ga^*_0$ be a compact set that is bounded away from the singular set.  Let $\psi \in C^\infty(X)$ be an eigenfunction of the full ring of invariant differential operators on $X$, with $\| \psi \|_2 = 1$ and spectral parameter $t \lambda$ where $t > 0$ and $\lambda \in B^*$.  We have

\begin{equation}
\label{mainfn}
\| \psi \|_p \ll_{B^*, p} \bigg\{ \begin{array}{ll} (\log t)^{1/2} t^{r \delta(n/r, p)}, & p = \tfrac{2(n+r)}{n-r}, \\
t^{r \delta(n/r, p)}, & p \neq \tfrac{2(n+r)}{n-r}, \end{array}
\end{equation}
where the function $\delta$ is as in (\ref{delta0}).  Moreover, these bounds are sharp up to the logarithmic factor in the case when $X$ is of compact type.

\end{thm}

A similar result was obtained in the Euclidean case by Mockenhaupt \cite{M}.  It will be apparent in the course of the proof of Theorem \ref{main} that when $X$ is the quotient of a product $S = S_1 \times \ldots \times S_d$ of irreducible symmetric spaces, the $L^p$ norm of an eigenfunction on $X$ is bounded by the product of the functions (\ref{mainfn}) for each irreducible factor of $S$.  Moreover, in the compact case this will again be sharp up to the logarithmic factors at the kink points.\\

\begin{figure}
\begin{center}
\includegraphics[width=330pt]{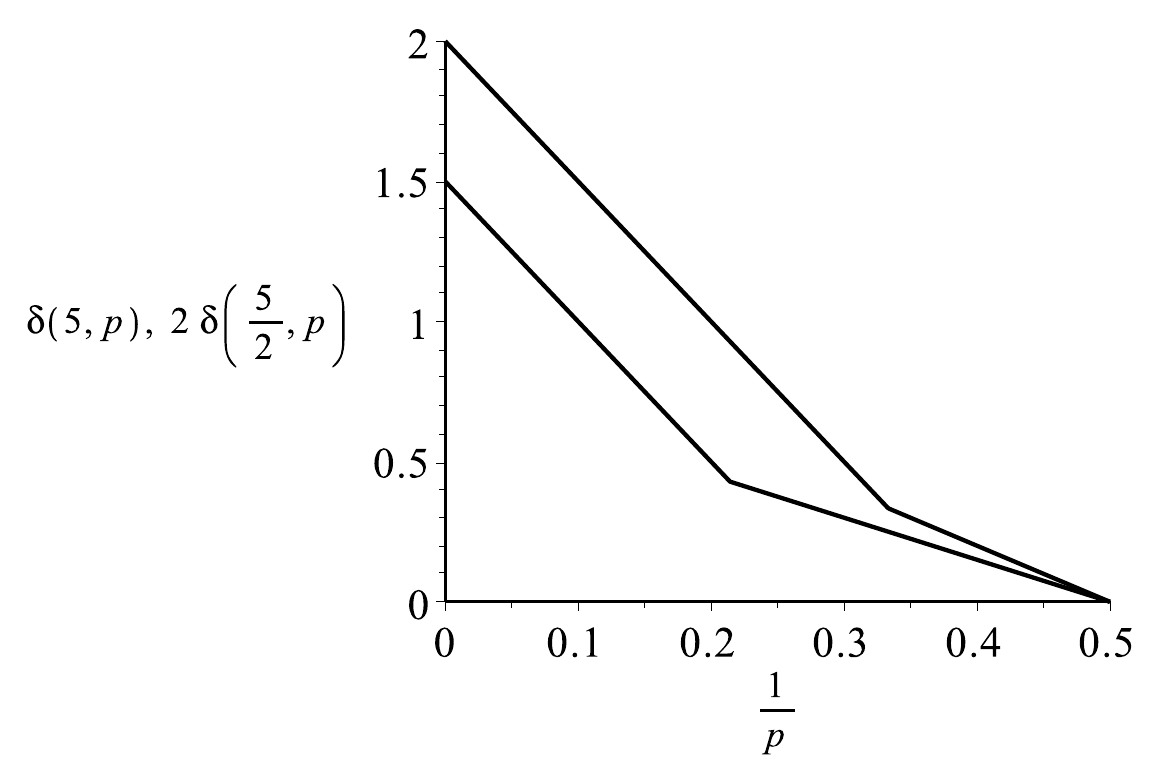}
\caption{Comparison of the two exponents $\delta(5,p)$ and $2\delta(5/2,p)$ appearing in Sogge's bound and Theorem \ref{main} in the case $G = SL(3,\R)$.}
\label{deltas}
\end{center}
\end{figure}

To give an example comparing the bound produced by Theorem \ref{main} with the classical bound (\ref{sogge}), let $X$ be a quotient of the globally symmetric space $SL(3,\R) / SO(3)$.  It was proven by Selberg \cite{Se} that the ring $R$ of invariant differential operators on $X$ is isomorphic to the free polynomial ring $\C[\Delta, D]$, where $D$ is an operator of degree 3.  Let $\psi$ be an eigenfunction of $R$, and assume that the spectral parameter of $\psi$ is restricted as in Theorem \ref{main}.  The two exponents $\delta(5,p)$ and $2\delta(5/2,p)$ appearing in Sogge's bound and Theorem \ref{main} are graphed together in Figure \ref{deltas}.  We see that by using the symmetry of $X$ in the form of its extra differential operators, we are able to significantly strengthen the bounds for $\| \psi \|_p$.

Let us take a moment to discuss the significance of the exponent in Theorem \ref{main}, and hopefully convince the reader that it is natural.  Suppose that $r | n$, and let $X$ be a product of $r$ compact manifolds $X_1 \times \ldots \times X_r$ of dimension $n/r$.  Let $\Delta_i$ be the Laplacian of $X_i$, and let $\psi = \psi_1 \times \ldots \times \psi_r$ be a joint eigenfunction of the Laplacians $\Delta_i$ on $X$.  Let $\Delta_i \psi = \lambda_i^2 \psi$, and assume that the ratios $\lambda_i / \lambda_j$ are all bounded by some constant.  By applying Sogge's bound (\ref{sogge}) to each $\psi_i$, we may show that

\begin{equation*}
\| \psi \|_p \ll \lambda^{r \delta(n/r, p)} \| \psi \|_2,
\end{equation*}
where $\lambda^2 = \lambda_1^2 + \ldots + \lambda_r^2$.  We may therefore summarize Theorem \ref{main} by saying that, from the point of view of the convex bound for $L^p$ norms of eigenfunctions, a locally symmetric space of dimension $n$ and rank $r$ whose universal cover is irreducible behaves like the product of $r$ general Riemannian manifolds of dimension $n/r$.

It would be interesting to know in which other cases this product behaviour occurs, that is when the $L^p$ bounds of Theorem \ref{main} hold for a more general compact manifold $M$ of dimension $n$ with $r$ commuting differential operators that are `independent' in some sense.  There are no nontrivial examples of this in the completely integrable case, as it was proven by Toth and Zelditch \cite{TZ} that if $M$ is a quantum completely integrable manifold and all joint eigenfunctions on $M$ are uniformly bounded then $M$ is a flat torus.\\

In proving Theorem \ref{main}, we shall in fact show that the same bounds hold for the $L^2 \rightarrow L^p$ norm of a spectral projector onto a ball of fixed radius about $\lambda$.  With this formulation, our bounds will be sharp up to the $\log$ in the case of both compact and noncompact type.  The fact that this bound is sharp for individual eigenfunctions in the compact case is due to the high multiplicity of the spectrum, so that by choosing the radius of our spectral projector to be sufficiently small we know that it will always pick out exactly one eigenvalue of high multiplicity.

In both cases, the bounds of Theorem \ref{main} are realised by simple wave packets which are the higher rank analogues of the zonal functions and Gaussian beams on a general Riemannian manifold.  We shall describe these packets on the globally symmetric space $S = G/K$, their analogues on $X$ being similar.  The cotangent bundle $T^*S$ of $S$ is isomorphic to the $K$-principal bundle $G \times_K \p^*$, which we recall is defined to be the quotient of the trivial bundle $G \times \p^*$ by the action

\begin{equation*}
(g, v)k = (gk, \Ad_k^{-1} v).
\end{equation*}
If $\lambda \in \ga^*$, we define $T^*_\lambda S \subset T^*S$ by

\bes
T^*_\lambda S = \{ (g, v) \in G \times_K \p^* | v \in \Ad_K (\lambda / \| \lambda \|) \}.
\ees
Saying that $\psi \in C^\infty(S)$ is an approximate eigenfunction of the ring of invariant differential operators on $S$ with parameter $\lambda$ then implies that the microlocal support of $\psi$ is concentrated on $T^*_\lambda S$; see \cite[$\mathsection$5.4]{SV}.

Let $o \in S$ correspond to the identity coset of $K$, and let $A$ be a maximal flat subspace containing $o$.  Define $T^*_\lambda A = A \times \lambda / \| \lambda \| \subset T^*_\lambda S$, and let $\Ll_\lambda \subset T^*_\lambda S$ be the orbit of $T^*_\lambda A$ under rotation by $K$ about $o$.  The $K$-biinvariant functions $k_t$ constructed in $\mathsection$\ref{Lp} and $\mathsection$\ref{cpctLp} saturate the $L^p$ norms on $S$ for $p$ above the kink point, and we believe that these functions should be microlocally concentrated on $\Ll_\lambda$.  The fibre of the projection map $\pi : \Ll_\lambda \rightarrow S$ at $s \in S$ can be identified with $\text{Stab}_K(s)$, so that this fiber is identified with $M$ for generic $s$ and with $K$ at $s = o$, and correspondingly $\psi$ will be strongly peaked at $o$ so that we may think of $\psi$ as an analogue of the usual zonal function on a Riemannian manifold.  Note that in the case of compact type we can prove that the spherical functions $\varphi_\lambda$ also saturate the $L^p$ bounds of Theorem \ref{main} for large $p$.

For $p$ below the kink point, the $L^p$ norms on $S$ are saturated by the higher rank analogue of a Gaussian beam, which is simply a wave packet concentrated on a maximal flat subspace, and whose microlocal support is concentrated on the set $T^*_\lambda A$.  These functions will be described more thoroughly in the case of compact type in $\mathsection$\ref{compactbeam}.\\

The methods we develop to prove Theorem \ref{main} also allow us to deduce the following result on the restrictions of eigenfunctions to flats in $X$.  We hope to extend this theorem to more general locally symmetric submanifolds in future.

\begin{thm}
\label{restrict}

With notations as in Theorem \ref{main}, let $E$ be an open ball in a maximal flat subspace of $X$.

\begin{enumerate}[(a)]

\item If $n > 3r$, the $L^p$ norms of $\psi|_E$ satisfy

\begin{equation*}
\| \psi|_E \|_p \ll_{B^*} t^{(n-r)/2 - r/p}.
\end{equation*}

\item If $n = 3r$, the $L^p$ norms of $\psi|_E$ satisfy

\begin{align*}
\| \psi|_E \|_p & \ll_{B^*,p} t^{(n-r)/2 - r/p}, \quad p > 2, \\
\| \psi|_E \|_2 & \ll_{B^*} (\log t)^{1/2} t^{n/2 - r}.
\end{align*}

\item If $n < 3r$, the $L^p$ norms of $\psi|_E$ satisfy

\bes
\| \psi|_E \|_p \ll_{B^*, p} \bigg\{ \begin{array}{ll} (\log t)^{1/2} t^{\delta(p)}, & p = \tfrac{4r}{n-r}, \\
t^{\delta(p)}, & p \neq \tfrac{4r}{n-r}, \end{array}
\ees
where $\delta(p)$ is the piecewise linear function

\bes
\delta(p) = \bigg\{ \begin{array}{ll} n - r - 2r/p, & 0 \le \tfrac{1}{p} \le \tfrac{n-r}{4r}, \\
(n-r)/2, & \tfrac{n-r}{4r} \le \tfrac{1}{p} \le \tfrac{1}{2}.
\end{array}
\ees

\end{enumerate}

Moreover, all of these bounds are sharp up to the logarithmic factor in the case of compact type.

\end{thm}

When $r = 1$, this is a slight weakening of a theorem of Burq, G\'erard and Tzvetkov \cite{BGT}.  We note that there are only finitely many globally symmetric spaces that fall under cases (b) and (c) of Theorem \ref{restrict}.  In case (b), these are the spaces associated to $SO(3,1)$, $SO(3,2)$, $SO(3,3)$, $SL_4(\R)$, and their compact duals, and in case (c) these are the spaces associated to $SL_2(\R)$, $SL_3(\R)$, and their compact duals.  Theorem \ref{restrict} will be proven in $\mathsection$\ref{flat}.

\subsection{Asymptotics for Spherical Functions}

In the course of proving Theorem \ref{main} we have found it necessary to develop sharp asymptotics for spherical functions of large eigenvalue on $G$, which we state here as separate theorems.  First let us assume that $G$ is semisimple and noncompact with finite center.  For $\lambda \in \ga^*_0$, let $\varphi_\lambda$ denote the standard spherical function with parameter $\lambda$, normalised so that $\varphi_\lambda(e) = 1$.  If $\alpha$ is a nonzero root of $\ga$ in $\g$, let $m(\alpha)$ denote its multiplicity.  Our result is the following:

\begin{thm}
\label{phithm}

Let $B \subset \ga_0$ and $B^* \subset \ga^*_0$ be compact sets, with $B^*$ bounded away from the singular set.  We have the upper bound

\begin{equation}
\label{phibd}
\varphi_{t\lambda}( \exp(H) ) \ll_{B,B^*} \prod_{\alpha \in \Delta^+} ( 1 + t |\alpha(H)| )^{-m(\alpha)/2}
\end{equation}
for $H \in B$ and $\lambda \in B^*$.

\end{thm}

Theorem \ref{phithm} is the strongest upper bound that can be given for $\varphi_{t\lambda}( \exp(H))$ when $H$ and $\lambda$ are bounded and $t$ grows, at least under the regularity assumption on $\lambda$ that we have made.  We have attempted to remove this assumption, but so far only have an approach to this in the case of rank 2.  We hope to carry this out in a future paper, and to use it to remove the regularity condition in Theorem \ref{main} in some cases.

Theorem \ref{phithm} is similar to results of Duistermaat, Kolk, and Varadarajan \cite[Corollary 9.3 and Theorem 11.1]{DKV}, and Blomer and Pohl \cite[Theorem 2]{BP}.  The result of Blomer and Pohl gives a bound for $\varphi_{t\lambda}( \exp(H) )$ which is not generally sharp, but which is uniform as $H$ and $\lambda$ vary in any compact subsets of $\ga_0$ and $\ga_0^*$.  The results of Duistermaat, Kolk, and Varadarajan are only uniformly sharp if $H$ is restricted to a compact equisingular set, but \cite[Theorem 11.1]{DKV} is uniformly sharp for $\lambda$ in any compact set.  In some sense, \cite[Theorem 11.1]{DKV} is complementary to Theorem \ref{phithm}, which requires $\lambda$ to be regular but is uniformly sharp in $H$.  Our proof of Theorem \ref{phithm} is similar to the proof of \cite[Theorem 11.1]{DKV}, with the main difference being that the phase function $\phi(k,H,\lambda)$ that appears in the oscillatory integrals is linear in $\lambda$, but nonlinear in the variable $H$ that we are allowing to degenerate.  Theorem \ref{phithm} will be derived from an analysis of stationary phase integrals in $\mathsection$\ref{phibounds}.

Our methods also allow us to strengthen the asymptotic formula for $\varphi_{t\lambda}(\exp(H))$ given in equation (9.10) of \cite{DKV}.  Let $\Vol_0(K)$ and $\Vol_0(M)$ be the volumes of $K$ and $M$ with respect to the metric induced from minus the Killing form on $\gk$, and for any $w \in W$ define

\begin{equation}
\label{index}
\sigma_w(H, \lambda) = - \sum_{\alpha \in \Delta^+} m(\alpha) \text{sgn}( \langle \lambda, \alpha \rangle \alpha(wH) ).
\end{equation}

\begin{thm}
\label{phiasympthm}

Let $\ga_r$ and $\ga_r^*$ denote the regular sets in $\ga_0$ and $\ga_0^*$ respectively.  Let $B \subset \ga_0$ and $B^* \subset \ga_r^*$ be compact sets.  If $H \in \ga_0$, let $\| H \|_s$ denote the Killing distance from $H$ to the singular set.  There are functions $f_w \in C^\infty(\ga_r \times \ga_r^* \times \R_{>0})$ for $w \in W$ such that

\bes
\left( \frac{\partial}{\partial H} \right)^a f_w(H, \lambda, t)  \ll_{B, B^*, a} \frac{1}{t \| H \|_s^{a+1}} \prod_{\alpha \in \Delta^+} (t |\alpha(H)|)^{-m(\alpha)/2},
\ees
and 

\begin{multline}
\label{phiasymp}
\varphi_{t \lambda}( \exp(H) ) = \prod_{\alpha \in \Delta^+} \left| \frac{ \langle \alpha, t\lambda \rangle}{ 2\pi} \sinh \alpha(H) \right|^{-m(\alpha)/2} \frac{ \Vol_0(M) }{ \Vol_0(K) } \sum_{w \in W} \exp(it \lambda(wH) + i \pi \sigma_w(H, \lambda) / 4) \\
+ \sum_{w \in W} \exp(it \lambda(wH)) f_w(H, \lambda,t) + O_{B, B^*,A}( (t \| H \|_s)^{-A}) \prod_{\alpha \in \Delta^+} (t |\alpha(H)|)^{-m(\alpha)/2}
\end{multline}
for $H \in B \cap \ga_r$ and $\lambda \in B^*$.

\end{thm}

We also have the following asymptotic, which is weaker than Theorem \ref{phiasympthm} but seems to be the most useful for our planned applications.

\begin{thm}
\label{phiasympthm2}

Let $\ga_r$ and $\ga_r^*$ denote the regular sets in $\ga_0$ and $\ga_0^*$ respectively.  Let $B \subset \ga_0$ and $B^* \subset \ga_r^*$ be compact sets.  If $H \in \ga_0$, let $\| H \|_s$ denote the Killing distance from $H$ to the singular set.  There are functions $f_w \in C^\infty(\ga_r \times \ga_r^* \times \R_{>0})$ for $w \in W$ such that

\bes
\left( \frac{\partial}{\partial H} \right)^a f_w(H, \lambda, t)  \ll_{B, B^*, a} \| H \|_s^{-a} \prod_{\alpha \in \Delta^+} (t |\alpha(H)|)^{-m(\alpha)/2},
\ees
and 

\begin{multline*}
\varphi_{t \lambda}( \exp(H) ) =  \sum_{w \in W} \exp(it \lambda(wH)) f_w(H, \lambda,t) + O_{B, B^*,A}( (t \| H \|_s)^{-A}) \prod_{\alpha \in \Delta^+} (t |\alpha(H)|)^{-m(\alpha)/2}
\end{multline*}
for $H \in B \cap \ga_r$ and $\lambda \in B^*$.

\end{thm}

We have a result analogous to Theorem \ref{phithm} in the case of compact type, but which is weakened by the requirement that the group variable be constrained to a small ball about the origin.  Let $U$ be a compact semisimple Lie group, and $K$ a subgroup with the property that $(U,K)$ is a Riemannian symmetric pair.  If $\mu$ is a spherical weight (defined in $\mathsection$\ref{compactnotn}), we let $\varphi_\mu$ be the $K$-spherical function on $U$ with parameter $\mu$, normalised so that $\varphi_\mu(e) = 1$.

\begin{thm}
\label{cpctphibd}
There exists a ball $B \subset i \ga_0$ about the origin such that for all compact sets $B^* \subset \ga_0^*$ that are bounded away from the singular set, we have

\begin{equation*}
\varphi_{t\mu}( \exp(H) ) \ll_{B^*} \prod_{\alpha \in \Delta^+} ( 1 + t | \alpha(H)|)^{-m(\alpha)/2}
\end{equation*}
for $H \in B$ and $\mu \in B^*$.

\end{thm}

Theorem \ref{cpctphibd} will be proved in $\mathsection$\ref{compactsph}.

\section{Bounds for $L^p$ norms in noncompact type}
\label{Lp}

We shall first prove Theorem \ref{main} in the case when $X$ is of noncompact type.  The proof in the case of compact type is similar, and we shall make the modifications to our argument that are required to treat it in $\mathsection$\ref{cpctLp}.

\subsection{Notation}
\label{Lpnotation}

\subsubsection{Lie algebras}
\label{sec211}

We shall denote real Lie algebras with a subscript $0$, and denote their complexifications by dropping this subscript.  Let $G$ be a connected noncompact semisimple real Lie group with finite center and Lie algebra $\g_0$.  In $\mathsection$\ref{Lp}--\ref{flat} we shall further assume that $G$ is almost simple, in the sense that $\g_0$ is simple over $\R$, or that $G$ does not factor after an isogeny.  Note that we shall only use this assumption when summing the bounds we obtain for truncated kernels.  We denote the Killing form on $\g$ by $\langle \, , \rangle$.  Let $\g_0 = \gk_0 + \p_0$ be a Cartan decomposition of $\g_0$, and $\theta$ the corresponding Cartan involution.  Let $K$ be the compact connected subgroup of $G$ with Lie algebra $\gk_0$, so that $S = G/K$ is a globally symmetric space of noncompact type.  Let

\begin{align*}
G & = NAK, \quad g = n(g) \exp( A(g)) k(g) \\
\g & = \gk + \ga + \gn
\end{align*}
be an Iwasawa decomposition of $G$.  Let $M'$ and $M$ be the normaliser and centraliser of $\ga$ in $K$, let $\gm$ be the Lie algebra of $M$, and let $W$ be the Weyl group $M' / M$.  We let $\Delta$ denote the set of roots of $\g$ with respect to $\ga$.  If $\alpha \in \Delta$ we denote the corresponding root space by $\g_\alpha$.

\begin{rem}

Note that we shall include 0 in $\Delta$, which is not standard notation, but it will be convenient for us.  In particular, $\g_\alpha = \ga + \gm$ when $\alpha = 0$.  To avoid confusion with the real Lie algebra $\g_0$, the expression $\g_\alpha$ with $\alpha = 0$ will only appear implicitly when we index over root spaces.

\end{rem}

We let $m(\alpha) = \dim \g_\alpha$ when $\alpha \neq 0$, and when $\alpha = 0$ we let $m(\alpha) = \dim \gm$.  We let $\Delta^+$ be a choice of positive roots, to which we associate the nilpotent subalgebra $\gn = \sum_{\alpha \in \Delta^+} \g_\alpha$ and closed positive Weyl chamber $\ga^+_0$.  We let $\ga^*_{0,+}$ denote the dual positive Weyl chamber.  We define $\rho = \frac{1}{2} \sum_{\alpha \in \Delta^+} m(\alpha) \alpha$.  If $\nu \in \ga^*$, $H_\nu \in \ga$ will be the vector dual to $\nu$ under the Killing form.  We emphasise the following piece of notation, as it is nonstandard and will be used frequently.

\begin{defin}

We let $\widetilde{\Delta}$ denote the multiset on $\Delta$ in which every $\alpha \in \Delta$ appears with multiplicity $m(\alpha)$.  If $R \subseteq \Delta$, we let $\widetilde{R}$ denote the corresponding subset of $\widetilde{\Delta}$.

\end{defin}

\subsubsection{The Harish-Chandra transform}
\label{sec212}

If $\lambda \in \ga^*$, we let $\varphi_\lambda$ denote the spherical function with parameter $\lambda$, defined by

\bes
\varphi_\lambda(x) = \int_K \exp( (\rho + i\lambda)(A(kx)) ) dk.
\ees
If $f \in C^\infty_0(S)$, we define its Harish-Chandra transform by

\bes
\widehat{f}(\nu) = \int_S f(x) \varphi_{-\nu}(x) dx.
\ees
If $f$ is $K$-biinvariant, we have the inversion formula

\bes
f(x) = \int_{\ga_0^* / W} \widehat{f}(\nu) \varphi_\nu(x) |c(\nu)|^2 d\nu,
\ees
where $c(\nu)$ is Harish-Chandra's $c$-function.  See \cite{GV} for more information about this transform.

\subsection{An Outline of the Proof}
\label{Lpoutline}

We shall assume that $B^*$ is contained in the positive dual Weyl chamber $\ga^*_{0,+}$.  We shall approach Theorem \ref{main} by the standard method of constructing a family of approximate spectral projectors $T_t$ onto a ball of radius 1 about $t \lambda$, and bounding the norms of $T_t$ from $L^2$ to $L^p$.  Note that all bounds we state will depend on $X$ and $B^*$ from now on, but will be uniform in $\lambda \in B^*$.

We shall construct $T_t$ using the Harish-Chandra transform, which will allow us to gain good control over the behaviour of the integral kernel of this operator.  Choose a function $h \in \cS( \ga^*_0 )$ of Paley-Wiener type that is real-valued and $\ge 1$ in a ball of radius 1 about the origin.  Let

\begin{equation*}
h_t(\nu) = \sum_{w \in W} h( w\nu + t\lambda),
\end{equation*}
and let $k^0_t$ be the $K$-biinvariant function on $S$ with Harish-Chandra transform $h_t$.  It is of compact support independent of $t\lambda$ by the Paley-Wiener theorem of \cite{Ga}.  Define $K^0_t$ to be the point pair invariant kernel on $S$ associated to $k^0_t$, given by $K_t^0(x,y) = k_t^0(x^{-1} y)$ for $x, y \in G$.  Let $T_t$ be the operator on $X$ with integral kernel

\be
\label{gammasum}
T_t(x,y) = \sum_{\gamma \in \Gamma} K^0_t(x, \gamma y).
\ee
As $T_t \psi = h_t( - t\lambda) \psi$ and $h_t(-t\lambda) \ge 1$, it will suffice to prove bounds for $\| T_t \psi \|_p$ of the form (\ref{mainfn}), uniformly for $\lambda \in B^*$.  As is common, we shall approach this by forming the adjoint square operator $T_t T_t^*$ and proving the bounds

\begin{equation}
\label{mainop}
\| T_t T_t^* f \|_p \ll_p \bigg\{ \begin{array}{ll} \log t \times t^{2 \delta(p)} \| f \|_{p'}, & p = \tfrac{2(n+r)}{n-r}, \\
t^{2 \delta(p)} \| f \|_{p'}, & p \neq \tfrac{2(n+r)}{n-r}, \end{array}
\end{equation}
for the operator norms of $T_t T_t^*$ using real interpolation.  Here, $p'$ and $p$ are dual exponents and $f \in C^\infty(X)$.  Note that $T_t$ is actually self-adjoint because $h_t$ is real, and so if we define the $K$-biinvariant function $k_t = k^0_t * k^0_t$, then $T_t T_t^*$ is associated to $k_t$ as in (\ref{gammasum}).  We define $B \subset \ga_0$ to be a ball about the origin such that $\text{supp}( k_t \circ \exp ) \subseteq B$.

\subsection{The Case of Rank One}
\label{Lprankone}

We begin by outlining the real interpolation argument used to prove (\ref{mainop}) when $r=1$, in which case it consists of a dyadic decomposition of $k_t$ in terms of its radial support.  Choose $g \in C^\infty_0(\R)$ to be a real, non-negative, even function that is identically 1 in a neighbourhood of 0, and for $m \in \Z_{\ge 0}$ let

\begin{equation*}
f_{t,m}(x) = \bigg\{ \begin{array}{ll} g(t x), & m = 0, \\
g( t e^{-m} x) - g( t e^{-m+1} x), & m > 0.
\end{array}
\end{equation*}
Fix an isomorphism between $\ga_0$ and $\R$, and pull the functions $f_{t,m}$ back to $\beta_{t,m}$ on $\ga_0$.  Define the $K$-biinvariant function $k_{t,m}$ by $k_{t,m}(\exp(H)) = \beta_{t,m}(H) k_t(\exp(H))$, let $K_{t,m}$ be the associated point pair invariant, and $T_{t,m}$ the integral operator on $X$ associated to $K_{t,m}$.  It may be shown that

\begin{align}
\label{rk1dyad1}
\| T_{t,m} f \|_\infty & \ll t^{n-1} e^{-m(n-1)/2} \| f \|_1, \\
\label{rk1dyad2}
\| T_{t,m} f \|_2 & \ll t^{-1} e^m \| f \|_2.
\end{align}
By interpolating between (\ref{rk1dyad1}) and (\ref{rk1dyad2}) we may prove the bound 

\begin{equation}
\label{rk1dyad3}
\| T_{t,m} f \|_p \ll t^{n(1-2/p) - 1} \exp\left( m \left( \tfrac{n+1}{p} - \tfrac{n-1}{2} \right) \right) \| f \|_{p'}, \quad 2 \le p \le \infty,
\end{equation}
and because the supports of $k_t$ were uniformly compact there is $C > 0$ such that $k_{t,m} = 0$ for $m > \log t + C$.  Summing over $m$ then gives

\begin{equation*}
\| T_t T_t^* f \|_p \ll t^{n(1-2/p) - 1} \sum_{0 \le m \le \log t + C} \exp\left( m \left( \tfrac{n+1}{p} - \tfrac{n-1}{2} \right) \right) \| f \|_{p'}, \quad 2 \le p \le \infty.
\end{equation*}
The sum is a geometric progression of length $\log t$ with extremal terms $1$ and $t^{(n+1)/p - (n-1)/2}$.  The bounds of Theorem \ref{main} follow immediately from this and the observation that

\begin{equation*}
2 \delta(p) = n(1-2/p) - 1 + \max \left\{ 0, \tfrac{n+1}{p} - \tfrac{n-1}{2} \right\}.
\end{equation*}

\subsection{Partitions of Unity}
\label{Lppartitions}

Our proof for higher rank groups works by applying a similar decomposition in terms of the Cartan $\ga_0$ co-ordinate to $k_t$.  We begin by defining the partition of unity that we shall use.

Consider a partition of $\Delta$ into three sets $R_0$ and $R_\pm$, and define $C$ to be the cone

\begin{multline*}
C = \{ v \in \ga_0 | \alpha(v) = 0, \alpha \in R_0 \} \cap \{ v \in \ga_0 | \alpha(v) > 0, \alpha \in R_+ \} \\
\cap \{ v \in \ga_0 | \alpha(v) < 0, \alpha \in R_- \}.
\end{multline*}
We let $\cC$ be the collection of nonempty cones obtained in this way, which form a partition of $\ga_0$.  We choose a point $p_C$ in the interior of every cone $C$.  We define a flag to be a sequence $\{ C_0, C_1, \ldots, C_r \}$ of elements of $\cC$ such that $C_i \subset \overline{C_{i+1}}$ and $\dim C_i = i$, and let the set of flags be $\cF$.  If $F = \{ C_0, C_1, \ldots, C_r \} \in \cF$, and $1 \le i \le r$, define $\Delta_{i,F}$ to be the set of roots that vanish on $C_{i-1}$ but not $C_i$.  For every $F \in \cF$ we define the closed cone $S(F)$ to be the positive linear span of the set $\{ p_C | C \in F \}$, so that $\ga_0 = \bigcup_{F \in \cF} S(F)$.  We may assume without loss of generality that $w S(F) = S(wF)$ for all $w \in W$.

Let $F = \{ C_0, C_1, \ldots, C_r \} \in \cF$, and let $\phi_F$ be the linear isomorphism $\phi_F: \ga_0 \simeq \R^r$ such that $\phi_F(p_{C_i})$ is the vector with $i$ 1's followed by $r-i$ 0's.  We see that $\phi_F$ maps $S(F)$ onto the cone $S_0 = \{ x_1 \ge x_2 \ldots \ge x_r \ge 0 \}$, and that $\phi_F(C_i) \subset \{ (x_1, \ldots, x_i, 0, \ldots 0) | x_j \in \R \}$.  If $\alpha \in \Delta$ is a root, we let $\phi_F^* \alpha$ be the pushforward of $\alpha$ to $\R^r$.

\begin{lem}
\label{phishape}

Let $\{ e_j | 1 \le j \le r \}$ be the standard basis of $\R^r$.  We may choose the points $p_C$ so that for all $\alpha \in \Delta$ and $F \in \cF$, $\phi_F^* \alpha$ is either nonpositive or nonnegative on the positive quadrant $\R_+^r$, and if $\alpha \in \Delta_{i,F}$ we have $\phi_F^* \alpha(e_j) = 0$ iff $j < i$.

\end{lem}

\begin{proof}

Let $F = \{ C_0, C_1, \ldots, C_r \}$.  We may assume without loss of generality that $\alpha$ is positive on $C_r$, and let $\alpha \in \Delta_{i,F}$.  We define a new collection of points $p'_C \in C$ by setting $p'_C = A^{\dim C} p_C$ for some large $A > 1$ to be chosen later.  It is equivalent to show that our new collection of points satisfies the condition $\phi_F^* \alpha(e_j) \ge 0$ with equality iff $j < i$.  As $e_j = \phi_F(p'_{C_j}) - \phi_F(p'_{C_{j-1}})$, we have

\begin{align*}
\phi_F^* \alpha(e_j) & = \alpha(p'_{C_j}) - \alpha(p'_{C_{j-1}}) \\
& = A^j \alpha(p_{C_j}) - A^{j-1} \alpha(p_{C_{j-1}}).
\end{align*}
Our assumptions on $\alpha$ imply that $\alpha(p_{C_k}) \ge 0$ with equality iff $k < i$, and so by choosing $A$ large enough we see that the same will be true for $\phi_F^* \alpha(e_j)$.  As there are only finitely many choices for $F$ and $\alpha$, some $A$ will work for all of them.

\end{proof}

Define $\cM_t$ to be the set

\bes
\cM_t = \{ {\bf m} \in \Z^r | \log t + 1 \ge m_1 \ge m_2 \ldots \ge 0 \}.
\ees
Choose a small $\delta >0$.  We define an equivalence relation on $\cM_t$ by setting ${\bf m} \sim {\bf m}'$ if and only if $m_1 = m_1'$, and for all $i$ with $m_i \neq m_i'$ we have $\max \{ m_i, m_i' \} \le \delta m_1$.  If we define

\bes
\cM_{t,\delta} = \{ \bm \in \cM_t \, | \, m_i = 0 \text{ or } m_i > \delta m_1, \; \forall \, i \, \},
\ees
then $\cM_{t,\delta}$ contains a representative for every equivalence class in $\cM_t$.

Let $g \in C^\infty_0(\R)$ be a real-valued function supported in $[-e,e]$ such that $g(x) = 1$ for $x \in [-1,1]$, and both $g(x)$ and $g(x) - g(ex)$ are nonnegative.  For ${\bf m} \in \cM_{t,\delta}$ and $1 \le i \le r$, define $f_{{\bf m}, i} \in C^\infty_0(\R)$ by

\be
\label{fdef}
f_{{\bf m}, i}(x) = \bigg \{ \begin{array}{ll} g( t e^{ -\lfloor \delta m_1 \rfloor - i} x_i) & \text{if } m_i = 0, \\
\chi_{[0,\infty)}(x_i) [g( t e^{-m_i - i} x_i) - g( t e^{-m_i+1 - i} x_i)] & \text{if } m_i > \lfloor \delta m_1 \rfloor, \end{array}
\ee
and define $f_{\bf m} \in C^\infty_0(\R^r)$ by

\bes
f_{\bf m}(x) = \prod_{i=1}^r f_{{\bf m}, i}(x_i) \ge 0.
\ees
Define $\overline{S}_0 \subset S_0$ to be the set $\{ 1 \ge x_1 \ge x_2 \ldots \ge x_r \ge 0 \}$.

\begin{lem}
\label{partition}

We have

\bes
\sum_{ {\bf m} \in \cM_{t,\delta} } f_{\bf m}(x) = 1
\ees
when $x \in \overline{S}_0$.

\end{lem}

\begin{proof}

When $r=1$, the result is obvious.  Assume $r > 1$, and define

\begin{align*}
\cN_t & = \{ {\bf n} \in \Z^{r-1} | \log t + 1 \ge n_1 \ge n_2 \ldots \ge 0 \}, \\
\cN_{t,\delta} & = \{ {\bf n} \in \cN_t \, | \, n_i = 0 \text{ or } n_i > \delta n_1, \; \forall \, i \, \}.
\end{align*}
If $\bm \in \cM_{t,\delta}$, let $\overline{\bm} \in \cN_{t,\delta}$ be its first $r-1$ entries.  We may define the function $f_{{\bf n},i}$ for ${\bf n} \in \cN_{t\delta}$ and $1 \le i \le r-1$ as in (\ref{fdef}), and write

\begin{align}
\notag
\sum_{ {\bf m} \in \cM_{t,\delta} } f_{\bf m}(x) & = \sum_{ {\bf n} \in \cN_{t,\delta} } \sum_{ \substack{\bm \in \cM_{t,\delta} \\ \overline{\bm} = {\bf n} }} f_{\bf m}(x) \\
\label{fsum}
& = \sum_{ {\bf n} \in \cN_{t,\delta} } \prod_{i=1}^{r-1} f_{ {\bf n}, i}(x_i) \sum_{ \substack{\bm \in \cM_{t,\delta} \\ \overline{\bm} = {\bf n} }} f_{\bm,r}(x_r).
\end{align}
If $n_{r-1} > \delta n_1$, we have

\begin{align*}
\sum_{ \substack{\bm \in \cM_{t,\delta} \\ \overline{\bm} = {\bf n} }} f_{\bm,r}(x_r) & = g( t e^{ -\lfloor \delta n_1 \rfloor - r} x_r) + \sum_{ \delta n_1 < m_r \le n_{r-1}} g( t e^{-m_r - r} x_r) - g( t e^{-m_r+1 - r} x_r) \\
& = g( t e^{-n_{r-1} - r} x_r),
\end{align*}
while if $n_{r-1} = 0$ we have

\bes
\sum_{ \substack{\bm \in \cM_{t,\delta} \\ \overline{\bm} = {\bf n} }} f_{\bm,r}(x_r) = g( t e^{ -\lfloor \delta n_1 \rfloor - r} x_r).
\ees
We may assume without loss of generality that $f_{ {\bf n}, r}(x_{r-1}) \neq 0$ in (\ref{fsum}).  If $n_{r-1} > \delta n_1$, this implies that

\bes
0 \le t e^{-n_{r-1} - r +1} x_{r-1} \le e.
\ees
Our assumption that $x \in \overline{S}_0$ implies that $x_{r-1} \ge x_r \ge 0$, so that $0 \le t e^{-n_{r-1} - r} x_r \le 1$ and $g(t e^{-n_{r-1} - r} x_r) = 1$.  Likewise, when $n_{r-1} = 0$ we also have $g( t e^{ -\lfloor \delta n_1 \rfloor - r} x_r) = 1$.  Applying this to (\ref{fsum}) gives

\bes
\sum_{ {\bf m} \in \cM_{t,\delta} } f_{\bf m}(x) = \sum_{ {\bf n} \in \cN_{t,\delta} } \prod_{i=1}^{r-1} f_{ {\bf n}, i}(x_i),
\ees
and proceeding inductively completes the proof.

\end{proof}

We now pull the functions $f_{\bf m}$ back to $\ga_0$ under $\phi_F$, and let the collection of functions we obtain be $\{ f_{{\bf m}, F} | {\bf m} \in \cM_{t,\delta} \}$.  We may assume without loss of generality that the set of functions we generate in this way is invariant under the Weyl group, i.e. that $f_{ {\bf m}, wF }(wH) = f_{ {\bf m}, F }(H)$ for $w \in W$.  By scaling the points $p_C$ if necessary we may assume that $\phi_F(2B) \subseteq [-1,1]^r$ for all $F$, and it follows from this and Lemma \ref{partition} that

\bes
\sum_{F \in \cF} \sum_{ {\bf m} \in \cM_{t,\delta} } f_{ {\bf m}, F}(H) \ge 1
\ees
for $H \in 2B$.  If we choose $f_\infty$ to be a smooth Weyl-invariant function that vanishes on $B$ and is equal to $1$ outside $2B$, we then have

\be
\label{GH}
G(H) = f_\infty(H) + \sum_{F \in \cF} \sum_{ {\bf m} \in \cM_{t,\delta} } f_{ {\bf m}, F}(H) \ge 1
\ee
for all $H \in \ga_0$.  We define the partition of unity $\{ \beta_{ {\bf m}, F} | {\bf m} \in \cM_{t,\delta}, F \in \cF \} \cup \beta_\infty$ on $\ga_0$ by setting

\bes
\beta_{ {\bf m}, F}(H) = f_{ {\bf m}, F}(H) / G(H), \quad \beta_\infty(H) = f_\infty(H) / G(H).
\ees
We have introduced the parameter $\delta$ so that we may prove the following lemma, which will allow us to prove that the Harish-Chandra transforms of our truncated kernels decay near the walls of $\ga^*_{0,+}$.

\begin{lem}
\label{betabox}

If $\partial^\alpha$ is a product of derivatives in the co-ordinate directions on $\ga_0$, we have

\be
\label{betadiff}
\partial^\alpha \beta_{\bm, F} \ll_\alpha t^{|\alpha|} e^{-\delta m_1 |\alpha|}.
\ee

\end{lem}

\begin{proof}

Each of the functions $f_{\bm,F}$ clearly satisfies the bound (\ref{betadiff}), and because there is some $N >0$ independent of $t$ such that each $H \in \ga_0$ lies in the support of at most $N$ of the functions $f_{\bm,F}$, we see that the function $G(H)$ in (\ref{GH}) also satisfies (\ref{betadiff}).  The lemma follows from this and the bound $G(H) \ge 1$.

\end{proof}

\subsection{Bounds for Truncated Kernels}
\label{Lpbounds}

We shall now use our partition of unity to decompse the $K$-biinvariant function $k_t$, and give bounds for the norms of the operators constructed from the truncated pieces.  For $\bm \in \cM_{t,\delta}$ and $F \in \cF$, the function $\beta_{\bm,[F]}(H) = \sum_{w \in W} \beta_{\bm,F}( w H)$ is Weyl-invariant and so we may define a $K$-biinvariant function $\tb_{\bm, F}$ by setting

\bes
\tb_{\bm, F}( \exp(H)) = \beta_{\bm,[F]}(H)
\ees
for $H \in \ga_0$.  We then define $k_{\bm, F} = \tb_{\bm, F} k_t$.  Clearly $k_{\bm, F} = k_{\bm, wF}$, and the condition that $\beta_\infty$ vanishes on $B$ implies that

\bes
k_t = \sum_{F \in \cF / W} \sum_{ \bm \in \cM_{t,\delta} } k_{\bm, F}.
\ees
As before, we let $K_{\bm,F}$ and $T_{\bm,F}$ be the point pair invariant and integral operator associated to $k_{\bm, F}$.  Let $L(F,p)$ be the linear functional

\bes
L(F,p)(x) = ( 1/(2p) - 1/4) \sum_{i=1}^r x_i | \tD_{i,F} | + 2/p \sum_{i=1}^r x_i.
\ees
We shall require the following bounds on $T_{\bm, F}$.

\begin{prop}
\label{truncbound}

There is a constant $N$ depending only on $G$, and a constant $C_1$ depending on $\phi$, such that if we define $\chi(\bm)$ by

\bes
\chi(\bm) = \Big\{ \begin{array}{ll} 1 & \text{if} \quad m_r < \delta m_1 + C_1 \\ 0 & \text{otherwise}, \end{array}
\ees
then we have

\be
\label{Lpbd}
\| T_{\bm,F} f \|_p \ll_{\delta} t^{n(1 - 2/p) - r} \exp \left( L(F,p)(\bm) + \chi(\bm) N \delta m_1 \right) \| f \|_{p'}.
\ee
for all $p \ge 2$ and $f \in C^\infty(X)$.  The implied constant is uniform in ${\bf m}$.

\end{prop}

\begin{proof}[Proof of Proposition \ref{truncbound} assuming Theorem \ref{phithm}]

We begin by establishing the following bounds for the values taken by the roots on the support of $\beta_{\bm,F}$.

\begin{lem}
\label{rootbound}

If $\alpha \in \Delta_{i,F}$, we have

\be
\label{rootupper}
\sup \{ |\alpha(H)| | H \in \textup{supp}( \beta_{\bm, F}) \} \ll t^{-1} \max \{ e^{m_i}, e^{\delta m_1} \},
\ee
and there are positive constants $C_1$ and $C_2$ such that if $m_i \ge \delta m_1 + C_1$, we have

\be
\label{rootlower}
\inf \{ |\alpha(H)| | H \in \textup{supp}( \beta_{\bm, F}) \} \ge C_2 t^{-1} e^{m_i}.
\ee

\end{lem}

\begin{proof}

Let $H \in \text{supp}( \beta_{\bm, F})$, so that $x = \phi_F( H) \in \text{supp}(f_\bm)$ and $\alpha(H) = \phi_F^* \alpha(x)$, and assume without loss of generality that $\alpha$ is positive on $C_r \in F$.  We know that all $x \in \text{supp}(f_\bm)$ satisfy

\ba
\label{xibound1}
|x_i| \le t^{-1} e^{ \lfloor \delta m_1 \rfloor + r+1} \quad \text{if $m_i$ = 0}, \\
\label{xibound2}
t^{-1} e^{m_i -r-1} \le x_i \le t^{-1} e^{m_i+r+1} \quad \text{otherwise}.
\ea
By Lemma \ref{phishape}, if we let the standard basis vectors of $\R^r$ be $e_i$ as before, we have

\bes
\phi_F^* \alpha(x) = \sum_{j=i}^r \phi_F^* \alpha( x_j e_j).
\ees
This implies that

\begin{align*}
|\phi_F^* \alpha(x)| & \le r \; \underset{j \ge i}{\max} \{ |\phi_F^* \alpha(e_i)| \} \; \underset{j \ge i}{\max} \{ |x_j| \} \\
& \ll \underset{j \ge i}{\max} \{ |x_j| \},
\end{align*}
and the bound (\ref{rootupper}) now follows from (\ref{xibound1}) and (\ref{xibound2}).  To prove (\ref{rootlower}), let $q$ be the smallest number such that $m_q = 0$.  If $q \le i$ then $m_i = 0$ and there is nothing to prove.  Otherwise, the inequalities $\phi_F^* \alpha(e_j) \ge 0$ from Lemma \ref{phishape} and $x_j \ge 0$ for $j < q$ from (\ref{xibound2}) give

\bes
\phi_F^* \alpha(x) \ge \phi_F^* \alpha(e_i x_i) + \sum_{j=q}^r \phi_F^* \alpha( x_j e_j).
\ees
Applying (\ref{xibound1}) and (\ref{xibound2}) gives a constant $C_3 > 0$ such that

\bes
\phi_F^* \alpha(x) \ge e^{-r-1} | \phi_F^* \alpha(e_i)| t^{-1} e^{m_i} - C_3 t^{-1} e^{\delta m_1}.
\ees
If we assume that $m_i \ge \delta m_1 + C_1$ for $C_1$ satisfying $e^{-r-1} | \phi_F^* \alpha(e_i)| - C_3 e^{-C_1} > 0$, then we have $\phi_F^* \alpha(x) \gg t^{-1} e^{m_i}$ as required.

\end{proof}

The second input we shall need is a bound on the pointwise norm of $k_t$.

\begin{lem}
\label{kpoint}

We have

\bes
k_t( \exp(H) ) \ll t^{n-r} \prod_{\alpha \in \tD^+} ( 1 + t |\alpha(H)| )^{-1/2}
\ees
for $H \in \ga_0$ and $\lambda \in B^*$.

\end{lem}

\begin{proof}

The Harish-Chandra transform of $k_t$ is equal to

\begin{align*}
\widehat{k_t}(\nu) & = h_t^2(\nu) \\
& = \left( \sum_{w \in W} h( w \nu + t\lambda) \right)^2 \\
& = \sum_{w \in W} h( w \nu + t\lambda)^2 + \sum_{ w_1 \neq w_2} h( w_1 \nu + t\lambda) h( w_2 \nu + t\lambda) \\
& = \sum_{w \in W} h( w \nu + t\lambda)^2 + s(\nu, t\lambda),
\end{align*}
where $s(\nu, t\lambda)$ satisfies 

\be
\label{sschwartz}
\| (1 + |\nu|)^k s( \nu, t\lambda) \|_{L^1(\ga^*_0)} \ll_{k,A} t^{-A}
\ee
as a function of $\nu$ for all $\lambda$.  It suffices to bound $| k_t( \exp(H) ) |$ with $H \in B$.  Inverting the Harish-Chandra transform as in $\mathsection$\ref{sec212} gives

\begin{align*}
k_t(\exp(H)) & = \int_{W \backslash \ga^*_0} \widehat{k_t}(\nu) \varphi_{\nu}(\exp(H)) | c(\nu) |^2 d\nu \\
& = \int_{\ga^*_0} h( \nu + t\lambda)^2 \varphi_{\nu}(\exp(H)) | c(\nu) |^2 d\nu + \int_{W \backslash \ga^*_0} s(\nu, t\lambda) \varphi_{\nu}(\exp(H)) | c(\nu) |^2 d\nu.
\end{align*}
By combining the bounds $| \varphi_{\nu}(\exp(H)) | \le 1$, $| c(\nu) |^2 \ll 1 + |\nu|^{n-r}$, and (\ref{sschwartz}), we may estimate the second integral by

\bes
\int_{W \backslash \ga^*_0} s(\nu, t\lambda) \varphi_{\nu}(\exp(H)) | c(\nu) |^2 d\nu \ll_{h, A} t^{-A}.
\ees
Let $B_1^* \subset \ga^*_{0,+}$ be a precompact open set that contains $B^*$ and is bounded away from the walls.  We divide the domain of the first integral into $-tB_1^*$ and $\ga^*_0 \setminus - tB_1^*$.  We know that $h( \nu + t\lambda)$ is rapidly decaying in $\nu$ and $t$ for $\nu \notin -t B_1^*$, and together with $| \varphi_{\nu}(\exp(H)) | \le 1$ this gives

\bes
\int_{\ga^*_0 \setminus - tB_1^*} h( \nu + t\lambda)^2 \varphi_{\nu}(\exp(H)) |c(\nu)|^2 d \nu \ll_{A} t^{-A},
\ees
so that

\be
\label{lastint}
k_t( \exp(H) ) = \int_{-tB_1^*} h( \nu + t\lambda)^2 \varphi_{\nu}(\exp(H)) | c(\nu) |^2 d\nu + O_{A}(t^{-A}).
\ee
When $\nu \in -tB_1^*$, we apply Theorem \ref{phithm} to the set $-B_1^*$ and $B$ as chosen here to obtain

\bes
\varphi_{\nu}(\exp(H)) \ll_{B,B_1^*} \prod_{\alpha \in \tD^+} ( 1 + t |\alpha(H)| )^{-1/2}.
\ees
Combining this with (\ref{lastint}) gives

\bes
k_t( \exp(H) ) \ll \prod_{\alpha \in \tD^+} ( 1 + t |\alpha(H)| )^{-1/2} \int_{-tB_1^*} h( \nu + t\lambda)^2 |c(\nu)|^2 d \nu + O_{A}(t^{-A}),
\ees
and the bound $| c(\nu) |^2 \ll 1 + |\nu|^{n-r}$ completes the proof.

\end{proof}

We shall prove (\ref{Lpbd}) by interpolating between the cases $p = \infty$ and $p = 2$.  To begin with the case $p = \infty$, proving a bound for the $L^1 \rightarrow L^\infty$ norm of $T_{\bm,F}$ is the same as proving a bound for $\| T_{\bm,F}( \cdot, \cdot) \|_\infty$.  If we assume that $B$ is sufficiently small that there is at most one nonzero term in the sum

\bes
T_{\bm,F}(x,y) = \sum_{\gamma} K_{\bm,F}(x,\gamma y)
\ees
for all $x$ and $y$, then we have $\| T_{\bm,F}( \cdot, \cdot) \|_\infty = \| k_{\bm,F} \|_\infty$.  By Lemma \ref{kpoint}, we have

\begin{align*}
\| k_{\bm,F} \|_\infty & \le \sup \{ | k_t( \exp(H) )| | H \in \text{supp}( \beta_{\bm,F} ) \} \\
& \ll t^{n-r} \sup \left\{ \prod_{\alpha \in \tD^+} ( 1 + t |\alpha(H)| )^{-1/2} | H \in \text{supp}( \beta_{\bm,F} ) \right\}.
\end{align*}
If $\alpha \in \Delta_{i,F}$ and $m_i \ge \delta m_1 + C_1$, we may apply Lemma \ref{rootbound} to obtain $1 + t |\alpha(H)| \gg e^{m_i}$, while if $m_i < \delta m_1 + C_1$ we have the trivial bound $1 + t |\alpha(H)| \ge 1 \gg e^{m_i - \delta m_1}$.  Combining these, we obtain

\bes
\prod_{\alpha \in \tD^+} ( 1 + t |\alpha(H)| )^{-1/2} \ll \exp \left( -\tfrac{1}{4} \sum_{i=1}^r m_i | \tD_{i,F} | + \delta m_1 \eta(\bm)/4 \right),
\ees
where $\eta(\bm)$ is given by

\be
\label{eta}
\eta(\bm) = \sum_{m_i < \delta m_1 + C_1} | \tD_{i,F} | \le |\tD|.
\ee
The bound (\ref{Lpbd}) with $p = \infty$ follows if we choose $N \ge |\tD|/4$.\\

To prove the case $p = 2$, we first note that the $L^2 \rightarrow L^2$ norm of $T_{\bm,F}$ is equal to $\sup \{ | \widehat{k}_{\bm,F}(-\nu)| | \nu \in \cS \}$, where $\cS \subset \ga^* / W$ is the set of spectral parameters of the joint eigenfunctions in $L^2(X)$.  It is known that if $\nu \in \cS$, then either $\nu \in \ga^*_0$, or $\text{Re}(\nu)$ is singular and $\| \text{Im}(\nu) \| \le \| \rho \|$, see \cite[Thm 8.1, $\mathsection$8, Ch. IV]{He2} or \cite[Thm. 16.6 and $\mathsection$16.5 ex. 7]{K}.  Let $B_1^*$ and $B_2^*$ be compact sets such that $B^* \subset B_1^* \subset B_2^* \subset \ga^*_{0,+}$, each set contains an open neighbourhood of the one preceeding it, and $B_2^*$ is bounded away from the walls.  The following two results allow us to reduce to the case in which $\nu \in t B_2^*$.

\begin{lem}
\label{offspec}

If $\nu \in \ga^*_{0,+} \setminus B_2^*$ and $\kappa \in \ga_0^*$ satisfies $\| \kappa \| \le \| \rho \|$, then $\widehat{k}_{\bm,F}( -t\nu - i \kappa) \ll_{\delta} t^{-r}$.  The implied constant is uniform in $\bm$.

\end{lem}

\begin{proof}

We have

\begin{align*}
\widehat{k}_{\bm,F}( -t\nu - i \kappa) & = \int_S k_{\bm,F}(x) \varphi_{t\nu + i\kappa}(x) dx \\
& = \int_S \tb_{\bm,F}(x) k_t(x) \varphi_{t\nu + i\kappa}(x) dx.
\end{align*}
As in the proof of Lemma \ref{kpoint}, we may invert the Harish-Chandra transform of $k_t$ to obtain

\bes
\widehat{k}_{\bm,F}( -t\nu - i \kappa) = \int_{-t B_1^*} \int_S \tb_{\bm,F}(x) \varphi_{\mu}(x) \varphi_{t\nu + i\kappa}(x) dx h(\mu - t\lambda)^2 | c(\mu)|^2 d\mu + O_{A}(t^{-A}).
\ees
The lemma now follows from Proposition \ref{phiphi} below and $|c(\mu)|^2 \ll 1 + |\mu|^{n-r}$.

\end{proof}

\begin{prop}
\label{phiphi}

If $\mu \in -B_1^*$ and $\nu \in \ga^*_{0,+} \setminus B_2^*$, then

\bes
\int_S \tb_{\bm,F}(x) \varphi_{t\mu}(x) \varphi_{t\nu + i\kappa}(x) dx \ll_{\delta} t^{-n},
\ees
where the implied constant is uniform in $\bm$.

\end{prop}

\begin{proof}

Unfolding the integrals over $K$ used to define $\varphi_{t\mu}$ and $\varphi_{t\nu + i\kappa}$, we have

\begin{align*}
\int_S \tb_{\bm,F}(x) \varphi_{t\mu}(x) \varphi_{t\nu + i\kappa}(x) dx & = \int_K \int_S \tb_{\bm,F}(x) \exp( (it \mu + \rho)( A(x)) \\
& \qquad \qquad + (it\nu - \kappa + \rho)(A(kx))) dx dk \\
& = \int_K \int_S \tb_{\bm,F}(x) a(k,x) \exp( it( \nu(A(kx)) + \mu(A(x)))) dx dk,
\end{align*}
where $a(k,x) = \exp( \rho(A(x)) + (\rho - \kappa)(A(kx)))$.  There is a natural identification of $T^*S$ with the principal bundle $G \times_K \p^*_0$.  We let $\cC_\nu \subset T^*S$ be the set of points of the form $(G,v)$, where $v$ is conjugate to $\nu$ under $K$, and define $\cC_\mu$ similarly.  We know that the differentials of $\nu(A(kx))$ and $\mu(A(x))$ with respect to $x$ lie in $\cC_\nu$ and $\cC_\mu$ respectively, and our assumption that $\nu \in \ga_{0,+}^* \setminus B_2^*$ and $-\mu \in B_1^*$ were separated implies that $\| \nabla \nu(A(kx)) + \nabla \mu(A(x)) \| \ge \epsilon$ for some $\epsilon > 0$ depending on $B_1^*$ and $B_2^*$.

We shall apply integration by parts with respect to $\nabla \nu(A(kx)) + \nabla \mu(A(x))$.  All derivatives of $\nabla \nu(A(kx)) + \nabla \mu(A(x))$ and $a(k,x)$ are bounded above.  It follows from Lemma \ref{betabox} that the $K$-biinvariant function $\tb_{\bm,F}$ satisfies the analogous bound to (\ref{betadiff}), i.e. that for any linear differential operator $D$ of degree $d$ on $S$ with continuous coefficients, we have

\bes
D \tb_{\bm,F} \ll_D t^d e^{- \delta m_1 d}.
\ees
To calculate the bound obtained by integration by parts, we shall begin by estimating the volume of the support of $\tb_{\bm,F}$ on $S$.  If we define $V(H) = \prod_{\alpha \in \tD^+} | \alpha(H) |$, the Weyl integration formula gives

\bes
\text{Vol}( \text{supp}(\tb_{\bm,F}) ) \ll \int_{ \text{supp}(\beta_{\bm,F}) } V(H) dH.
\ees
On the support of $\beta_{\bm,F}$, Lemma \ref{rootbound} implies (as in the proof of our bound for $\| k_{\bm,F} \|_\infty$) that

\be
\label{Vbound}
V(H) \ll t^{-(n-r)} \exp \left( \tfrac{1}{2} \sum m_i | \tD_{i,F} | + \delta m_1 \eta(\bm)/2 \right),
\ee
with $\eta(\bm)$ as in (\ref{eta}).  It follows from our construction of $\beta_{\bm,F}$ that

\be
\label{volBbound}
\text{Vol}( \text{supp}(\beta_{\bm,F}) ) \ll t^{-r} \exp \left( \sum m_i + \delta m_1 q(\bm) \right),
\ee
where $q(\bm)$ is the number of zero entries in $\bm$, and combining these gives

\begin{align*}
\int_S \tb_{\bm,F}(x) a(k,x) e^{ it( \nu(A(x)) + \mu(A(x)))} dx & \ll \text{Vol}( \text{supp}(\tb_{\bm,F}) ) \\
& \ll t^{-n} \exp \bigg( \sum m_i (1 + | \tD_{i,F} |/2) \\
& \qquad \qquad + \delta m_1 \eta(\bm)/2 + \delta m_1 q(\bm) \bigg).
\end{align*}

Each partial integration produces a factor of $t^{-1}$, and a factor of $t e^{- \delta m_1}$ from differentiating $\tb_{\bm,F}$.  Performing this $A$ times therefore gives

\begin{multline*}
\int_S \tb_{\bm,F}(x) a(k,x) \exp( it( \nu(A(x)) + \mu(A(x))) dx \\
\ll_{A} t^{-n} \exp \left( \sum m_i (1 + | \tD_{i,F} |/2)  + \delta m_1( -A + \eta(\bm)/2 + q(\bm)) \right).
\end{multline*}
If we choose $A$ to be large enough, the exponential expression will be less than 1.  We therefore have

\bes
\int_S \tb_{\bm,F}(x) a(k,x) \exp( it( \nu(A(x)) + \mu(A(x))) dx \ll_{\delta} t^{-n}
\ees
as required.

\end{proof}

We now estimate $\widehat{k}_{\bm,F}(-\nu)$ for $\nu \in tB_2^*$.  The Weyl integration formula gives

\begin{align*}
\widehat{k}_{\bm,F}(-\nu) & = \int_S \tb_{\bm,F}(x) k_t(x) \varphi_{\nu}(x) dx \\
& \ll \int_{\ga_0} \beta_{\bm,F}(H) |k_t( \exp(H)) \varphi_{\nu}( \exp(H))| V(H) dH,
\end{align*}
where $V(H)$ is as above.  If we assume that $H \in \text{supp}(\beta_{\bm,F})$, then reasoning as in the proof of our bound for $\| k_{\bm,F} \|_\infty$ gives

\bes
k_t( \exp(H)) \varphi_{\nu}( \exp(H)) \ll t^{n-r} \exp \left( -\tfrac{1}{2} \sum m_i | \tD_{i,F} | + \delta m_1 \eta(\bm)/2 \right),
\ees
and combining this with (\ref{Vbound}) we have

\bes
\widehat{k}_{\bm,F}(-\nu) \ll e^{\delta m_1 \eta(\bm)} \int_{\ga_0} \beta_{\bm,F}(H) dH.
\ees
Equation (\ref{volBbound}) then gives

\bes
\widehat{k}_{\bm,F}(\nu) \ll t^{-r} \exp \left( \sum m_i + \delta m_1 q(\bm) + \delta m_1 \eta(\bm) \right).
\ees
If we choose $N \ge r + |\tD|$, this completes the proof of (\ref{Lpbd}) when $p = 2$, and of Proposition \ref{truncbound} with $C_1$ as in Lemma \ref{rootbound}.

\end{proof}

\subsection{Summation of $L^p$ bounds}
\label{Lpsummation}

We now sum the bound of Proposition \ref{truncbound} over $\bm$ and $F$ to obtain a bound for $T_t T^*_t$.  We begin with summation over $\bm$.  Define

\bes
T_F = \sum_{\bm \in \cM_{t, \delta}} T_{\bm,F}.
\ees
We have

\begin{align*}
\| T_F f \|_p & \le \sum_{ \bm \in \cM_{t, \delta} } \| T_{\bm,F} f \|_p \\
& \ll \| f \|_{p'} t^{n(1 - 2/p) - r} \sum_{ \bm \in \cM_{t, \delta} } \exp \left( L(F,p)(\bm) + \chi(\bm) N \delta m_1 \right).
\end{align*}
Dividing the sum into the terms with $\chi(\bm)$ equal to 0 and 1 gives

\be
\label{TFsum}
\| T_F f \|_p \ll \| f \|_{p'} t^{n(1 - 2/p) - r} \sum_{ \substack{ \bm \in \cM_{t, \delta} \\ \chi(\bm) = 0 } } e^{ L(F,p)(\bm)}
+ \| f \|_{p'} t^{n(1 - 2/p) - r} \sum_{ \substack{ \bm \in \cM_{t, \delta} \\ \chi(\bm) = 1 } } e^{ L(F,p)(\bm) + N \delta m_1}.
\ee
We define $\cM_t^0$ to be the subset of $\cM_t$ with $m_r = 0$, and let $\pi : \{ \cM_{t,\delta} | \chi(\bm) = 1 \} \rightarrow \cM^0_t$ be the projection obtained by setting $m_r = 0$.  If $\chi(\bm) = 1$, there are constants $N'$ and $D$ so that

\bes
( 1/(2p) - 1/4) m_r | \tD_{r,F} | + ( 2 - 2/p) m_r \le N' \delta m_1 + D.
\ees
Because the fibers of $\pi$ have at most $C_1 + 2$ elements, we may restrict the second sum in (\ref{TFsum}) from $\{ \cM_{t,\delta} | \chi(\bm) = 1 \}$ to $\cM_t^0$ (while increasing $N$ if necessary) and enlarge the first sum to $\cM_t$ to obtain

\be
\label{TFsum2}
\| T_F f \|_p \ll \| f \|_{p'} t^{n(1 - 2/p) - r} \sum_{ \bm \in \cM_t } e^{ L(F,p)(\bm)} + \| f \|_{p'} t^{n(1 - 2/p) - r} \sum_{ \bm \in \cM_t^0 } e^{ L(F,p)(\bm) + N \delta m_1}.
\ee
If $C_j \in F$ is a cell, we define the function $L(C_j, p)$ by

\begin{align*}
L(C_j,p) & = ( 1/(2p) - 1/4) \sum_{ i \le j} | \tD_{i,F} | + (2/p) j \\
& = ( 1/(2p) - 1/4) ( |\tD| - | \tD_{C_j} |) + (2/p) j,
\end{align*}
where $\Delta_{C_j}$ is the set of roots that vanish on $C_j$.  The function $L(C_j,p)$ is the value of $L(F,p)(x)$ at the vertex of $\overline{S}_0$ corresponding to $C_j$.  Note that $L(C_j,p)$ depends only on the Weyl orbit of $C_j$, and not $F$.  The first sum in (\ref{TFsum2}) is the generalised geometric progression obtained by summing $e^{L(F,p)(x)}$ over the integer points in the simplex $(\log t + 1) \overline{S}_0$, and the second sum is (up to the $N \delta m_1$ term) the sum over one its boundary faces.  The following proposition will allow us to estimate these sums.

\begin{prop}
\label{rootcount}

Define

\bes
M(p) = \max \{ L(C_r,p), L(C_0,p) \} = \max \{ (n+r)/p - (n-r)/2, 0 \}.
\ees
If $C \in \cC$ is a cone with $\dim C \notin \{ 0,r \}$, then $L(C,p) < M(p)$ for all $p \ge 2$.

\end{prop}

\begin{proof}

It suffices to prove the analogous statement for the linear functions $K(C,x) = (1/4 - x/2) | \tD_C| + 2x \dim C$ for $x \in [0, 1/2]$.  For $0 \le s \le r$, define

\bes
D(s) = \max \{ | \tD_C | | \dim C = s \} \quad \text{and} \quad K_0(s,x) = (1/4 - x/2) D(s) + 2s x.
\ees
We then have $K(C,x) \le K_0(\dim C,x)$.  The linear function $K_0(s,x)$ interpolates between the points $(0, D(s)/4)$ and $(1/2, s)$, and to show that this collection of functions is dominated by $K_0(0,x)$ and $K_0(r,x)$, it suffices to show that $D(s)$ is strictly concave up as a function of $s$.

The cones $C$ are in bijection with Levi subgroups $M$ of $G$ satisfying $A \subseteq M$, in such a way that $\exp(C)$ generates a maximal $\R$-split torus in the center of $M$ and $\Delta_C$ are the restricted roots of $M$.  Using our assumption that $\g_0$ was simple over $\R$, and Cartan's classification of globally symmetric spaces, it is then easy to check that $D(s)$ is concave up as required.

\end{proof}

Note that $M(p)$ is a piecewise linear function of $1/p$, with a kink point at $p = 2(n+r)/(n-r)$.  Proposition \ref{rootcount} implies that the function $L(F,p)(x)$ attains its maximum on $\overline{S}_0$ at either $(0, \ldots,0)$ (if $p > 2(n+r)/(n-r)$), $(1, \ldots,1)$ (if $p < 2(n+r)/(n-r)$), or on the edge joining them (if $p = 2(n+r)/(n-r)$).  It follows that the first sum in (\ref{TFsum2}) satisfies the estimate

\bes
\sum_{ \bm \in \cM_t } e^{L(F,p)(\bm)} \ll (\log t) t^{M(p)}
\ees
uniformly for $p \ge 2$, and that if $p \neq 2(n+r)/(n-r)$ we have

\bes
\sum_{ \bm \in \cM_t } e^{L(F,p)(\bm)} \ll_p t^{M(p)}.
\ees
To estimate the second sum, let $\partial \overline{S}_0$ be the boundary face of $\overline{S}_0$ on which $x_r = 0$.  If $\delta$ is chosen sufficiently small, Proposition \ref{rootcount} implies that there will be an $\epsilon > 0$ such that the linear functional $L(F,p)(x) + N \delta x_1$ attains its maximum on $\partial \overline{S}_0$ at $x = (0, \ldots, 0)$ for all $p > 2(n+r)/(n-r) - \epsilon$.  Moreover, for $p \le 2(n+r)/(n-r) - \epsilon$ and $\delta$ small we will have

\bes
\sup \{ L(F,p)(x) + N \delta x_1 | x \in \partial \overline{S}_0 \} < M(p).
\ees
Combining these gives

\bes
\sum_{ \bm \in \cM_t^0 } e^{L(F,p)(\bm) + N \delta m_1} \ll t^{M(p)}.
\ees
If we observe that $n(1 - 2/p) - r + M(p) = 2 \delta(p)$, then we have

\bes
\| T_F f \|_p \ll (\log t) t^{2 \delta(p)} \| f \|_{p'}
\ees
and

\bes
\| T_F f \|_p \ll_p t^{2 \delta(p)} \| f \|_{p'}
\ees
for $p \neq 2(n+r)/(n-r)$.  Theorem \ref{main} now follows by summing over $F$.

\section{Restrictions to maximal flat subspaces}
\label{flat}

Theorem \ref{restrict} may be proven using a slight modification of the method used to prove Theorem \ref{main}.  We shall assume that we are in the case of noncompact type.  The proof in the case of compact type is similar, and may be deduced from the results of $\mathsection$\ref{compact} and $\mathsection$\ref{compact2}.

We continue to use the notation of $\mathsection$\ref{Lp}, including the $K$-biinvariant kernel $k_t$, operator $T_t$, and the collection of flags $\cF$ and simplices $S(F)$.  We define $R$ to be the operator of restriction to $E$, and let $a \in C^\infty_0(E)$ be a real-valued cutoff function.  It suffices to bound the operator norms of

\begin{equation*}
a R T_t : L^2(X) \rightarrow L^p(E),
\end{equation*}
and if we let $\phi_1 \in C^\infty_0(E)$ and $\phi_2 \in C^\infty(X)$ be arbitrary functions with $\| \phi_1 \|_{p'} = \| \phi_2 \|_2 = 1$, it suffices to bound $\langle \phi_1, a R T_t \phi_2 \rangle$.  By taking adjoints and applying Cauchy-Schwarz, we have the inequality

\begin{equation*}
\langle \phi_1, a R T_t \phi_2 \rangle \le \langle a \phi_1, R T_t T_t^* R^* a \phi_1 \rangle.
\end{equation*}
Embed $E$ isometrically inside $\ga_0$, and let $P_t$ be the integral operator on $\ga_0$ with translation-invariant kernel

\begin{equation*}
P_t(H_1, H_2) = k_t( \exp(H_1 - H_2) ).
\end{equation*}
If we assume that the supports of $a$ and $k_t$ are small enough, we have $R T_t T_t^* R^* a \phi_1 = P_t a \phi_1$.  We therefore have

\begin{equation*}
\langle a \phi_1, R T_t T_t^* R^* a \phi_1 \rangle = \langle a \phi_1, P_t a \phi_1 \rangle,
\end{equation*}
and so it suffices to estimate the $L^{p'} \rightarrow L^p$ norms of $P_t$.\\

We do this by combining a dyadic decomposition of the kernel $P_t$ with an interpolation between $L^2 \rightarrow L^2$ and $L^1 \rightarrow L^\infty$ bounds as before.  The decomposition we make is simpler in this case, as we do not need to introduce the modified index set $\cM_{t,\delta}$.  If $g$ is as in $\mathsection$\ref{Lppartitions}, and ${\bf m} \in \cM_t$ and $1 \le i \le r$, we therefore define $f_{{\bf m}, i} \in C^\infty_0(\R)$ by

\bes
f_{{\bf m}, i}(x) = \bigg \{ \begin{array}{ll} g( t e^{-i} x_i), & m_i = 0, \\
\chi_{[0,\infty)}(x_i) [g( t e^{-m_i -i} x_i) - g( t e^{-m_i+1 - i} x_i)], & m_i > 0, \end{array}
\ees
and define $f_{\bf m} \in C^\infty_0(\R^r)$ by

\bes
f_{\bf m}(x) = \prod_{i=1}^r f_{{\bf m}, i}(x_i).
\ees
We let $\{ \beta_{\bm, F} | \bm \in \cM_t, F \in \cF \} \cup \beta_\infty$ be the partition of unity on $\ga_0$ derived from the functions $f_{\bf m}$ as in $\mathsection$\ref{Lppartitions}, define $k_{\bm, F} \in C^\infty_0( \ga_0)$ by

\bes
k_{\bm, F}(H) = k_t(\exp(H)) \beta_{\bm, F}(H),
\ees
and let $P_{\bm, F}$ be the operator with kernel

\bes
P_{\bm, F}(H_1, H_2) = k_{\bm, F}( H_1 - H_2)
\ees
so that

\bes
P_t = \sum_{ \substack{ \bm \in \cM_t \\ F \in \cF } } P_{\bm, F}.
\ees

The $L^1 \rightarrow L^\infty$ and $L^2 \rightarrow L^2$ norms of $P_{\bm, F}$ are bounded by

\begin{align*}
\| P_{\bm, F} f \|_\infty & \le \| k_{\bm, F} \|_\infty \| f \|_1 \\
\| P_{\bm, F} f \|_2 & \le \| k_{\bm, F} \|_1 \| f \|_2.
\end{align*}
If we define $J(F,p)(x)$ to be the linear functional

\bes
J(F,p)(x) = - 1/4 \sum_{i=1}^r x_i | \tD_{i,F} | + 2/p \sum_{i=1}^r x_i,
\ees
then we may prove the following bound for the $L^p \rightarrow L^{p'}$ norm of $P_{\bm, F}$ by bounding $\| k_{\bm, F} \|_1$ and $\| k_{\bm, F} \|_\infty$ using Theorem \ref{phithm} as in $\mathsection$\ref{Lpbounds}.

\begin{prop}

We have the bound $\| P_{\bm, F} f \|_p \ll t^{n - r - 2r/p} \exp( J(F,p)(\bm)) \| f \|_{p'}$.

\end{prop}

If $C_j \in F$ is a cell, we define

\bes
J(C_j, p) = - 1/4 \sum_{i \le j} | \tD_{i,F} | + 2j/p,
\ees
which is the value of $J(F,p)(x)$ at the vertex of $\overline{S}_0$ corresponding to $C_j$.  The conclusion of Theorem \ref{restrict} in the case $n > 3r$ may be deduced from the following lemma as in $\mathsection$\ref{Lpsummation}.

\begin{lem}

If $n > 3r$, we have $J(C_0,p) > J(C_j,p)$ for all $j > 0$ and $p \ge 2$.

\end{lem}

\begin{proof}

It is clear that $J(C_0, \infty) > J(C_j, \infty)$ for $j > 0$, so it suffices to show the same for $p = 2$.  Because

\bes
-|\tD|/4 + D(j)/4 + j \ge - 1/4 \sum_{i \le j} | \tD_{i,F} | + j = J(C_j,2),
\ees
where $D(j)$ is as in the proof of Proposition \ref{rootcount}, it suffices to show that

\be
\label{Jineq}
J(C_0,2) \ge -|\tD|/4 + D(j)/4 + j
\ee
with equality iff $j = 0$.  We know that $D(j)$ is concave up as a function of $j$, and so $D(j)/4 + j$ is also.  We know that equality holds in (\ref{Jineq}) when $j = 0$, and so it suffices to prove that strict inequality holds when $j = r$.  However, this is equivalent to our assumption that $n > 3r$.

\end{proof}

In cases (b) and (c), Theorem \ref{restrict} follows by examining the functions $J(C_j, p)$ for the finte number of globally symmetric spaces to which these cases apply.  The sharpness of the upper bounds in the case of compact type follows from the remarks of $\mathsection$\ref{compactsharp}.

\section{Bounds for spherical functions on noncompact groups}
\label{phibounds}

We shall prove Theorems \ref{phithm} and \ref{phiasympthm} by studying the expression

\be
\label{Kint1}
\varphi_{t\lambda}( \exp(H)) = \int_K \exp( (\rho + it\lambda)( A( k \exp(H)) ) ) dk
\ee
for $\varphi_\lambda$ as an oscillatory integral over $K$.  We define

\bes
\phi(k, H, \lambda) = \lambda( A( k \exp(H)) )
\ees
to be the phase of this integral, so that we may rewrite (\ref{Kint1}) as

\be
\label{Kint3}
\varphi_{t\lambda}( \exp(H)) = \int_K b(k,H) e^{it \phi(k,H,\lambda) } dk
\ee
where $b(k,H) = \exp( \rho( A( k \exp(H)) ) )$ is a function with all derivatives uniformly bounded.  We shall prove a uniformisation theorem for $\phi$ in $\mathsection$\ref{phiuniform}, which will reduce Theorems \ref{phithm} and \ref{phiasympthm} to an application of stationary phase to the integral (\ref{Kint3}) in $\mathsection$\ref{phiphase}.

\subsection{The critical set of $\phi$}
\label{phinotation}

We begin by recalling some properties of the critical point set of $\phi$, taken from \cite{DKV}.  Note that we shall always talk about the critical points of $\phi$ with respect to the variable $k$ only.  Let $\Delta_0^+ = \Delta^+ \cup \{ 0 \}$, and for every $\alpha \in \tD_0^+$, choose a vector $Y_\alpha \in (\g_{\alpha} + \g_{-\alpha}) \cap \gk_0$ so that $\{ Y_\alpha | \alpha \in \tD_0^+ \}$ is an orthonormal basis of $\gk_0$ with respect to $-\langle \; , \, \rangle$.  Note that when $\alpha = 0$, we are chosing a basis for $\gm_0$.  We also let $Y_\alpha$ denote the corresponding left-invariant vector fields on $K$.  Fix a point $l \in K$, and define $V \subseteq \ga_0$ to be the subspace

\begin{equation*}
V = \Ad_l^{-1} \ga_0 \cap \ga_0.
\end{equation*} 
We recall that a Levi subgroup of $G$ is called semi-standard if it contains $A$.

\begin{lem}
\label{Levidef}

There is a semi-standard Levi subgroup $A \subseteq L \subseteq G$ with real Lie algebra $\gl_0$ such that $V$ is the center $\ga_{0,L}$ of $\gl_0$.

\end{lem}

\begin{proof}

Pick $H \in V$ generic, in the sense that if $\alpha(H) = 0$ then $\alpha(V) = 0$ for $\alpha \in \Delta$.  Let $L$ be the connected centraliser of $H$ in $G$, which is a semi-standard Levi subgroup whose Lie algebra $\gl_0$ is the centraliser of $H$ in $\g_0$.  We have

\bes
\gl_0 = \bigoplus_{ \substack{ \alpha \in \Delta \\ \alpha(H) = 0 } } \g_{0,\alpha}, \quad \ga_{0,L} = \bigcap_{ \substack{ \alpha \in \Delta \\ \alpha(H) = 0 } } \ker \alpha,
\ees
so that $V \subseteq \ga_{0,L}$.  As $\gl_0$ is stable under $\theta$ we may decompose $\gl_0$ as $(\p_0 \cap \gl_0) + (\gk_0 \cap \gl_0) = \p_{0,L} + \gk_{0,L}$.  The group $K_L = K \cap L$ is maximal compact in $L$, as it is compact with Lie algebra $\gk_{0,L}$.  The subspaces $\ga_0$ and $\Ad_l^{-1} \ga_0 \subseteq Z_{\g_0}(H) = \gl_0$ are maximal abelian in $\p_{0,L}$, and so there exists $l_0 \in K_L$ such that $\Ad_{l_0} \Ad_l^{-1} \ga_0 = \ga_0$.  This implies that $\Ad_l^{-1} \ga_0 = \Ad_{l_0}^{-1} \ga_0$, so that $\ga_{0,L} \subseteq \Ad_l^{-1} \ga_0$ and $\ga_{0,L} \subseteq V$.  This completes the proof.

\end{proof}

\begin{defin}

It follows from the proof of Lemma \ref{Levidef} that $l \in M' K_L$, and we fix a decomposition $l = w l_0$ with $w \in M'$ and $l_0 \in K_L$ for the remainder of $\mathsection$\ref{phibounds}.

\end{defin}

  We define $X_\alpha = \Ad_l^{-1} Y_\alpha$ for $\alpha \in \tD_0^+$, and also let $X_\alpha$ denote the corresponding left-invariant vector field on $K$.  Decompose $\ga$ as an orthogonal direct sum $\ga = \ga_L + \ga^L$.  We let $\Delta_L$ be the set of roots that vanish on $\ga_L$, which is exactly the root system of $L$, and let $\Delta^L = \Delta \setminus \Delta_L$ be its complement.  We let $\Delta_L^+ = \Delta_L \cap \Delta^+$ and $\Delta^L_+ = \Delta^L \cap \Delta^+$.

\begin{prop}[Proposition 5.4 of \cite{DKV}]
\label{critical}

Fix $H \in \ga_0$ and regular $\lambda \in \ga^*_0$, and let $K_H$ be the stabiliser of $H$ in $K$.  The function $\phi(k, H, \lambda)$ is right-invariant under $K_H$, and its critical point set is equal to $M' K_H$.

\end{prop}

\begin{lem}
\label{al}

If $\lambda \in \ga_0^*$ is regular, $l$ is a critical point of $\phi(k, H, \lambda)$ if and only if $H \in \ga_{0,L}$.

\end{lem}

\begin{proof}

As $l \in M' K_L$, we clearly have $l \in M' K_H$ if $H \in \ga_{0,L}$.  For the converse, suppose that $H \in \ga_0$ is such that $l = w' k_H$ for $w' \in M'$ and $k_H \in K_H$.  We then have $\Ad_l^{-1} w' H = \Ad_{k_H}^{-1} H = H$, so that $H \in \Ad_l^{-1} \ga_0 \cap \ga_0 = \ga_{0,L}$.

\end{proof}

The following result is stated as Proposition 6.5 of \cite{DKV}, however we have included a derivation to avoid any possible error in converting the result to our notation.

\begin{prop}
\label{hessian}

When $H \in \ga_{0,L}$, the Hessian of $\phi$ with respect to the vector fields $\{ X_\alpha | \alpha \in \tD_0^+ \}$ at $l$ is diagonal, and satisfies

\begin{equation}
(D\phi)_{\alpha \alpha} = \tfrac{1}{2} \langle \lambda, \alpha \rangle ( 1 - e^{2\alpha(wH)} ).
\end{equation}

\end{prop}

\begin{proof}

We have

\bes
l \exp(t X_\alpha) \exp(H) = \exp( t Y_\alpha) \exp(wH) l,
\ees
and so if we consider the Iwasawa decomposition

\bes
\exp( t Y_\alpha) \exp(wH) = n(t) \exp(V(t)) k(t)
\ees
then we have $X_\alpha^2 \phi(l, H, \lambda) = \lambda(V''(0))$.  As $n(0) = k(0) = e$ and $V(0) = wH$, we may write

\bes
n(t) = \exp(tN_1 + t^2 N_2 + O(t^3)), \qquad k(t) = \exp(tK_1 + t^2 K_2 + O(t^3)),
\ees
so that

\begin{align}
\notag
\exp( t Y_\alpha) \exp(wH) & = \exp(tN_1 + t^2 N_2) \exp( wH + t^2 V''(0)/2) \exp(tK_1 + t^2 K_2) + O(t^3) \\
\label{hessianIwasawa}
\exp( t Y_\alpha) & = \exp(tN_1 + t^2 N_2) \exp( t^2 V''(0)/2) \\
\notag
& \qquad \qquad \exp(t \Ad_{\exp(wH)}K_1 + t^2 \Ad_{\exp(wH)}K_2) + O(t^3).
\end{align}
Equating first order terms gives

\bes
Y_\alpha = N_1 + \Ad_{\exp(wH)}K_1,
\ees
and if we write $Y_\alpha = V_\alpha + V_{-\alpha}$ with $V_{\pm \alpha} \in \g_{0,\pm \alpha}$ we may solve this to obtain

\bes
N_1 = (1 - e^{2\alpha(wH)}) V_\alpha, \qquad K_1 = e^{\alpha(wH)} Y_\alpha.
\ees
Applying the Baker-Campbell-Hausdorff formula in (\ref{hessianIwasawa}) and equating second order terms gives

\bes
0 = N_2 + V''(0)/2 + \Ad_{\exp(wH)}K_2 + [N_1, \Ad_{\exp(wH)}K_1]/2.
\ees
Because $N_2 + \Ad_{\exp(wH)}K_2 \in \ga^\perp$, this implies that

\begin{align*}
V''(0) & = -\text{proj}_\ga [N_1, \Ad_{\exp(wH)}K_1] \\
& = (e^{2\alpha(wH)} - 1) \langle V_\alpha, V_{-\alpha} \rangle H_\alpha.
\end{align*}
Our assumption that $\langle Y_\alpha, Y_\alpha \rangle = -1$ implies that $\langle V_\alpha, V_{-\alpha} \rangle = -1/2$, so that

\bes
X_\alpha^2 \phi(l, H, \lambda) = \lambda(V''(0)) = \tfrac{1}{2} \langle \lambda, \alpha \rangle (1 - e^{2\alpha(wH)})
\ees
as required.  The proof that the off-diagonal terms vanish is similar, and omitted.

\end{proof}

\subsection{Notation for complexifiying $\phi$}
\label{complexification}

When uniformising $\phi$, we will use different methods in the cases $l \in M'$ and $l \notin M'$.  Both cases involve analytically continuing $\phi$ into a complex domain, but the second case also involves blowing up the $\ga$-coordinate of this domain along the edges of a flag.  We treat the case $l \notin M'$ first, as it is the more difficult of the two.  We establish the notation used for doing this here.  By passing to an isogenous group if necessary, we may assume that $G$ is an analytic subgroup of a complex Lie group $G_\C$ with real Lie algebra $\g$, and that there is a closed complex subgroup $K_\C \subset G_\C$ with real Lie algebra $\gk$ such that $K = K_\C \cap G$.

\subsubsection{Generalities on complex germs}

We shall use the language of local complex spaces and holomorphic germs, for which we refer to \cite{GLS} for definitions.  All local complex spaces we shall work with will be regular, and we shall denote them by $(M,p)$, where $M$ is a complex manifold and $p \in M$ is a point.  We denote the ring of holomorphic germs on $(M,p)$ by $\cO(M,p)$.  All the local complex spaces we work with will have a natural complex conjugation, which we will denote by $c$ in all cases.  If $f$ is a holomorphic function we let $Z_f$ denote its zero divisor.

\subsubsection{Blowing up $\ga$}
\label{sec422}

Let $F = \{ C_0, \ldots, C_r \}$ be a flag as in $\mathsection$\ref{Lppartitions}.  Choose points $p_i \in C_i$ for each $i$, and let $J$ be the non-negative linear span of the $p_i$.  Let $V_i \subseteq \ga$ be the complex subspace spanned by $C_i$.  Let $\{ x_i | 0 \le i \le r-1 \}$ be the unique linear functions on $\ga$ such that $x_i(p_j) = 0$ if $i \ge j$ and $1$ otherwise, which form a co-ordinate system.  We define $\cA$ to be $\C^r$ with the standard linear co-ordinates $\{ z_i | 0 \le i \le r-1 \}$, and define $\pi_\cA : \cA \rightarrow \ga$ to be the blow-down map given by

\bes
\pi_\cA^* x_j = \prod_{i \le j} z_i.
\ees
The space $\cA$ is then a Zariski-open subset of the blowup of $\ga$ along the subspaces $V_0, \ldots, V_{r-2}$.  If we denote the interior of $J$ by $J^0$ it may be seen that

\bes
\pi_\cA^{-1}(J^0) = \{ (0,\infty) \times (0,1)^{r-1} \subset \R^r \subset \cA \},
\ees
and we define

\bes
\cJ = \overline{\pi_\cA^{-1}(J^0)} = \{ [0,\infty) \times [0,1]^{r-1} \subset \R^r \subset \cA \}.
\ees

\subsubsection{Blowing up $K_\C \times \ga \times \ga^*$}

Define the complex manifolds

\bes
S = \cA \times \ga^*, \quad X = K_\C \times S.
\ees
We shall denote points in $S$ by $s' = (u', \lambda')$.  We shall think of all roots $\alpha \in \Delta$ as holomorphic functions on $\cA$ by pullback, and let $\alpha_X$ denote the pullback of $\alpha$ to a function on $X$ under the natural projection.  We let $\pi_S : X \rightarrow S$ and $\pi_\lambda : X \rightarrow \ga^*$ be the natural projections.  We let $X_\alpha^+$ and $X_\alpha^-$ be the unique holomorphic and antiholomorphic vector fields on $X$ such that $X_\alpha^+ + X_\alpha^- = X_\alpha$ on the real submanifold $K \times \R^r \times \ga_0^*$ of $X$, and likewise for $Y_\alpha$.

\subsubsection{Germs of $\phi$}

We define $\cA_L = \pi_\cA^{-1}(\ga_L)$, and $S_L = \cA_L \times \ga^*$.  Let $p$ be the largest integer such that $V_p \subseteq \ga_L$.  We see that $Q = z_0 \ldots z_{p}$ is a defining function for $\cA_L$.  We let $Q_X$ denote the pullback of $Q$ to $X$.

We choose a point $s = (u, \lambda) \in (\cJ \cap \cA_L) \times \ga_0^*$ with $\lambda$ regular.  We let $x = (l, s) \in X$.  As $\phi$ is an analytic function on $K \times \ga_0 \times \ga_0^*$, we may complexify it and pull it back to obtain a germ in $\cO(X,x)$.  It follows from Lemma \ref{al} that $(l, u', \lambda') \in (X,x)$ is a critical point of $\phi$ exactly when $u' \in \cA_L$. 

\subsubsection{Divisors}

We let $D_i = \{ z \in \cA | z_{i} = 0 \}$ for $0 \le i \le r-1$ be the co-ordinate divisors on $\cA$, and define $S_i = D_i \times \ga^* \text{ for } 0 \le i \le p$.  We have $\pi_\cA(D_{i}) \subset V_i$, and $\cA_L = \bigcup_{0 \le i \le p} D_i$.  Let $q$ be the largest integer with $q \le p$ and $u \in D_{q}$.  We shall think of the divisors $D_i$ as subspaces of $(\cA,u)$ from now on, so that $D_i$ is empty if $u \notin D_i$ (and in particular if $i > q$).

\begin{lem}
\label{rootzero}

Recall that $\Delta_{j,F} = \{ \alpha \in \Delta | \; \alpha|_{V_j} \neq 0, \alpha|_{V_{j-1}} = 0 \}$ for $1 \le j \le r$.  If $\alpha \in \Delta_{j,F}$, we have $(Z_\alpha,u) = \sum_{0 \le i \le j-1} D_i$.

\end{lem}

\begin{proof}

Assume without loss of generality that $\alpha$ is nonnegative on $J$.  If $\alpha$ vanishes on $C_{j-1}$ but not $C_j$, then $\alpha/x_{j-1}$ must satisfy $C > \alpha/x_{j-1} > c > 0$ on $J^0$.  After pulling back to $\cA$, we see that the function $(z_0 \ldots z_{j-1})^{-1} \alpha(z)$ satisfies $C > (z_0 \ldots z_{j-1})^{-1} \alpha(z) > c > 0$ on $J^0$, and so it extends to an invertible function in $\cO(\cA,u)$.  The result now follows.

\end{proof}

\subsection{Uniformisation of $\phi$}
\label{phiuniform}

The uniformisation theorem for $\phi$ that we shall use is as follows.  We define $\Sigma_i = \{ \alpha \in \Delta^+ | \; w^{-1} \alpha|_{V_i} \neq 0 \}$.  Let $d = | \tE_q|$ and $d' = \dim K - d$, and identify $\C^d$ with $\C^{\tE_q}$ so that $\{ z_\alpha | \alpha \in \tE_q \}$ form a system of co-ordinates on this space.

\begin{thm}
\label{uniform1}

There is an isomorphism $f$

\bes
\xymatrix{ (X, x) \ar[rr]^f \ar[dr]_{\pi_S}&&  (\C^d \times \C^{d'}, 0) \times (S, s)  \ar[dl]^{0 \times \textup{id}}\\& (S, s)},
\ees
a function $\phi_S \in \cO( S, s)$, and a non-constant affine-linear map $L : \C^{d'} \longrightarrow \C$, such that $f$, $\phi_S$ and $L$ all commute with $c$, and such that

\bes
f_* \phi(z, z', s') = \phi_S(s') - \sum_{\alpha \in \tE_q} \langle \lambda', \alpha \rangle w^{-1}\alpha(u') z_\alpha^2 + Q(u')L(z').
\ees

\end{thm}

In other words, this expresses $\phi$ as the sum of a quadratic form on $\C^d$, and a linear function on $\C^{d'}$ that is zero exactly when $z' \in \cA_L$.  Proposition \ref{uniform2} and Corollary \ref{uniformcor} below carry out the uniformisation in the first set of co-ordinates $\C^d$.  They work by constructing a smooth subspace $(Y,x) \subset (X,x)$ that projects regularly to $(S,s)$ (see Definition 1.112 of \cite{GLS}), so that the fibers $Y_{s'}$ are smooth, and performing a change of variables that fixes $Y_{s'}$ and converts $\phi$ to a quadratic form transversally to $Y_{s'}$.  Proposition \ref{uniform2} builds $(Y,x)$ by induction on its codimension, and Corollary \ref{uniformcor} summarises the end result.

The main idea of the induction is as follows.  Let $\alpha \in \tE_q$.  The derivative $X^+_\alpha \phi$ vanishes on $Z_{w^{-1} \alpha}$, and so we may divide to obtain the holomorphic function $(w^{-1} \alpha)^{-1} X^+_\alpha \phi$.  The divisor of this function gives us our first submanifold $(Y,x)$, and we may repeat this process to decrease its dimension.

\begin{rem}

The argument we use does not require complexification.  We have written it in this way because we originally thought it was necessary in order to apply complex stationary phase in $\mathsection$\ref{sec61}, and because we felt the constructions were more familiar in a complex setting.

\end{rem}

For $0 \le i \le p$, let $K_i$ be the centraliser of $V_i$ in $K$.  Because $V_i \subseteq \ga_L$ for $0 \le i \le p$, we have $K_L \subseteq K_i$.  We have

\begin{align}
\notag
\text{Lie}(K_i) & = \text{span} \{ Y_\alpha | \alpha \in \tD^+_0, \; \alpha|_{V_i} = 0  \} \\
\label{lieKj}
& = \text{span} \{ X_\alpha | \alpha \in \tD^+_0 \setminus \tE_i \},
\end{align}
where the second equality follows from the fact that $l_0 \in K_L \subseteq K_i$.

\begin{prop}
\label{uniform2}

Let $\tE_{j-1} \subseteq R \subset \tE_j$ with $1 \le j \le q$ be given, and suppose that there exists a subspace $(Y,x) \subset (X,x)$ and an isomorphism

\bes
\xymatrix{ (X,x) \ar[rr]^f \ar[dr]_{\pi_S}&& (Y,x) \times (\C^R, 0) \ar[dl]^{\pi_S \times 0}\\& (S,s)}
\ees
with the following properties:

\begin{enumerate}[(a)]

\item
\label{hypd}

$f|_Y$ is the identity.

\item
\label{hypa}
$(Y,x)$ is invariant under $c$, and $f$ commutes with $c$.

\item
\label{hypb}
The projection $(Y,x) \rightarrow (S,s)$ is regular (see Definition 1.112 of \cite{GLS}).

\item
\label{hypf}
We have $Y_{s'} \subseteq l K_{i,\C}$ when $s' \in S_i$ with $0 \le i < j$, and $l K_{i,\C} \subseteq Y_{s'}$ when $s' \in S_i$ with $j \le i \le q$, where $Y_{s'}$ is the fiber of $Y$ above $s' \in S$.

\item
\label{hypc}
When $s' \in S_i$ with $j \le i \le q$, we have $l \in Y_{s'}$ and

\bes
T_l^{(1,0)} Y_{s'} = \textup{span} \{ X_\alpha^+ | \alpha \in \tD_0^+ \setminus R \}.
\ees

\item
\label{hype}

We have

\be
\label{transverse}
f_* \phi(y,z) = \phi(y) - \sum_{\alpha \in R} \langle \pi_\lambda(y), \alpha \rangle w^{-1}\alpha_X(y) z_\alpha^2.
\ee

\end{enumerate}

Then if $\beta \in \tE_j \setminus R$, there exists a subspace $(Y',x)$ and an isomorphism $f'$ having the same properties with respect to $R \cup \{ \beta \}$.

\end{prop}

\begin{proof}

We first note that property (\ref{hypb}) and the regularity of $(S,s)$ imply that both $(Y,x)$ and $(Y_s, l)$ are regular.  Let $R' = R \cup \{ \beta \}$.  Define $\phi_1 = f_* \phi$, and push the vector fields $X_\alpha^\pm$ forward under $f$ to obtain fields on $Y \times \C^R$, which we also denote $X_\alpha^\pm$.  Let $V_\alpha^\pm$ be the vector fields on $Y$ obtained by applying the natural projection $TY \times T \C^R \rightarrow TY$ to $X_\alpha^\pm$.  Hypothesis (\ref{hypc}) implies that when $s' \in S_i$ with $j \le i \le q$, and $\alpha \in \tD_0^+ \setminus R$, we have

\be
\label{VX}
V_\alpha^+ = X_\alpha^+ \in T_{(l,s')} Y,
\ee
and (\ref{transverse}) implies that when $y \in Y$ and $\alpha \in \tD_0^+$ we have

\be
\label{VX2}
V_\alpha^+ \phi_1(y,0) = X_\alpha^+ \phi(y).
\ee

\begin{lem}
\label{Xvanish}

If $y \in Y$ and $w^{-1}\beta_X(y) = 0$, we have $X_\beta^+ \phi(y) = 0$.

\end{lem}

\begin{proof}

Our assumption that $\beta \in \tE_j \setminus \tE_{j-1}$ implies that $w^{-1}\beta \in \tD_{j,F}$.  If $w^{-1}\beta_X(y) = 0$, Lemma \ref{rootzero} implies that $\pi_S(y) \in \bigcup_{0 \le i \le j-1} S_i$.  It follows that the image of $y$ under projection to $\cA$ and then blow-down by $\pi_\cA$ lies in $V_{j-1}$.  Proposition \ref{critical} implies that $\phi$ is right-invariant under $K_{j-1}$, and the lemma follows from (\ref{lieKj}).

\end{proof}

Lemma \ref{Xvanish} and (\ref{VX2}) imply that we also have $V_\beta^+ \phi_1(y,0) = 0$ when $w^{-1}\beta_X(y) = 0$.  We can therefore define an analytic function $\psi \in \cO(Y, x)$ by $\psi = (w^{-1}\beta_X)^{-1} V_\beta^+ \phi_1$, and define $Y'$ to be the zero locus of $\psi$.

We now establish (\ref{hypf}) for $Y'$.  The first inclusion $Y'_{s'} \subseteq l K_{i, \C}$ for $s' \in S_i$ and $0 \le i < j$ follows from $Y' \subset Y$.  To establish the second inclusion, let $j \le i \le q$ and assume that $(S_i,s)$ is nonempty.  Proposition \ref{critical} and the inclusion $\pi_\cA(D_i) \subseteq V_{i}$ imply that $l K_{i, \C}$ lies in the critical locus of $\phi$ when $s' \in S_i$, so that $V_\beta^+ \phi_1$ vanishes on $l K_{i, \C} \times S_i$.  It follows that $l K_{i, \C} \subseteq Y'_{s'}$ when $s' \in S_i$ and $w^{-1}\beta(u') \neq 0$, and because $j \le i$ and $S_i$ is irreducible, $w^{-1}\beta(u')$ is nonzero on an open dense subset of $S_i$.  The result then follows by continuity.  In particular, $x \in Y'$ and so $(Y',x)$ is a subspace of $(Y,x)$.  The following lemma implies that $(Y',x)$ and $(Y'_s, l)$ are both regular, and that

\bes
T_l^{(1,0)} Y'_{s'} = \textup{span} \{ X_\alpha^+ | \alpha \in \tD_0^+ \setminus R' \}
\ees
for $s' \in S_i$ with $j \le i \le q$ so that $Y'$ satisfies (\ref{hypc}).  Moreover, we see that $Y'$ satisfies (\ref{hypb}) by combining the regularity of the fiber $Y'_s$ with Proposition 1.85 and Theorem 1.115 of \cite{GLS}.

\begin{lem}
\label{Vq}

We have $V_\beta^+ \psi(l, s) \neq 0$, and $V_\alpha^+ \psi(l, s') = 0$ for all $\alpha \in \tD_0^+ \setminus R'$ and all $s' \in S_i$ with $j \le i \le q$.

\end{lem}

\begin{proof}

Let $\alpha \in \tD_0^+ \setminus R$, choose $j \le i \le q$ with $(S_i,s)$ nonempty, and assume that $s' \in S_i$ with $w^{-1}\beta(u') \neq 0$.  Equation (\ref{VX2}) implies that

\begin{align*}
V_\alpha^+ \psi(l, s') & = (w^{-1}\beta(u'))^{-1} V_\alpha^+ V_\beta^+ \phi_1(l, s') \\
& = (w^{-1}\beta(u'))^{-1} V_\alpha^+ X_\beta^+ \phi_1(l,s'),
\end{align*}
and (\ref{VX}) then gives

\bes
V_\alpha^+ \psi(l, s') = (w^{-1}\beta(u'))^{-1} X_\alpha^+ X_\beta^+ \phi(l,s').
\ees
We may apply Proposition \ref{hessian} to obtain

\bes
V_\beta^+ \psi(l, s') = -\langle \lambda', \beta \rangle e^{w^{-1}\beta(u')} \frac{ \sinh(w^{-1}\beta(u'))}{ w^{-1}\beta(u')},
\ees
and $V_\alpha^+ \psi(l, s') = 0$ for $\alpha \neq \beta$, and the result follows by continuity and that fact that $w^{-1}\beta$ is nonzero on an open dense subset of $S_i$.

\end{proof}

Lemma \ref{Vq} implies that the vector field $V_\beta^+$ is transverse to $Y'$.  Integrating along the flow of $V_\beta^+$ gives the following.

\begin{lem}

There is a unique isomorphism

\bes
\xymatrix{ (Y,x) \ar[rr]^g \ar[dr]_{\pi_S}&& (Y' ,x) \times (\C, 0) \ar[dl]^{\pi_S \times 0}\\& (S,s)}
\ees
with the properties that $g|_{Y'}$ is the identity, and if we let $(y', z)$ be the co-ordinates on $Y' \times \C$ then $g_* V_\beta^+ = \partial / \partial z$.

\end{lem}

We let $\phi_2 = g_* \phi_1$, and define $\phi_2'(y',z) = \phi_2(y',z) - \phi_2(y',0)$.  We know that $\phi_2'(y',z)$ vanishes to second order along $Y' \times 0$ by the definition of $Y'$, and we have $\partial \phi_2' / \partial z = V_\beta^+ \phi_1 = X_\beta^+ \phi$ so that $\phi_2'(y',z)$ vanishes identically when $w^{-1}\beta_X(y') = 0$ by Lemma \ref{Xvanish}.  We can therefore define $\psi_0(y',z) = (w^{-1}\beta_X(y'))^{-1} \phi_2'(y',z) \in \cO(Y' \times \C, (x,0))$, which also vanishes to second order on $Y' \times 0$ by continuity.

\begin{lem}

We have $\partial^2 \psi_0 / \partial z^2(x,0) \neq 0$.

\end{lem}

\begin{proof}

As in Lemma \ref{Vq}, we may calculate $\partial^2 \psi_0 / \partial z^2$ at $(l,s',0) \in Y' \times 0$ for $s' \in S_q$.  When $w^{-1}\beta(u') \neq 0$ we have

\begin{align*}
\partial^2 \psi_0 / \partial z^2 \big|_{(l,s',0)} & = (w^{-1}\beta(u'))^{-1} \partial^2 \phi_2' / \partial z^2 \big|_{(l,s',0)} \\
& = (w^{-1}\beta(u'))^{-1} (V_\beta^+)^2 \phi_1 \big|_{(l,s',0)}.
\end{align*}
By combining (\ref{VX}) and (\ref{VX2}) as before we may rewrite this as

\begin{align}
\notag
\partial^2 \psi_0 / \partial z^2 \big|_{(l,s',0)} & = (w^{-1}\beta(u'))^{-1} (X_\beta^+)^2 \phi \big|_{(l,s')} \\
\label{ybetareal}
& = -\langle \lambda', \beta \rangle e^{w^{-1}\beta(u')} \frac{ \sinh(w^{-1}\beta(u'))}{ w^{-1}\beta(u')},
\end{align}
and continuity gives the result.

\end{proof}

We may therefore define $z_\beta = \sqrt{ -\psi_0 / \langle \pi_\lambda(y'), \beta \rangle}$, which is an element of $\cO(Y' \times \C, (x,0))$ that satisfies $z_\beta(y',0) = 0$ and $\partial z_\beta / \partial z (y',0) \neq 0$ for all $y' \in (Y',x)$.  We define $g'$ to be the automorphism

\bes
(Y',x) \times (\C,0) \longrightarrow (Y',x) \times (\C,0), \quad (y',z) \mapsto (y', z_\beta),
\ees
and let $\phi_3 = g'_* \phi_2$.  The definitions of $\psi_0$ and $z_\beta$ imply that

\begin{align}
\notag
\phi_3(y', z_\beta) & = \phi_2(y',z) \\
\notag
& = \phi_2(y',0) + w^{-1} \beta_X(y') \psi_0(y',z) \\
\label{ybeta}
& = \phi_3(y',0) - \langle \pi_\lambda(y'), \beta \rangle w^{-1}\beta_X(y') z_\beta^2.
\end{align}
We define $f'$ to be the composition $f \circ g \circ g'$.  Equations (\ref{transverse}) and (\ref{ybeta}) imply that $\phi$ satisfies (\ref{hype}) with respect to $f'$ and $R'$, and $f'$ clearly acts as the identity on $Y'$.

It remains to establish (\ref{hypa}).  We first show that the function $\psi$ commutes with $c$.  

\begin{lem}

We have $\overline{\psi(y)} = \psi(\overline{y})$.

\end{lem}

\begin{proof}

We have $c_*(X_\beta^+) = X_\beta^-$ and $c_*(V_\beta^+) = V_\beta^-$, which implies that

\begin{align*}
\overline{ V_\beta^+ \phi(y)} & = V_\beta^- \overline{\phi}(y) \\
& = V_\beta^-( \phi \circ c)(y) \\
& = (c_*( V_\beta^-) \phi)(\overline{y}) \\
& = V_\beta^+ \phi(\overline{y}),
\end{align*}
and this implies the lemma.

\end{proof}

It follows that $Y'$ is invariant under $c$.  It can be shown that the conditions that define the map $g$ are also satisfied by $\overline{g(\overline{y})}$, and so by uniqueness we must have $g(\overline{y}) = \overline{g(y)}$.  This implies that $\phi_2$ and $\psi_0$ commute with $c$, and because $\langle \pi_\lambda(y'), \beta \rangle$ does also we have $z_\beta \circ c = \pm \overline{z_\beta}$.  Equation (\ref{ybetareal}) implies that $\partial z_\beta / \partial z (x)$ is real and nonzero, which means that we in fact have $z_\beta \circ c = \overline{z_\beta}$.  This completes the proof.

\end{proof}

Applying Proposition \ref{uniform2} inductively, we obtain

\begin{cor}
\label{uniformcor}

There exists a subspace $(Y, x) \subset (X, x)$ and an isomorphism

\be
\label{localiso1}
\xymatrix{ (X, x) \ar[rr]^f \ar[dr]_{\pi_S}&& (Y,x) \times (\C^d, 0) \ar[dl]^{\pi_S \times 0}\\& (S, s)}
\ee
with the following properties:

\begin{enumerate}[(a)]

\item
\label{cord}

$f|_Y$ is the identity.

\item
\label{cora}
$(Y, x)$ is invariant under $c$, and $f$ commutes with $c$.

\item
\label{corb}
The projection $(Y,x) \rightarrow (S, s)$ is regular.

\item
\label{corf}
We have $Y_{s'} = l K_{q,\C}$ when $s' \in S_q$, and $Y_{s'} \subseteq l K_{i,\C}$ when $s' \in S_i$ with $0 \le i < q$.

\item
\label{core}

We have

\be
\label{transverse2}
f_* \phi(y,z) = \phi(y) - \sum_{\alpha \in \tE_q} \langle \pi_\lambda(y), \alpha \rangle w^{-1}\alpha_X(y) z_\alpha^2.
\ee

\end{enumerate}

\end{cor}

\begin{proof}

We only need to describe how to change indices from $j$ to $j+1$ in the induction argument.  The only thing that requires explanation is how to pass from the inclusion $l K_{j,\C} \subseteq Y_{s'}$ when $s' \in S_{j}$ and $R \subset \tE_{j}$ to $Y_{s'} \subseteq l K_{j,\C}$ when $s' \in S_{j}$ and $\tE_{j} \subseteq R$.  In the boundary case when $R' = \tE_{j}$, the subspace $Y'$ produced by Proposition \ref{uniform2} satisfies $l K_{j,\C} \subseteq Y'_{s'}$ when $s' \in S_j$, and because $\dim Y'_{s'} = \dim K - | \tE_{j}| = \dim K_{j,\C}$ and both spaces are smooth we must in fact have $l K_{j,\C} = Y'_{s'}$ when $s' \in S_j$.

\end{proof}

\begin{proof}[Proof of Theorem \ref{uniform1}]

Because $Y$ is regular over $S$, there is a commutative diagram

\bes
\xymatrix{ (Y, x) \ar[rr]^i \ar[dr]_{\pi_S}&& ( Y_s, l) \times (S, s) \ar[dl]^{0 \times \text{id}}\\& (S, s)}
\ees
with $i$ an isomorphism.  It may be seen that we can choose $i$ to commute with $c$, for instance by choosing the generators $t_i$ of $\cO(Y_s, l)$ and their lifts $h_i$ to be real in the proof of Theorem 1.115 of \cite{GLS}.  Moreover, by condition (\ref{corf}) of Corollary \ref{uniformcor} and the fact that $(Y_s,l)$ is smooth, we may choose $i$ to satisfy

\be
\label{icondition}
i(l,s') = (l,s') \quad \text{for} \quad s' \in S_q.
\ee
We define the function $\xi: (S, s) \longrightarrow (Y, x)$ by $\xi(s') = i^{-1}( l \times s' )$, and define $\phi_0(y) = \phi(y) - \phi(\xi \circ \pi_S(y)) \in \cO(Y,x)$.  Proposition \ref{critical} implies that $\phi$ is right-invariant under $K_{i}$ when $s' \in S_i$, and it follows from this and condition (\ref{corf}) of Corollary \ref{uniformcor} that $\phi_0(y)$ vanishes when $s' \in S_L$.  We may therefore define $\psi = Q_X^{-1} \phi_0 \in \cO(Y, x)$, so that $\phi(y) = \phi( \xi \circ \pi_S(y)) + Q_X(y) \psi(y)$.  Transfer the fields $\{Y_\alpha^+ | \alpha \in \tD_L^+ \}$ to $Y \times \C^d$ via the map $f$ of Corollary \ref{uniformcor}, and let $W_\alpha^+$ be the projections of $Y_\alpha^+$ to $TY$ along $Y$. We wish to show that $\nabla_{Y_s} \psi |_l \neq 0$, where $\nabla_{Y_s}$ denotes the gradient along $Y_s$, and this will follow from knowing that $W_\alpha^+ \psi(x) \neq 0$ for some $\alpha \in \tD_L^+$.  We begin with the following lemma.

\begin{lem}
\label{Qdiff}

There exist $H_L \in \ga_{0,L}$ and $H^L \in \ga_0 \setminus \ga_{0,L}$ such that for all $\alpha \in \tD_L^+$, we have

\be
\label{Wpsi}
W_\alpha^+ \psi(x) = (\partial / \partial t) Y_\alpha^+ \phi(l, H_L + tH^L, \lambda) |_{t=0}.
\ee

\end{lem}

\begin{rem}

Vectors $H_L$ and $H^L$ in $\ga_0$ do not have to be orthogonal for the rest of $\mathsection$\ref{phiuniform}.

\end{rem}

\begin{proof}

Let $\alpha \in \tD_L^+$.  We have

\begin{align*}
W_\alpha^+ (Q_X \psi) & = W_\alpha^+ \phi \in \cO(Y,x)\\
\psi W_\alpha^+ Q_X + Q_X W_\alpha^+ \psi & = W_\alpha^+ \phi,
\end{align*}
and because $W_\alpha^+ Q_X = 0$ this gives $Q_X W_\alpha^+ \psi = W_\alpha^+ \phi$.  Let $\overline{u} = (\overline{u}_0, \ldots \overline{u}_{r-1}) \in D_q$ be a generic point near $u$, so that $\overline{u}_i = 0$ iff $i = q$, and let $s'(t) = (\overline{u} + te_{q}, \lambda)$.  Substituting $\xi(s'(t))$ into $Q_X W_\alpha^+ \psi = W_\alpha^+ \phi$ gives

\bes
t \prod_{ \substack{ i \le p \\ i \neq q} } \overline{u}_i W_\alpha^+ \psi(\xi(s'(t))) = W_\alpha^+ \phi(\xi(s'(t))),
\ees
and because $W_\alpha^+ \phi(y) = Y_\alpha^+ \phi(y)$ for $y \in Y$ by equation (\ref{transverse2}) we may rewrite this as

\bes
t \prod_{ \substack{ i \le p \\ i \neq q} } \overline{u}_i W_\alpha^+ \psi(\xi(s'(t))) = Y_\alpha^+ \phi(\xi(s'(t))).
\ees
Taking $\partial / \partial t$ of both sides and setting $t = 0$, and noting that $\xi(s'(0)) = (l, s'(0))$ by condition (\ref{icondition}) and our assumption that $\overline{u} \in D_q$, we obtain

\begin{align*}
\prod_{ \substack{ i \le p \\ i \neq q} } \overline{u}_i W_\alpha^+ \psi(l, s'(0)) = (\partial / \partial t) Y_\alpha^+ \phi(\xi(s'(t))) \big|_{t = 0}.
\end{align*}
We have $(\partial / \partial t) \xi( s'(t)) \big|_{t = 0} = \partial s'/ \partial t \big|_{t=0} + V \in T_{(l,s'(0))} X$ for some $V \in T_l K_\C$.  Because $\alpha \in \tD_L^+$, Proposition \ref{critical} implies that $Y^+_\alpha \phi(y,s')$ vanishes for all $y$ when $s' \in S_L$, and because $s'(0) \in S_L$, we have $V Y_\alpha^+ \phi(l, s'(0)) = 0$.  We therefore have

\begin{align*}
\prod_{ \substack{ i \le p \\ i \neq q} } \overline{u}_i W_\alpha^+ \psi(l, s'(0)) & = (\partial s'/ \partial t \big|_{t=0} + V) Y_\alpha^+ \phi(l, s'(0)) \\
& = (\partial / \partial t) Y_\alpha^+ \phi(l, s'(t)) |_{t=0}. \\
\end{align*}
We may rewrite this and let $\overline{u} \rightarrow u$ to obtain

\bes
W_\alpha^+ \psi(x) =  \underset{\overline{u} \rightarrow u}{\lim} \; \underset{t \rightarrow 0}{\lim} \; Q^{-1}(\overline{u} + te_{q})Y_\alpha^+ \phi(l, \overline{u} + te_{q}, \lambda).
\ees
As the function $Y^+_\alpha \phi(l,u',\lambda) \in \cO(\cA,u)$ vanishes on $\cA_L$, we know that $Q^{-1}(u') Y^+_\alpha(l,u',\lambda)$ extends to a function in $\cO(\cA,u)$ so that we may rewrite the limit more simply as

\bes
W_\alpha^+ \psi(x) = \underset{u' \rightarrow u}{\lim} Q^{-1}(u')Y_\alpha^+ \phi(l, u', \lambda).
\ees

Let $\cA^1$ be $\C^r$ with the standard linear co-ordinate functions $w_0, \ldots, w_{r-1}$, and define the maps $\cA \overset{\pi_1}{\longrightarrow} \cA^1 \overset{\pi_2}{\longrightarrow} \ga$ by

\begin{align*}
\pi_1^* w_j & = \Bigg\{ \begin{array}{ll} \prod_{i \le j} z_i & j \le p \\ \prod_{p < i \le j} z_i & p < j \end{array} \\
\pi_2^* x_j & = \Bigg\{ \begin{array}{ll} w_j & j \le p \\ w_{p} w_j & p < j. \end{array}
\end{align*}
We then have $\pi_\cA = \pi_2 \circ \pi_1$, and $\pi_1^* w_p = Q$.  We may naturally think of $\cA^1$ as a Zariski-open set in the blow-up of $\ga$ along $\ga_L$. The function $w_{p}^{-1}(w') Y_\alpha^+ \phi(l, w', \lambda)$ extends to a holomorphic germ in $\cO(\cA^1, \pi_1(u))$, and we have

\bes
Q^{-1}(u')Y_\alpha^+ \phi(l, u', \lambda) = \pi_1^* ( w_{p}^{-1}(w') Y_\alpha^+ \phi(l, w', \lambda)) \in \cO(\cA, u).
\ees
Write $u = (u_0, \ldots, u_{r-1})$.  Define $H_L = \pi_\cA(u) \in \ga_{0,L}$, and $H^L \in \ga_0$ by

\bes
x_j(H^L) = \Bigg\{ \begin{array}{ll} 0 & j \le p \\
            1 & j = p \\
	    \prod_{p < i \le j} u_j & p < j.
           \end{array}
\ees
We then have $H^L \in \ga_0 \setminus \ga_{0,L}$ and $\pi_2( \pi_1(u) + t e_{p}) = H_L + tH^L$.  As $u_{q} = 0$ we have 

\bes
w_{p}( \pi_1(u) + te_{p}) = w_{p}( te_{p}) = t,
\ees
so that

\begin{align*}
W_\alpha^+ \psi(x) & = \underset{u' \rightarrow u}{\lim} Q^{-1}(u')Y_\alpha^+ \phi(l, u', \lambda) \\
& = \underset{w' \rightarrow \pi_1(u)}{\lim} w_{p}^{-1}(w') Y_\alpha^+ \phi(l, w', \lambda) \\
& = \underset{t \rightarrow 0}{\lim} \; t^{-1} Y_\alpha^+ \phi(l, \pi_1(u) + t e_{p}, \lambda) \\
& = \underset{t \rightarrow 0}{\lim} \; t^{-1} Y_\alpha^+ \phi(l, H_L + tH^L, \lambda) \\
& = (\partial / \partial t) Y_\alpha^+ \phi(l, H_L + tH^L, \lambda) |_{t=0}
\end{align*}
as required.

\end{proof}

We now have to prove that the RHS of (\ref{Wpsi}) is nonzero for some choice of $\alpha \in \tD_L^+$.  We begin by simplifying the expression as follows.  If $\alpha \in \tD^+_L$, we write $Y_\alpha = V_\alpha + V_{-\alpha}$ with $V_{\pm \alpha} \in \g_{\pm \alpha}$.

\begin{lem}
\label{Wpsi1}

Let $\alpha \in \tD_L^+$, $H_L \in \ga_{0,L}$, and $H^L \in \ga_0 \setminus \ga_{0,L}$.  We have

\bes
(\partial / \partial t) Y_\alpha \phi( l, H_L + tH^L, \lambda) \big|_{t = 0} = \alpha(H^L) \langle \Ad^{-1}_{l} H_\lambda, V_{-\alpha} - V_\alpha \rangle,
\ees
where $H_\lambda \in \ga_0$ is dual to $\lambda$ under the Killing form.

\end{lem}

\begin{proof}

For $t \in \R$, we have the Iwasawa decomposition $l \exp( H_L + tH^L) = n(t) a(t) k(t)$.  If we write the first-order approximation to the Iwasawa decomposition of $l \exp( s Y_\alpha) \exp( H_L + tH^L)$ as

\begin{multline*}
l \exp( s Y_\alpha) \exp( H_L + tH^L) = n(t) \exp( s N_1 + O(s^2)) \\
a(t) \exp( s A_1 + O(s^2) ) k(t) \exp( s K_1 + O(s^2) ),
\end{multline*}
then we have $Y_\alpha \phi( l, H_L + tH^L, \lambda) = \lambda(A_1)$.  Moving the terms involving $s$ to the right and equating first-order parts gives

\begin{align*}
e^{ -t \alpha(H^L)} V_\alpha + e^{t \alpha(H^L)} V_{-\alpha} & = \Ad_{a(t) k(t)}^{-1} N_1 + \Ad_{k(t)}^{-1} A_1 + K_1 \\
\Ad_{k(t)} ( e^{ -t \alpha(H^L)} V_\alpha + e^{t \alpha(H^L)} V_{-\alpha} ) & = \Ad_{a(t)}^{-1} N_1 + A_1 + \Ad_{k(t)} K_1.
\end{align*}
We have $\Ad_{a(t)}^{-1} N_1 + \Ad_{k(t)} K_1 \in \ga^\perp$, and so

\begin{align*}
\lambda(A_1) & = \langle H_\lambda, \Ad_{k(t)} ( e^{ -t \alpha(H^L)} V_\alpha + e^{t \alpha(H^L)} V_{-\alpha} ) \rangle \\
Y_\alpha \phi( l, H_L + tH^L, \lambda) & = \langle \Ad^{-1}_{k(t)} H_\lambda, e^{ -t \alpha(H^L)} V_\alpha + e^{t \alpha(H^L)} V_{-\alpha} \rangle,
\end{align*}
Differentiating at $t = 0$ gives

\begin{align*}
(\partial / \partial t) Y_\alpha \phi( l, H_L + tH^L, \lambda) \big|_{t = 0} = \partial / \partial t \langle \Ad^{-1}_{k(t)} H_\lambda, Y_\alpha \rangle \big|_{t=0} + \alpha(H^L) \langle \Ad^{-1}_{k(0)} H_\lambda, V_{-\alpha} - V_\alpha \rangle.
\end{align*}
Because $\Ad^{-1}_{k(t)} H_\lambda \in \p$ and $Y_\alpha \in \gk$, the first term vanishes.  Because $l \exp(H_L) = \exp(w^{-1} H_L) l$, we have $k(0) = l$, which completes the proof.

\end{proof}

\begin{lem}
\label{Wpsi2}

There is $\alpha \in \tD_L^+$ such that $\alpha(H^L) \langle \Ad^{-1}_l H_\lambda, V_{-\alpha} - V_\alpha \rangle \neq 0$.

\end{lem}

\begin{proof}

The Lie algebra of $L$ is given by $\gl = \bigoplus_{\alpha \in \Delta_L} \g_\alpha$, and our choice of $Y_\alpha = V_\alpha + V_{-\alpha}$ implies that the Cartan $-1$-eigenspace $\p_L \subset \gl$ is given by

\be
\label{Lcartan}
\p_L = \text{span} \{ V_\alpha - V_{-\alpha} | \alpha \in \tD^+_L \}.
\ee
Suppose that $\alpha(H^L) \langle \Ad^{-1}_l H_\lambda, V_{-\alpha} - V_\alpha \rangle = 0$ for all $\alpha \in \tD^+_L$.  Because $\Ad_l^{-1} H_\lambda \in \p_L$, (\ref{Lcartan}) implies that

\bes
\Ad_l^{-1} H_\lambda \in \bigoplus_{ \substack{ \alpha \in \Delta_L \\ \alpha(H^L) = 0 } } \g_\alpha.
\ees
The RHS is the Lie algebra $\gl'$ of a semi-standard Levi subgroup $L' \subset L$, where the inclusion is proper because $H^L \notin \ga_{0,L}$.  We let $K_{L'} = K \cap L'$, which is maximal compact in $L'$.  Let $\ga_{0,L'}$ be the center of $\gl_0'$.  There is $l' \in K_{L'}$ such that $\Ad_{l'} \Ad_l^{-1} H_\lambda \in \ga_0$, and because $H_\lambda$ is regular this implies that $l \in M' l' \subset M' K_{L'}$.  It follows that $\Ad_l^{-1} \ga_0 \cap \ga_0$ contains $\ga_{0,L'}$, which contradicts Lemma \ref{Levidef}.

\end{proof}

Combining Lemmas \ref{Wpsi1} and \ref{Wpsi2} and applying the holomorphy of $\phi$, we see this also implies that there is $\alpha \in \tD^+_L$ such that

\bes
(\partial / \partial t) Y_\alpha^+ \phi(l, H_L + tH^L, \lambda) |_{t=0} \neq 0,
\ees
and by Lemma \ref{Qdiff} this gives $W_\alpha^+ \psi(x) \neq 0$ as required.

We have $\phi(y) = \phi( \xi \circ \pi_S(y)) + Q_X(y) \psi(y)$, so that $i_* \phi(y_s, s') = i_* \phi(l, s') + Q(u') i_* \psi(y_s, s')$.  Because $\nabla_{Y_s} \psi(l, s) \neq 0$ and $\psi$ commutes with $c$, there is an isomorphism

\bes
\xymatrix{ (Y_s, l) \times (S,s) \ar[rr]^{i'} \ar[dr]_{0 \times \text{id}}&& (\C^{d'}, 0) \times (S,s) \ar[dl]^{0 \times \text{id}}\\& (S, s)}
\ees
such that $(i' \circ i)_* \psi$ is a non-constant affine-linear function $L$, and such that $i'$ and $L$ both also commute with $c$.  Defining $\phi_S(s') = i_* \phi(l, s')$ completes the proof.

\end{proof}

\subsection{Proof of Theorem \ref{phithm}}
\label{phiphase}

We now use Theorem \ref{uniform1} to bound the contribution to the integral (\ref{Kint3}) from points away from $M'$.  Throughout $\mathsection$\ref{phiphase}, $H = H_L + H^L$ will denote the orthogonal decomposition of $H$ corresponding to the decomposition $\ga = \ga_L + \ga^L$ associated to a semi-standard Levi subgroup $L$.

\begin{prop}
\label{statphase}

Let $B$ and $B^*$ be as in the statement of Theorem \ref{phithm}.  Let $l \in K$ with $l \notin M'$.  Recall the notation of $\mathsection$\ref{phinotation} associated to $l$, including the Levi $L$ and decomposition $H = H_L + H^L$.  There is an open set $l \in U \subset K$ such that for all $b_0 \in C^\infty_0(U)$ and all $(H, \lambda) \in B \times B^*$, we have

\be
\label{Kint2}
\int_K b_0(k) b(k,H) e^{it \phi(k, H, \lambda)} dk \ll (1 + \| t H^L \|)^{-A} \prod_{\alpha \in \tD^+ } ( 1 + t | \alpha(H) | )^{-1/2}.
\ee
The implied constant depends on $A$, $B$, $B^*$, $l$, and $b_0$.

\end{prop}

\begin{proof}

Assume that the collection of cones $J(F)$ associated to $F \in \cF$ in $\mathsection$\ref{sec422} satisfies $\ga_0 = \cup_{F \in \cF} J(F)$.  Choose $F \in \cF$, and recall the notation associated to $F$ in $\mathsection$\ref{complexification}.  Define $\cB = \pi_\cA^{-1}(B) \cap \cJ$ and $\cB_L = \pi_\cA^{-1}(B) \cap \cJ \cap \cA_L$, so that $\pi_\cA(\cB) = B \cap J$ and $\pi_\cA(\cB_L) = B \cap J \cap \ga_L$.  For each $s' \in \cB_L \times B^*$, let $V_{s'} \subset K_\C$ and $W_{s'} \subset \cA \times \ga^*$ be open neighbourhoods of $l$ and $s'$ respectively such that $U_{s'} = V_{s'} \times W_{s'}$ realises the isomorphism $f$ of Theorem \ref{uniform1}.  We also assume that $W_{s'}$ intersects only the divisors $S_i$ that contain $s'$.  Let $V^0_{s'} \subset V_{s'}$ and $W^0_{s'} \subset W_{s'}$ be smaller open neighbourhoods such that $\overline{V^0_{s'}} \subset V_{s'}$ and $\overline{W^0_{s'}} \subset W_{s'}$.  By compactness, there exists a finite collection of points $\{ s_i \}$ such that $W^0_{s_i}$ cover $\cB_L \times B^*$.  We define $U = \cap_{s_i} V^0_{s_i} \cap K$, and let $U_\cB \subset \cB$ be a relatively open neighbourhood of $\cB_L$ in $\cB$ such that $U_\cB \times B^* \subset \cup_{s_i} W^0_{s_i}$.

Fix a point $s_i = (u_i, \lambda_i)$.  Let $q$ be the largest integer such that $u_i \in D_q$ as in $\mathsection$\ref{complexification}, and let $\Sigma_q$ be as in $\mathsection$\ref{phiuniform}.  Fix $s' = (u', \lambda') \in (U_\cB \times B^*) \cap W^0_{s_i}$, and let $H = \pi_\cA(u')$.  Applying Theorem \ref{uniform1} and restricting to the fibre above $s'$, we obtain an open set $U' \subset \R^d \times \R^{d'}$ and a real-analytic diffeomorphism $f: U \longrightarrow U'$ such that

\bes
f_* \phi(x, x', s') = \phi_S(s') - \sum_{\alpha \in \tE_q} \langle \lambda', \alpha \rangle w^{-1} \alpha(H) x_\alpha^2 + Q(u') L(x').
\ees
Making this change of co-ordinates in the integral (\ref{Kint2}) gives

\begin{multline*}
\int_K b_0(k) b(k,H) e^{it \phi(k, H, \lambda')} dk = e^{it \phi_S(s')} \int_{U'} c(x,x') \\
\exp \bigg( it \bigg[ -\sum_{\alpha \in \tE_q} \langle \lambda', \alpha \rangle w^{-1} \alpha(H) x_\alpha^2 + Q(u')L(x') \bigg] \bigg) dx dx',
\end{multline*}
where $c \in C^\infty_0(U')$ is the product of $b_0$, $b$, and the determinant of the Jacobian of $f$.  Because $f$ extends to a complex analytic function on the set $U_{s_i}$, which contains $\overline{U} \times \overline{W^0_{s_i}}$, all derivatives of $c$ with respect to $x$ and $x'$ are bounded for $(x,x') \in U'$, uniformly for $s' \in (U_\cB \times B^*) \cap W^0_{s_i}$.  Application of van der Corput and the bound $|Q(u')| \gg \| H^L \|$ therefore gives

\begin{align}
\label{linearised}
& \int_{U'} c(x,x') \exp \bigg( it \bigg[ -\sum_{\alpha \in \tE_q} \langle \lambda', \alpha \rangle w^{-1} \alpha(H) x_\alpha^2 + Q(u') L(x') \bigg] \bigg) dx dx' \\
\notag
& \qquad \qquad \ll_{A,B^*} (1 + t |Q(u')|)^{-A} \prod_{\alpha \in \tE_q } ( 1 + t | w^{-1} \alpha(H) | )^{-1/2} \\
\label{linearised2}
& \qquad \qquad \ll_{A,B^*} (1 + t |Q(u')|)^{-A} \prod_{ \substack{ \alpha \in \tD^+ \\ \alpha|_{V_q} \neq 0 } } ( 1 + t | \alpha(H) | )^{-1/2}
\end{align}
We now pass from (\ref{linearised2}) to the RHS of (\ref{Kint2}).  We first apply the following lemma.

\begin{lem}
\label{QH}

If $u' \in \cB$ and $H = \pi_\cA(u')$, we have $|Q(u')| \sim \| H^L \|$.

\end{lem}

\begin{proof}

We have $Q = \pi_\cA^* x_p$, where $x_p$ is as in $\mathsection$\ref{sec422}.  Because $\ga_{0,L} \cap J = C_p$ and $H \in J$, we also have $|x_p(H)| \sim \| H^L \|$.

\end{proof}

It remains to enlarge the product in (\ref{linearised2}) to one over $\tD^L_+$, and then $\tD^+$.  If $\alpha \in \tD^L_+$ satisfies $\alpha|_{V_q} = 0$, our assumption on $W_{s_i}$ intersecting only those divisors $S_j$ that contain $s_i$ implies that $\alpha / Q$ is holomorphic on $W_{s_i}$.  This implies that $|\alpha(u')| \ll |Q(u')| \ll \| H^L\|$, and so $(1 + t|\alpha(H)|)^{-1/2} \gg (1 + \| t H^L \|)^{-1/2}$.  We may therefore enlarge the product in (\ref{linearised}) to one over $\tD^L_+$, which is the same as the bound

\bes
(\ref{linearised}) \ll_{A,B^*} (1 + \| t H^L \|)^{-A} \prod_{\alpha \in \tD^L_+ } ( 1 + t | \alpha(H) | )^{-1/2}.
\ees
If $\alpha \in \tD^+_L$, we have $\alpha(H) \ll \| H^L \|$ and so we are free to enlarge the product further to $\tD^+$.  Applying this bound for each set $W^0_{s_i}$, we obtain the inequality (\ref{Kint2}) for all $s' \in U_\cB \times B^*$.

We may therefore assume that $s' \in (\cB \setminus U_\cB) \times B^*$, which is equivalent to assuming that $|Q(u')| > \delta > 0$, or that $\| H^L \| > \delta > 0$.  After possibly shrinking $U$, this implies that $\| \nabla_k \phi \| > \epsilon > 0$ on $U$, and so we have

\bes
\int_K b_0(k) b(k,H) e^{it \phi(k, H, \lambda')} dk \ll_{A,B,B^*} t^{-A}.
\ees
As $\pi_\cA(\cB) = J \cap B$, applying this argument for every $F \in \cF$ completes the proof.

\end{proof}

We now bound the contribution to the integral (\ref{Kint3}) from a neighbourhood of $M'$.  It will be convenient to reduce the integral to one over $R = M \backslash K$, which may be done as $\phi$ and $b_0$ are both left-invariant under $M$.  We shall use the uniformisation of $\phi$ at the points $W \in R$ given by Proposition \ref{uniform3} below, which may be proved in exactly the same way as Proposition \ref{uniform2}.  

Let $R_\C = M_\C \backslash K_\C$ be the complexification of $R$, and let $S = \ga \times \ga^*$ and $X = R_\C \times \ga \times \ga^*$.  Let $\pi_S : X \rightarrow S$ be the natural projection.  Choose $w \in W$, $H \in B$ and $\lambda \in B^*$, and let $s = (H, \lambda) \in S$ and $x = (w,s) \in X$.  We extend $\phi$ to a holomorphic germ in $\cO(X,x)$.

\begin{prop}
\label{uniform3}

There is an isomorphism $f$

\bes
\xymatrix{ (X, x) \ar[rr]^f \ar[dr]_{\pi_S}&&  (\C^{n-r}, 0) \times (S, s)  \ar[dl]^{0 \times \textup{id}}\\& (S, s)}
\ees
which commutes with $c$ and such that

\be
\label{weyluniform}
f_* \phi(z, H', \lambda') = \phi(w,H', \lambda') - \sum_{\alpha \in \tD^+} \langle \lambda', \alpha \rangle \alpha(wH') z_\alpha^2.
\ee

\end{prop}

As in the proof of Proposition \ref{statphase}, Proposition \ref{uniform3} implies that there is a neighbourhood $W \in U \subset R$ such that for all $b \in C^\infty_0(U)$ and all $(H, \lambda) \in B \times B^*$, we have

\bes
\int_R b(r) b_0(r,H) e^{it\phi(r,H,\lambda)} dr \ll_{B,B^*} \prod_{\alpha \in \tD^+} (1 + t |\alpha(H)))^{-1/2}.
\ees
Combined with Proposition \ref{statphase}, this completes the proof of Theorem \ref{phithm}.

\subsection{Proof of Theorems \ref{phiasympthm} and \ref{phiasympthm2}}

We prove Theorems \ref{phiasympthm} and \ref{phiasympthm2} by a more detailed analysis of the contribution to the integral

\be
\label{rint}
\varphi_{t\lambda}(\exp(H)) = \int_R b_0(r,H) e^{it\phi(r,H,\lambda)} dr
\ee
from the Weyl points.  We begin with four lemmas that provide a form of the stationary phase asymptotic adapted to the integrals we wish to study.

\begin{lem}
\label{statphase1}

Let $b \in C^\infty_0(\R \times \R_{\neq 0})$, and define $f \in C^\infty(\R_{ \neq 0} \times \R_{> 0} )$ by

\be
\label{fdef1}
\int b(x,y) e^{-it y x^2} dx = b(0,y) \pi^{1/2} |ty|^{-1/2} \exp( - i \pi \textup{sgn}(y) /4) + f(y,t).
\ee
If $g \in C^\infty_0(\R)$ and $k \ge 0$, we define

\bes
\| g \|_{C^k} = \sup \{ |g^{(j)}(x)| \; | \; x \in \R, 0 \le j \le k \}.
\ees
We then have

\be
\label{jdiffbd}
\frac{\partial^k f}{\partial y^k}(y,t) \ll |ty|^{-3/2} \sum_{i + j = k} |y|^{-i} \| (\partial^j b / \partial y ^j)(\cdot, y) \|_{C^{k+3}}
\ee
for all $k \ge 0$ and all $(y, t) \in \R_{ \neq 0} \times \R_{> 0}$, where the implied constant depend only on $k$ and a bound for the support of $b$ in the $x$-variable.

\end{lem}

\begin{proof}

We use induction.  Suppose that (\ref{jdiffbd}) is known for some $k \ge 0$ and all $b \in C^\infty_0(\R \times \R_{\neq 0})$.  Note that the base case $k = 0$ is given by the stationary phase asymptotic, see for instance \cite[Lemma 7.7.3]{Ho2}.  Differentiating the LHS of (\ref{fdef1}) and integrating by parts gives

\begin{align*}
\left( \frac{\partial}{\partial y} \right) \int b(x,y) e^{-it y x^2} dx & = \int (\partial b/ \partial y)(x,y) e^{-it y x^2} dx + \int b(x,y) (-it x^2) e^{-it y x^2} dx \\
& = \int ( (\partial b/ \partial y)(x,y) - (\partial / \partial x)( x b(x,y) )/2y ) e^{-it y x^2} dx.
\end{align*}
Comparing this with the derivative of the RHS of (\ref{fdef1}) gives

\begin{multline*}
\int ( (\partial b/ \partial y)(x,y) - (\partial / \partial x)( x b(x,y) )/2y ) e^{-it y x^2} dx \\
= ( (\partial b/ \partial y)(0,y) - b(0,y)/2y ) \pi^{1/2} |ty|^{-1/2} \exp( - i \pi \textup{sgn}(y) /4) + (\partial f / \partial y)(y,t).
\end{multline*}
Applying the inductive hypothesis to the two functions $\partial b/ \partial y$ and $(\partial / \partial x)( x b(x,y) )/2y$ separately gives the result.

\end{proof}

The following lemma may be proved in exactly the same way as Lemma \ref{statphase1}, with the base case provided by van der Corput's Lemma, see for instace \cite[Ch. VIII, $\mathsection$1.2, Corollary]{St}.

\begin{lem}
\label{statphase2}

Let $b \in C^\infty_0(\R \times \R_{\neq 0})$, and define $f \in C^\infty(\R_{ \neq 0} \times \R_{> 0} )$ by

\bes
\int b(x,y) e^{-it y x^2} dx = f(y,t).
\ees
We then have

\bes
\frac{\partial^k f}{\partial y^k}(y,t) \ll |ty|^{-1/2} \sum_{i + j = k} |y|^{-i} \| (\partial^j b / \partial y ^j)(\cdot, y) \|_{C^{k+1}}
\ees
for all $k \ge 0$ and all $(y, t) \in \R_{ \neq 0} \times \R_{> 0}$, where the implied constant depend only on $k$ and a bound for the support of $b$ in the $x$-variable.

\end{lem}

We recall the definition of $\ga_r$ and $\ga_r^*$ as the regular sets in $\ga_0$ and $\ga_0^*$, and $\| H \|_s$ as the distance from $H \in \ga_0$ to the singular set.

\begin{lem}
\label{errorterminduct}

Let $b \in C^\infty_0(\R^{n-r} \times \ga_0 \times \ga_r^*)$.  Let $R \subset \tD^+$, and let $d = |\tD^+ \setminus R|$.  We define

\bes
\sigma_w(R, \lambda, H) = - \sum_{\alpha \in R} \textup{sgn} ( \langle \lambda, \alpha \rangle \alpha(wH) ).
\ees
If $x \in \R^{n-r}$, we write $x = (x_R, x^R)$ with $x_R \in \R^{n - r - d} \simeq \R^R$ and $x^R \in \R^d \simeq \R^{\tD^+ \setminus R}$.  Suppose that there exists a function $f \in C^\infty(\R^d \times \ga_r \times \ga_r^* \times \R_{>0})$ with the following properties.

\begin{enumerate}[(a)]

\item
\label{errorterma}

If $H \in \ga_r$, $\lambda \in \ga_r^*$, and $t > 0$, we have

\begin{multline*}
\int b(x_R, x^R, H, \lambda) \exp \bigg( -it \sum_{\alpha \in \tD^+} \langle \lambda, \alpha \rangle \alpha(wH) x_\alpha^2 \bigg) dx_R = \exp \bigg( -it \sum_{\alpha \in \tD^+ \setminus R} \langle \lambda, \alpha \rangle \alpha(wH) x_\alpha^2 \bigg) \\
\left( \pi^{|R|/2} \prod_{\alpha \in R} (t | \langle \alpha, \lambda \rangle \alpha(wH)|)^{-1/2} \exp( i\pi \sigma_w(R, \lambda, H)/4) b(0, x^R, H, \lambda) + f(x^R, H, \lambda, t) \right).
\end{multline*}

\item
\label{errortermb}

The function $f$ satisfies

\begin{align*}
\left( \frac{\partial}{\partial x^R} \right)^p \left( \frac{\partial}{\partial H} \right)^q f & \ll \prod_{\alpha \in R} (t |\alpha(wH)|)^{-1/2} \frac{1}{t \| H \|_s^{q+1}},
\end{align*}
where the implied constant depends on $p$, $q$, and $b$.

\item
\label{errortermc}

There is a compact set $B \subset \R^d \times \ga_0 \times \ga^*_r$ such that $\textup{supp}(f) \subset B \times \R_{> 0}$.

\end{enumerate}

If $\beta \in \tD^+ \setminus R$ and we define $R' = R \cup \{ \beta \}$, there exists $f'$ satisfying the same conditions with respect to $R'$.

\end{lem}

\begin{proof}

We write $x = (x_{R'}, x^{R'})$ in the same way as $x = (x_{R}, x^{R})$.  If we apply property (\ref{errorterma}), we see that

\bes
\int b(x_{R'}, x^{R'}, H, \lambda) \exp \bigg( -it \sum_{\alpha \in \tD^+} \langle \lambda, \alpha \rangle \alpha(wH) x_\alpha^2 \bigg) dx_{R'}
\ees
is the sum of

\begin{multline}
\label{phasepart1}
\exp \bigg( -it \sum_{\alpha \in \tD^+ \setminus {R'}} \langle \lambda, \alpha \rangle \alpha(wH) x_\alpha^2 \bigg) \pi^{|R|/2} \prod_{\alpha \in R} (t | \langle \alpha, \lambda \rangle \alpha(wH)|)^{-1/2} \exp( i\pi \sigma_w(R, \lambda, H)/4) \\
\int b(0, x_\beta, x^{R'}, H, \lambda) \exp( -it \langle \lambda, \beta \rangle \beta(wH) x_\beta^2) dx_\beta
\end{multline}
and

\be
\label{phasepart2}
\exp \bigg( -it \sum_{\alpha \in \tD^+ \setminus {R'}} \langle \lambda, \alpha \rangle \alpha(wH) x_\alpha^2 \bigg) \int f(x_\beta, x^{R'}, H, \lambda,t) \exp( -it \langle \lambda, \beta \rangle \beta(wH) x_\beta^2) dx_\beta.
\ee

We deal with (\ref{phasepart1}) by writing $\ga_0$ as a direct sum of the kernel of $w^{-1} \beta$ and any transverse subspace, and applying Lemma \ref{statphase1} with $y = \langle \lambda, \beta \rangle \beta(wH)$.  Note that we are free to truncate the support of $b$ away from the set $\beta(wH) = 0$ to ensure that the hypotheses of Lemma \ref{statphase1} are satisfied.  This implies that

\begin{multline*}
\int b(0, x_\beta, x^{R'}, H, \lambda) \exp( -it \langle \lambda, \beta \rangle \beta(wH) x_\beta^2) dx_\beta \\
= b(0, x^{R'}, H, \lambda) ( t | \langle \lambda, \beta \rangle \beta(wH) | / \pi)^{-1/2} \exp( -i\pi \text{sgn}( \langle \lambda, \beta \rangle \beta(wH) ) / 4) + f_1'(x^{R'}, H, \lambda, t),
\end{multline*}
where $f_1'$ satisfies

\be
\label{f1diff}
\left( \frac{\partial}{\partial x^{R'}} \right)^p \left( \frac{\partial}{\partial H} \right)^q f_1' \ll (t |\beta(wH)|)^{-3/2} |\beta(wH)|^{-q} \ll (t |\beta(wH)|)^{-1/2} t^{-1} \| H \|_s^{-q-1}.
\ee
To deal with (\ref{phasepart2}), we define

\bes
f_2'(x^{R'}, H, \lambda, t) = \int f(x_\beta, x^{R'}, H, \lambda, t) \exp( -it \langle \lambda, \beta \rangle \beta(wH) x_\beta^2) dx_\beta.
\ees
We may show that $f_2'$ satisfies (\ref{errortermb}) with respect to $R'$ in the same way as we proved (\ref{f1diff}), by truncating $f$ away from the singular set in $H$, and applying Lemma \ref{statphase2} and the assumption that $f$ satisfied (\ref{errortermb}).  If we define

\begin{multline*}
f'(x^{R'}, H, \lambda, t) = \pi^{|R|/2} \prod_{\alpha \in R} (t | \langle \alpha, \lambda \rangle \alpha(wH)|)^{-1/2} \exp( i\pi \sigma_w(R, \lambda, H)/4) f_1'(x^{R'}, H, \lambda, t)\\
 + f_2'(x^{R'}, H, \lambda, t),
\end{multline*}
it may be seen that $f'$ satisfies the conditions of the lemma with respect to $R'$.

\end{proof}

As the conditions of Theorem \ref{errorterminduct} with $R = \emptyset$ are satisfied, we may proceed by induction to obtain

\begin{lem}
\label{errortermcor}

Let $b \in C^\infty_0(\R^{n-r} \times \ga_0 \times \ga_r^*)$.  There exists a function $f \in C^\infty(\ga_r \times \ga_r^* \times \R_{>0})$ with the following properties.

\begin{enumerate}[(a)]

\item
\label{errortermcora}

If $H \in \ga_r$, $\lambda \in \ga_r^*$, and $t > 0$, we have

\begin{multline}
\label{errortermunif}
\int b(x, H, \lambda) \exp \bigg( -it \sum_{\alpha \in \tD^+} \langle \lambda, \alpha \rangle \alpha(wH) x_\alpha^2 \bigg) dx = 
\pi^{(n-r)/2} \exp( i\pi \sigma_w(H, \lambda)/4) b(0, H, \lambda) \\
\times \prod_{\alpha \in \tD^+} (t | \langle \alpha, \lambda \rangle \alpha(H)|)^{-1/2}  + f(H, \lambda, t).
\end{multline}

\item
\label{errortermcorb}

The function $f$ satisfies

\be
\label{fdiff}
\left( \frac{\partial}{\partial H} \right)^p f \ll \prod_{\alpha \in \tD^+} (t |\alpha(H)|)^{-1/2} \frac{1}{t \| H \|_s^{p+1}},
\ee
where the implied constant depends on $p$ and $b$.
\end{enumerate}

\end{lem}

\begin{proof}[Proof of Theorem \ref{phiasympthm}]

Let $b \in C^\infty(R)$ be supported in a neighbourhood of $W$, and equal to 1 on a smaller neighbourhood of $W$.  We write

\bes
\varphi_{t\lambda}(\exp(H)) = \int_R b(r) b_0(r,H) e^{it\phi(r,H,\lambda)} dr + \int_R (1 - b(r)) b_0(r,H) e^{it\phi(r,H,\lambda)} dr.
\ees
Proposition \ref{statphase} implies that the second term may be absorbed into the error term in (\ref{phiasymp}).  By applying a partition of unity in the variables $(H, \lambda)$ and shrinking the support of $b$, we may use Proposition \ref{uniform3} to write the first term as a finite sum of integrals of the form (\ref{errortermunif}) multiplied by $\exp(i t \lambda(wH))$ for some $w \in W$.  Applying Lemma \ref{errortermcor} gives functions $\{ c_w \in C^\infty(\ga \times \ga_r^*) | w \in W \}$ and functions $\{ f_w \in C^\infty( \ga_r \times \ga_r^* \times \R_{> 0}) | w \in W \}$ satisfying (\ref{fdiff}) such that

\begin{multline*}
\int b(r) b_0(r,H) e^{it\phi(r,H,\lambda)} dr = \pi^{(n-r)/2} \prod_{\alpha \in \tD^+} (t | \langle \alpha, \lambda \rangle \alpha(H)|)^{-1/2} \\
\times \sum_{w \in W} \exp( it\lambda(wH) + i\pi \sigma_w(H, \lambda)/4) c_w(H, \lambda) + \sum_{w \in W} \exp( it\lambda(wH) )f_w(H, \lambda, t)
\end{multline*}
for $H \in B$ and $\lambda \in B^*$.  This gives an asymptotic for $\varphi_{t\lambda}$ of the same type as Theorem \ref{phiasympthm}, but with the presence of factors $c_w(H,\lambda) \in C^\infty(\ga \times \ga_r^*)$ in the main terms.  These may be calculated by comparison with the formula (9.10) of \cite{DKV} when $H \in B$ and $\lambda \in B^*$, which completes the proof.

\end{proof}

The proof of Theorem \ref{phiasympthm2} follows from a similar induction, with Lemma \ref{statphase2} used instead of Lemma \ref{statphase1}.

\section{Symmetric Spaces of Compact Type}
\label{compact}

We now consider the case in which $X$ is a locally symmetric space of compact type.  We assume without loss of generality that $X$ is a simply connected globally symmetric space $S = U/K$.  As in the noncompact case, most of the work in proving Theorem \ref{main} lies in establishing a sharp pointwise bound for the kernel of an approximate spectral projector, and the bound we shall use is exactly that of Theorem \ref{cpctphibd} for the spherical function $\varphi_\mu$ on $S$.  We shall prove this bound using the method of the previous sections, after first deriving an expression for $\varphi_\mu$ as an average of plane waves which is an analogue of the usual expression for $\varphi_\lambda$ as a $K$-integral in the noncompact case.

\subsection{Notation}
\label{compactnotn}

Let $(\gu_0, \theta)$ be a semisimple orthogonal symmetric Lie algebra of the compact type.  Let $\gu_0 = \gk_0 + i\p_0$ be the associated Cartan decomposition.  Let $(U,K)$ be the unique Riemannian symmetric pair associated to $(\gu_0, \theta)$ with $U$ simply connected and $K$ the connected subgroup with Lie algebra $\gk_0$.  Let $\g$ be the complexification of $\gu_0$.  Let $(\g_0, s)$ be the orthogonal symmetric Lie algebra dual to $(\gu_0, \theta)$, so that $\g_0 \subset \g$ is a real form of $\g$ and $s$ is the restriction of $\theta$ to $\g_0$, and the Cartan decomposition of $\g_0$ is $\gk_0 + \p_0$.  Let $G$ be a connected Lie group with real Lie algebra $\g_0$ and finite center.  After an isogeny, we may assume that $U$ and $G$ are both analytic subgroups of the simply connected group $G_\C$ with real Lie algebra $\g$.  We denote the Killing form on $\g$ by $\langle \; , \, \rangle$.  Let $\q_0$ be the orthogonal compliment of $\ga_0$ in $\p_0$ with respect to the Killing form, so that we have

\bes
\g_0 = \gk_0 + \p_0 = \gk_0 + \ga_0 + \q_0.
\ees
Let

\begin{align*}
G & = NAK, \quad g = n(g) \exp( A(g)) k(g) \\
\g & = \gk + \ga + \gn
\end{align*}
be an Iwasawa decomposition of $G$.  Let $M'$ and $M$ be the normaliser and centraliser of $\ga$ in $K$, and let $W$ be the Weyl group $M' / M$.  We let $M_0$ be the connected component of the identity in $M$, and $\gm$ be its Lie algebra.  Let $\Delta$ denote the set of roots of $\g$ with respect to $\ga$.  Note that we assume that $0 \in \Delta$ as in $\mathsection$\ref{sec211}.  We let $\Delta^+$ be a choice of positive roots, to which we associate the Lie algebra $\gn = \sum_{\alpha \in \Delta^+} \g_\alpha$ and positive Weyl chamber $\ga^+_0$.  We let $\ga_{0,+}^*$ denote the positive dual Weyl chamber.  We shall let $\tD$ denote the multiset on $\Delta$ in which each root is counted with multiplicity $m(\alpha)$, with $m(\alpha)$ as in $\mathsection$\ref{sec211}, and likewise for any subset of $\Delta$.  We let $\Delta_0^+ = \Delta^+ \cup \{ 0 \}$, and for every $\alpha \in \tD_0^+$ choose $Y_\alpha \in (\g_\alpha + \g_{-\alpha}) \cap \gk_0$ so that $\{ Y_\alpha | \alpha \in \tD_0^+ \}$ is an orthonormal basis of $\gk_0$ with respect to $-\langle \; , \, \rangle$.  Extend $\ga_0$ to a Cartan subalgebra $\gh_0$ of $\g_0$.  Define $T$ to be the connected subgroup of $U$ with Lie algebra $i \ga_0$, and let $\overline{T}$ be the image of $T$ in $S$ so that $\overline{T}$ is a maximal flat subspace of $S$.

\subsubsection{Spherical functions}

Define $\Lambda$ to be the set

\begin{equation*}
\Lambda = \left\{ \mu \in \ga^* : \frac{ \langle \mu, \alpha \rangle}{ \langle \alpha, \alpha \rangle} \in \Z^+ \text{ for } \alpha \in \Delta^+ \right\}.
\end{equation*}
For each $\mu \in \Lambda$, we extend $\mu$ to a linear functional on $\gh$ that is $0$ on $\gh \cap \gk$, and let $(\pi_\mu, V_\mu)$ denote the irreducible representation of $G_\C$ with highest weight $\mu$.  By \cite[Thm. 4.12, $\mathsection$4, Ch. II]{He}, the set of irreducible representations of $G_\C$ whose restriction to $G$ is spherical (that is, has a $K$-fixed vector) is exactly $\{ \pi_\mu | \mu \in \Lambda \}$.  Let $\langle \; , \, \rangle_\pi$ be a $\pi_\mu(U)$-invariant Hermitian inner product on $V_\mu$, and let $d(\mu)$ be the dimension of $V_\mu$.  We let $s^*$ be the automorphism of $\ga^*$ such that $\pi_\mu$ and $\pi_{s^* \mu}$ are contragredient, which is given by composing the map $\mu \mapsto -\mu$ with the long element of the Weyl group.

Let $e_\mu \in V_\mu$ belong to the weight $\mu$ and let $v_\mu \in V_\mu$ be a unit vector fixed under $K$.  We define the two functions $\varphi_\mu$ and $b_\mu$ in $C^\infty(G_\C)$ by

\begin{align*}
\varphi_\mu(g) & = \langle \pi_\mu(g^{-1}) v_\mu, v_\mu \rangle_\pi, \\
b_\mu(g) & = \langle \pi_\mu(g^{-1}) e_\mu, v_\mu \rangle_\pi.
\end{align*}
The restriction of $\varphi_\mu$ to $U$ is the spherical function with spectral parameter $\mu$, and we shall see in $\mathsection$\ref{compactbeam} that the restriction of $b_\mu$ to $U$ may be thought of as a higher rank Gaussian beam.

If $f \in C^\infty(U)$ is a $K$-biinvariant function, we define its spherical transform $\widehat{f}$ by

\bes
\widehat{f}(\mu) = \int_U f(u) \varphi_{s^* \mu}(u) du, \quad \mu \in \Lambda.
\ees
The following inversion formula for the spherical transform is a consequence of the Peter-Weyl theorem, see for instance \cite[Prop. 9.1, $\mathsection$9, Ch. III]{He}.

\begin{prop}
\label{sphericalinversion}

We have

\bes
f(u) = \sum_{\mu \in \Lambda} \varphi_\mu(u) \widehat{f}(\mu) d(\mu).
\ees

\end{prop}

\subsubsection{Complex Iwasawa coordinates}

The mapping

\begin{equation*}
(X, H, J) \rightarrow \exp X \exp H \exp J \qquad (X \in \gn, H \in \ga, J \in \gk)
\end{equation*}
is a holomorphic diffeomorphism of of a neighbourhood of $0$ in $\g$ onto a neighbourhood $U_\C^0$ of $e$ in $G_\C$.  We can therefore analytically continue the map $A : G \rightarrow \ga_0$ to a map $U^0_\C \rightarrow \ga$ by defining

\begin{equation*}
A : \exp X \exp H \exp J \rightarrow H.
\end{equation*}
As $b_\mu$ is holomorphic on $G_\C$, $N$-invariant on the left and $K$-invariant on the right, we have

\begin{equation*}
b_\mu( \exp X \exp H \exp J ) = b_\mu( \exp H) = e^{-\mu(H)} b_\mu(e).
\end{equation*}
It follows that $b_\mu(e) \neq 0$, and we shall always normalise $b_\mu$ by $b_\mu(e) = 1$ so that

\begin{equation}
\label{beamdef}
b_\mu(u) = e^{-\mu(A(u))}, \qquad u \in U_\C^0.
\end{equation}
We shall need the following invariance property of $A$.

\begin{lem}
\label{MA}

If $u \in U_\C^0$ and $m \in M_0$ satisfy $mu \in U_\C^0$, then $A(u) = A(mu)$.

\end{lem}

\begin{proof}

Because $e_\mu$ is fixed by $M_0$, we have $b_\mu(mu) = b_\mu(u)$ for all $\mu$ and $m \in M_0$.  If we define the lattice $\ga_\Lambda = \{ H \in i \ga_0 | \mu(H) \in 2\pi i \Z \text{ for all } \mu \in \Lambda \}$, this implies that $A(mu) - A(u) \in \ga_\Lambda$.  Shrinking $U_\C^0$ if necessary, the lemma follows.

\end{proof}

The following lemma allows us to extend the representation (\ref{beamdef}) to $T U_\C^0$.

\begin{lem}
\label{extendA}

We may extend $A$ to a function $A : T U_\C^0 \rightarrow \ga / \ga_\Lambda$.

\end{lem}

\begin{proof}

If $g = \exp(H) u \in T U_\C^0$ with $H \in i \ga_0$ and $u \in U_\C^0$, we define $A(\exp(H)u) = H + A(u) \in \ga / \ga_\Lambda$.  To show that this is well defined, assume that $g = \exp(H_1) u_1 = \exp(H_2) u_2$.  We have

\bes
b_\mu(g) = e^{-\mu(H_1 + A(u_1))} = e^{-\mu(H_2 + A(u_2))},
\ees
so that

\bes
\mu(H_1 + A(u_1)) - \mu(H_2 + A(u_2)) \in 2 \pi i \Z.
\ees
As this holds for all $\mu \in \Lambda$, we have $H_1 + A(u_1) - H_2 + A(u_2) \in \ga_\Lambda$ as required. 

\end{proof}

\subsection{Approach to Proving Theorem \ref{cpctphibd}}

Let $B^* \subset \ga_{0,+}^*$ be a compact set that is bounded away from the singular set, and assume that $\mu \in B^*$ and $t > 0$ satisfy $t \mu \in \Lambda$.  We have

\begin{equation}
\label{beamint}
\varphi_{t\mu}(\exp(H)) = \int_K b_{t\mu}(k\exp(H)) dk, \quad H \in i\ga_0.
\end{equation}
We see from equation (\ref{beamdef}) that there is a strong formal analogy between this expression for $\varphi_{t\mu}$ and the standard representation of $\varphi_\lambda$ as a $K$-integral in the noncompact case, and so one would hope to be able to prove Theorem \ref{cpctphibd} by applying the techniques of $\mathsection$\ref{phibounds} to this integral.  This works when $H$ is regular, and we use this approach to prove an asymptotic expansion for $\varphi_{t\mu}$ in Lemma \ref{Hregular}.  However, the fact that $b_{t\mu}$ is sharply concentrated along a flat subspace (in particular, that its absolute value has large derivatives) makes it difficult to prove bounds for $\varphi_{t\mu}(\exp(H))$ using the representation (\ref{beamint}) that are uniform as $H$ degenerates.  We get around these difficulties by observing that the terms in the expansion of Lemma \ref{Hregular} behave much more like plane waves on $G/K$ than the function $b_{t\mu}$, as their absolute values are not changing rapidly.  As a result, we may prove Theorem \ref{cpctphibd} by first averaging $b_{t\mu}$ under the action of a small open neighbourhood of the identity in $K$ to generate a plane wave on some open set in $S$, and then expressing $\varphi_{t\mu}$ as an average of the plane wave under rotation about a point in this set.

\subsection{The Structure of Gaussian Beams}
\label{compactbeam}

To begin this approach, we shall prove that $b_{t\mu}$ is localised around $\overline{T}$ at scale $t^{-1/2}$, making it a higher rank analogue of a Gaussian beam.  By Lemma \ref{extendA}, we may define

\begin{align*}
A^0 : T U^0_\C \rightarrow \ga_0 \\
g \mapsto \text{Re}(A(g)).
\end{align*}
It follows from Lemma \ref{extendA} that $A^0$ is left-invariant under $T$, so that $V A^0(e) = 0$ for $V \in i \ga_0$, and it may likewise be seen that $V A^0(e) = 0$ for $V \in i \q_0$.  This implies that when we restrict $\mu \circ A^0$ to $S$ it has a critical point at $e$, and hence along $\overline{T}$.  The following lemma shows that this critical point is negative definite transversally to $\overline{T}$.

\begin{lem}
\label{beamhess}
There are positive constants $C_1$ and $C_2$ depending only on $B^*$ such that for all $V \in i \q_0$, we have

\begin{equation*}
-C_1 \langle V, V \rangle \ge \frac{d^2}{dt^2} \mu( A^0( \exp( t V) ) ) \Big|_{t=0} \ge -C_2 \langle V, V \rangle \ge 0
\end{equation*}

\end{lem}

\begin{proof}

Let

\begin{align*}
V & = \sum_{\alpha \in \Delta^+} i c_\alpha( X_\alpha - X_{-\alpha}) \in i \q_0, \\
V^\pm & = \sum_{\alpha \in \Delta^+} i c_\alpha X_{\pm \alpha},
\end{align*}
where $c_\alpha \in \R$, $X_{\pm \alpha} \in \g_{0, \pm \alpha}$ and $X_{-\alpha} = \theta X_\alpha$.  Write the second order approximation to the Iwasawa decomposition of $\exp( t V)$ in terms of unknowns $V_1$, $V_2$ and $V_3$ as

\begin{equation*}
\exp(  tV) = \exp( 2t V^+ + t^2 V_1) \exp( t^2 V_2) \exp( -t( V^+ + V^-) + t^2 V_3 ) + O(t^3).
\end{equation*}
After applying the Baker-Campbell-Hausdorff formula to the RHS we have

\begin{equation*}
\exp(  t V) = \exp( t V - t^2 [ V^+, V^- ] + t^2 (V_1+V_2+V_3)) + O(t^3).
\end{equation*}
Equating coefficients gives $V_1+V_2+V_3 = [ V^+, V^- ]$, so that $V_2$ is the projection of $[ V^+, V^- ]$ to $\ga$.  Calculating this projection using the formula $[X_\alpha, X_{-\alpha}] = \langle X_\alpha, X_{-\alpha} \rangle H_\alpha$ gives

\bes
V_2 = - \sum_{\alpha \in \Delta^+} c_\alpha^2 \langle X_\alpha, X_{-\alpha} \rangle H_\alpha.
\ees
It follows that

\bes
\frac{d^2}{dt^2} \mu( A^0( \exp( t V) ) ) \Big|_{t=0} = -2 \sum_{\alpha \in \Delta^+} c_\alpha^2 \langle X_\alpha, X_{-\alpha} \rangle \langle \mu, \alpha \rangle,
\ees
and our assumption that $\mu \in B^*$ implies that $\langle \mu, \alpha \rangle \sim_{B^*} 1$.  Combining this with

\bes
\langle V, V \rangle = 2 \sum_{\alpha \in \Delta^+} c_\alpha^2 \langle X_\alpha, X_{-\alpha} \rangle
\ees
completes the proof.

\end{proof}

It follows from Lemma \ref{beamhess} that $b_{t\mu}$ has Gaussian decay at scale $t^{-1/2}$ transversally to $\overline{T}$, which implies that $b_{t\mu} \in L^2(S)$ has norm $\| b_{t\mu} \|_2 \gg t^{-(n-r)/4}$.  We next show that $b_{t\mu}$ decays rapidly in sets that are bounded away from $\overline{T}$ by an argument involving pseudodifferential operators.

\begin{prop}
\label{beamloc}
If $D \subset S$ is any compact set that does not intersect $\overline{T}$, we have

\begin{equation*}
| b_{t \mu}(x)| \ll_{D,A} t^{-A}, \qquad x \in D.
\end{equation*}

\end{prop}

\begin{proof}

Let $\Delta$ be the positive Laplacian on $S$ associated to the metric $- \langle \;, \, \rangle$ on $i\p_0$, which is equal to the restriction of the Casimir operator on $U$ to the space of right $K$-invariant functions.  Let $\mu_0 = \mu / \langle \mu, \mu \rangle^{1/2}$, and let $\partial H_\mu$ be the vector field on $S$ whose value at $uK$ is $(\partial / \partial t) \exp(itH_{\mu_0}) uK |_{t=0}$.  Under the isomorphism $TS \simeq U \times_K i \p_0$, $\partial H_\mu$ is given by $(u, \text{proj}_{i\p} ( \Ad_u^{-1} iH_{\mu_0}) )$.  The actions of $\Delta$ and $i \partial H_\mu$ on $b_{t\mu}$ are

\begin{align*}
\Delta b_{t\mu} & = \langle t \mu, t\mu \rangle, \\
i \partial H_\mu b_{t\mu} & = \langle t \mu, \mu_0 \rangle,
\end{align*}
and we shall prove the proposition by comparing these.  As we have already established that $\| b_{t\mu} \|_2 \gg t^{-(n-r)/4}$, it suffices to prove the proposition after first rescaling $b_{t\mu}$ to have $L^2$ norm one.

\begin{lem}

The principal symbol $p_0(x,\xi)$ of the operator $P_0 = \Delta - (i\partial H_\mu)^2$ satisfies $p_0(x,\xi) \ge 0$, and if $p_0(x,\xi) = 0$ then $x \in \overline{T}$ or $\xi = 0$.

\end{lem}

\begin{proof}

We shall denote the principal symbols of the operators $\Delta$ and $i \partial H_\mu$ by $p_\Delta$ and $p_\mu$.  Under the isomorphism $T^*S \simeq U \times_K i \p_0^*$, the formulas for $p_\Delta$ and $p_\mu$ are

\begin{align*}
p_\Delta & : (u, V) \mapsto -\langle V, V \rangle, \\
p_X & : (u, V) \mapsto V( \text{proj}_{i\p} ( \Ad_u^{-1} iH_{\mu_0})) = \langle V, \Ad_u^{-1} i \mu_0 \rangle.
\end{align*}
We then have

\begin{align*}
p_0(u,V) & = -\langle V, V \rangle - \langle V, \Ad_u^{-1} i \mu_0 \rangle^2,
\end{align*}
so that Cauchy-Schwarz implies that $p_0(u,V) \ge 0$ with equality iff $V = c \Ad_u^{-1} i\mu_0$ for some $c$.  Suppose that $p_0(u,V) = 0$, and assume that $0 \neq V \in i \ga_0^*$ so that $\Ad_u^{-1} i\mu \in i \ga_0^*$.  By following the proof of Proposition 8.8 (ii) in Chapter VII of \cite{He1}, and observing that that torus $T \subset P_*$ introduced there must be contained in our torus $T$ as $\mu$ is regular, we see that we must have $u = kt$ for $k \in K$ and $t \in T$.  We then have $\Ad_k^{-1} i \mu \in i \ga_0^*$, so that $k \in M'$ and $u = kt \in TK$ as required.

\end{proof}

As $P_0 b_{t\mu} = 0$ and $P_0$ is elliptic away from $\overline{T}$, it is a general principle of semiclassical analysis that $b_{t\mu}$ is rapidly decaying away from $\overline{T}$ as $t \rightarrow \infty$.  We shall give a quick proof of this fact.  Let $D \subset U_1 \subset U_2$ be open neighbourhoods of $D$ with $\overline{U}_1 \subset U_2$ and $\overline{U}_2 \cap \overline{T} = \emptyset$, and choose non-negative functions $a_1$ and $a_2$ in $C^\infty(S)$ satisfying

\begin{align*}
a_1(x) & = 1, \quad x \in \overline{T} \\
a_2(x) & = 1, \quad x \in U_2 \\
a_1 a_2 & \equiv 0.
\end{align*}
Define the operator $P$ by

\begin{equation*}
P = (1 + a_1) \Delta - ( i \partial H_\mu)^2
\end{equation*}
so that $P$ is elliptic on $S$.  If we define $P(a_2 b_{t\mu}) = \rho$ then we have $\text{supp}(\rho) \cap U_2 = \emptyset$.  As $P$ is an elliptic differential operator it has a parametrix $E$ such that $E P = I + \cS$ for some smoothing operator $\cS$ \cite[VI $\mathsection$4, 3.5]{St}, and applying $E P$ to $a_2 b_{t\mu}$ gives

\begin{equation}
\label{parametrix}
E \rho = a_2 b_{t\mu} + \cS(a_2 b_{t\mu}).
\end{equation}
As $E$ is also a pseudodifferential operator it is local up to smoothing, and because $\text{supp}(\rho) \cap U_2 = \emptyset$ this means there is a second smoothing operator $\cS_1$ such that $E\rho(x) = \cS_1\rho(x)$ for $x \in \overline{U}_1$.  Combining this with (\ref{parametrix}) for $x \in \overline{U}_1$ gives

\begin{align*}
(a_2 b_{t\mu})(x) & = E\rho(x) - \cS(a_2 b_{t\mu})(x) \\
b_{t\mu}(x) & = \cS_1 P(a_2 b_{t\mu})(x) - \cS(a_2 b_{t\mu})(x) \\
& = (\cS_1 P a_2 - \cS a_2)(b_{t\mu})(x)
\end{align*}
As $\cS_1 P a_2 - \cS a_2$ is a smoothing operator, this implies that the $L^2$ norm of $b_{t\mu}$ restricted to $U_1$ is rapidly decaying.  The standard methods of bounding $L^\infty$ norms of Laplace eigenfunctions in terms of their $L^2$ norms then imply that $|b_{t\mu}(x)| \ll_{D,A} t^{-A}$ for $x \in D$, which concludes the proof.

\end{proof}

Combining Lemma \ref{beamhess} and Proposition \ref{beamloc} on $U \cap T U^0_\C$, we obtain

\begin{cor}
\label{realA}

We have $\mu( A^0(u)) \le 0$ on $U \cap T U^0_\C$, with equality iff $u \in TK$.

\end{cor}

\subsection{Sharpness of Theorem \ref{main} in the Compact Case}
\label{compactsharp}

We may now prove that Theorem \ref{main} is sharp up to the logarithmic factor in the case of compact type.

The spherical function $\varphi_{t\mu}$ saturates the $L^p$ bounds for $p$ above the kink point.  To see this, first observe that $b_{t\mu}$ is roughly constant in a ball of radius $\gg t^{-1}$ about the identity in $S$ by (\ref{beamdef}), and so the expression (\ref{beamint}) for $\varphi_{t\mu}$ implies that $| \varphi_{t\mu}(s)| \gg 1$ in the same ball.  Moreover, the Weyl dimension formula implies that $\| \varphi_{t\mu} \|_2 \sim t^{-(n-r)/2}$.  These two facts imply that $t^{(n-r)/2} \varphi_{t\mu}$ has $L^2$ norm $\sim 1$, and has absolute value $\gg t^{(n-r)/2}$ on a set of measure $\gg t^{-n}$, so that $\| t^{(n-r)/2} \varphi_{t\mu} \|_p \gg t^{n(1/2 - 1/p) - r/2}$ as required.

Lemma \ref{beamhess} and Proposition \ref{beamloc} imply that the functions $t^{(n-r)/4} b_{t\mu}$ saturate the bounds of Theorem \ref{main} for $p$ below the critical point.  Indeed, by Proposition \ref{beamloc} it suffices to understand the behaviour of $b_{t\mu}$ in the open neighbourhood $U \cap T U_\C^0$ of $T$, and Lemma \ref{beamhess} implies that $|b_{t\mu}|$ is essentially the characteristic function of a ball of radius $t^{-1/2}$ around $\overline{T}$ in $S$.  It easily follows that the $L^p$ norm of $t^{(n-r)/4} b_{t\mu}$ is approximately $t^{(n-r)(1/2-1/p)/2}$.

\section{Bounds for spherical functions on compact groups}
\label{compact2}

In this section we shall derive Theorem \ref{cpctphibd} from the results of $\mathsection$\ref{compactbeam}, before using Theorem \ref{cpctphibd} to prove Theorem \ref{main}.

\subsection{Plane Waves and Integral Representations}
\label{sec61}

We begin by averaging $b_{t\mu}$ over rotations by a small neighbourhood of the identity in $K$ to generate a plane wave on $S$.  Let $B_1 \subset B \subset K$ be two open balls around $e$ that satisfy $B_1 = B_1^{-1}$ and $\overline{B_1} \subset B$.  Let $b_1 \in C^\infty_0(B)$ be a non-negative function that is equal to 1 on $B_1$, and define the function $\varphi_{t\mu}^0 \in C^\infty(U)$ by

\begin{equation*}
\varphi_{t\mu}^0(u) = \int_K b_1(k) b_{t\mu}(ku) dk.
\end{equation*}
To state the asymptotic we require for $\varphi_{t\mu}^0$, we introduce Cartan coordinates on $S$.  We define the map

\begin{align*}
\Phi : K/M \times i \ga_0 & \rightarrow S \\
(kM, H) & \mapsto k \exp(H).
\end{align*}
Define the diagram $D(U,K)$ and the regular set $\ga_r$ by

\bes
D(U,K) = \{ H \in i \ga_0 | \alpha(H) \in \pi i \Z \text{ for some } \alpha \in \Delta^+ \}, \qquad \ga_r = i\ga_0 \setminus D(U,K).
\ees
The regular set $\ga_r$ is a union of open simplices, and we choose one such simplex $P_0$ whose closure contains the origin.  It is known (see Theorem 3.3, Chapter VII of \cite{He1}) that $\Phi(K/M, D(U,K))$ is an analytic set of codimension at least 2 in $S$, and we define the regular set $S_r$ to be $S \setminus \Phi(K/M, D(U,K))$.

\begin{prop}
\label{cartan}

We have $\Phi(K/M, P_0) = S_r$, and the map $\Phi: K/M \times P_0 \rightarrow S_r$ is a covering map.  Moreover, if we have

\be
\label{cartancoincide}
u = k_1 \exp(H_1) k_2 = l_1 \exp(H_2) l_2
\ee
with $H_i \in P_0$ and $k_i, l_i \in K$, then $H_1 = H_2$, $k_1 M' = l_1 M'$, and $M' k_2 = M' l_2$.

\end{prop}

\begin{proof}

The assertions that $\Phi(K/M, P_0) = S_r$, $\Phi$ is a covering map, and $H_1 = H_2$ in (\ref{cartancoincide}) are proven in Lemma 8.1 and Theorem 8.6, Chapter VII of \cite{He1}.  To prove that $k_1 M' = l_1 M'$, consider

\bes
u (\theta u)^{-1} = k_1 \exp(2H_1) k_1^{-1} = l_1 \exp(2H_1) l_1^{-1},
\ees
so that if we let $k = l_1^{-1} k_1$ then we have $k \exp(2H_1) k^{-1} = \exp(2H_1) = s$.  Following the proof of Lemma 8.7 of \cite{He1}, Chapter VII, we let $Z_s$ be the centraliser of $s$ in $U$.  Our assumption that $H_1$ was in $\ga_r$ implies that the Lie algebra of $Z_s$ is exactly $\gm_0 + i \ga_0$, and so reasoning as in the proof of that lemma we see that $k \in M'$ as required.  The claim that $M' k_2 = M' l_2$ follows in the same way.

\end{proof}

Proposition \ref{cartan} implies that we may define the Cartan $A$-coordinate $\rho: S_r \rightarrow P_0$ by $k \exp(H) \mapsto H$.  Assume that $B B_1 \cap M' \subseteq M_0$.  If we choose $B$ and $B_1$ to be sufficiently small, there exists an open set $Q_0 \subset P_0$ with the property that $B B_1 \exp(Q_0) \subset U_\C^0$.  Define $V = \Phi( B_1, Q_0) \subset S_r$.  We may assume that $B_1$ and $B$ are small enough that $\Phi$ provides a diffeomorphism $V \simeq B_1 M / M \times Q_0$.

\begin{lem}
\label{planewavelemma}

There is a function $a \in C^\infty( \R \times V \times B^*)$ with an asymptotic expansion

\bes
a(t, s, \mu) = \sum_{i=0}^\infty t^{-i} a_i(s, \mu)
\ees
that converges locally uniformly, with $a_i \in C^\infty( V \times B^*)$ and $a_0$ nonvanishing, such that we have

\begin{equation}
\label{planewave}
\varphi_{t\mu}^0(s) = t^{-(n-r)/2} a(t, s, \mu) e^{-t\mu(\rho(s))}, \quad s \in V.
\end{equation}

\end{lem}

\begin{proof}

Our asumption that $B B_1 \exp(Q_0) \subset U^0_\C$ implies that we may define $\phi \in C^\infty( B B_1 \times Q_0 \times B^*)$ by $\phi(k, H, \mu) = -\mu( A( k \exp(H)) )$.  Corollary \ref{realA} implies that $\text{Re} \; \phi \le 0$ with equality iff $k \exp(H) \in \overline{T}$.  Theorem 8.3 (iii), Section 8, Chapter VII of \cite{He1} implies that $\overline{T} = \Phi(M', \overline{P_0})$, and Proposition \ref{cartan} and our assumption that $B B_1 \cap M' \subseteq M_0$ then imply that $k \exp(H) \in \overline{T}$ iff $k \in M_0$.  If $k_1 \in B_1$ and $H \in Q_0$, we may write

\begin{align*}
\varphi_{t\mu}^0(k_1 \exp(H)) & = \int_K b_1(k) e^{ -t\mu( A( k k_1 \exp(H)))} dk \\
& = \int_K b_1(k k_1^{-1}) e^{t \phi(k, H, \mu)} dk.
\end{align*}
It will be convenient to reduce this integral to one with a single critical point.  Lemma \ref{MA} implies that $\phi$ is left-invariant under $M_0$, so that we may define $R = M_0 \backslash K$ and reduce $\phi$ to a function on $M_0 \backslash M_0BB_1 \subset R$, which we continue to denote by $\phi$.  We also define $b_1' : R \rightarrow \R$ by

\bes
b_1'(M_0k) = \int_{M_0} b_1( m k ) dm.
\ees
The support of $b_1'$ is contained in $M_0 \backslash M_0B$, and our assumption that $b_1$ equals 1 on $B_1 = B_1^{-1}$ implies that $b_1'(M_0 k_1^{-1}) > 0$ for all $k_1 \in B_1$.  We have

\begin{equation}
\label{phi0R}
\varphi_{t\mu}^0(k_1 \exp(H) ) = \int_R b_1'(r k_1^{-1}) e^{t \phi(r, H, \mu)} dr.
\end{equation}
If $D\phi$ is the Hessian of $\phi$ at $e$, we may calculate $D\phi$ with respect to the basis $\{ Y_\alpha | \alpha \in \tD^+ \}$ of $T_e R$ as in Proposition \ref{hessian} to be the diagonal matrix

\begin{equation*}
(D \phi)_{\alpha \alpha} = \tfrac{1}{2} \langle \mu, \alpha \rangle ( e^{2\alpha(H)} -1 ).
\end{equation*}
We have $\alpha(H) \in i \R \setminus \pi i \Z$ for all $\alpha$ when $H \in \ga_r$, so that $\text{Re}( D \phi_{\alpha \alpha} ) < 0$ on $Q_0 \times B^*$.  We may apply the stationary phase method for complex phases (see for instance \cite[Theorem 7.75]{Ho2}) to obtain a function $a \in C^\infty( \R \times B_1 \times Q_0 \times B^*)$ with an asymptotic expansion

\bes
a(t, k, H, \mu) = \sum_{i=0}^\infty t^{-i} a_i(k, H, \mu)
\ees
that converges locally uniformly, with $a_i \in C^\infty( B_1 \times Q_0 \times B^*)$, such that we have

\bes
\varphi_{t\mu}^0( k \exp(H) ) = t^{-(n-r)/2} a(t, k, H, \mu) e^{-t\mu(\rho((s))}, \quad k \in B_1, H \in Q_0.
\ees
Moreover, the condition that $b_1'(M_0 k_1^{-1}) > 0$ for all $k_1 \in B_1$ implies that $a_0$ is nonvanishing.  The functions $a$ and $a_i$ must be right-invariant under $M$, so we may push them forward under $\Phi$ to obtain functions $a \in C^\infty( \R \times V \times B^*)$ and $a_i \in C^\infty( V \times B^*)$ as required.

\end{proof}

Arguing in the same way allows us to prove an asymptotic expansion for $\varphi_{t\mu}(\exp(H))$ when $H$ is regular.  We now let $B \subset i\ga_0$ denote a ball around the origin such that $\exp(\Ad_K B) \subset U_\C^0$, and let $B_r  = B \cap \ga_r$.

\begin{lem}
\label{Hregular}

There is a function $a \in C^\infty( \R \times B_r \times B^* \times W)$ with an asymptotic expansion

\bes
a(t, H, \mu, w) = \sum_{i=0}^\infty t^{-i} a_i(H, \mu, w)
\ees
that converges locally uniformly, with $a_i \in C^\infty( B_r \times B^* \times W)$ and $a_0$ nonvanishing, such that we have

\begin{equation*}
\varphi_{t\mu}(\exp(H)) = t^{-(n-r)/2} \sum_{w \in W} a(t, H, \mu, w) e^{-t\mu(wH)} \quad \text{for } H \in B_r.
\end{equation*}

\end{lem}

\begin{proof}

We write

\begin{align}
\notag
\varphi_{t\mu}(\exp(H)) & = \int_K b_{t\mu}(\exp(\Ad_k H)) dk \\
\label{Hregularint}
& = \int_K e^{-\mu(A(\exp(\Ad_k H)))} dk.
\end{align}
Our assumption on $B$ implies that we may define $\phi'(k,H,\mu) = -\mu(A(\exp(\Ad_k H)))$.  The function $\phi'$ is clearly right-invariant under $M$, so that we may reduce this integral to one on $K / M$, and it satisfies $\text{Re}(\phi') \le 0$ with equality iff $k \exp(H) k^{-1} \in TK$.  When $H \in \ga_r$ we again have $k \exp(H) k^{-1} \in TK$ iff $k \in M'$, and so by Proposition \ref{beamloc} it suffices to consider neighbourhoods of the Weyl points in the integral (\ref{Hregularint}).  The lemma now follows from stationary phase as before.

\end{proof}

Lemma \ref{planewavelemma} shows that there is a clear similarity between $\varphi_{t\mu}^0$ and the plane waves $e^{i\lambda(A(g))}$ on $G/K$.  We will make use of this by choosing $H_0 \in Q_0$, letting $h = \exp(H_0) \in V$, and expressing $\varphi_{t\mu}$ as an average of $\varphi_{t\mu}^0$ under rotation about $h$.  The fact that $a_0$ in Lemma \ref{planewavelemma} was nonvanishing means that, for $t$ sufficiently large, we may normalise $\varphi_{t\mu}^0$ by setting $\varphi_{t\mu}^0(h) = 1$.  We then have

\be
\label{planewavephi}
\varphi_{t\mu}( \exp(H) ) = \int_K \varphi_{t\mu}^0(h k \exp(H) ) dk,
\ee
and when $s \in V$ we have an asymptotic

\bes
\varphi_{t\mu}^0(s) = a(t, s, \mu) e^{-t\mu(\rho(s))}
\ees
with $a$ as in Lemma \ref{planewavelemma}.  If $B \in \ga_0$ is a ball such that $h K \exp(iB) \subset V$, we may define $\psi \in C^\infty(K \times B \times B^*)$ by

\bes
\psi(k, H, \mu) = i \mu( \rho(h k \exp(iH) ) ).
\ees
We choose to multiply by $i$ in this way so that $\psi$ and its domain are both real.  When $H \in B$, this allows us to rewrite (\ref{planewavephi}) as

\bes
\varphi_{t\mu}( \exp(iH) ) = \int_K a(t, h k \exp(iH), \mu) e^{it \psi(k, H, \mu)} dk.
\ees
By applying the asymptotic expansion of $a$, we see that Theorem \ref{phithm} will follow from

\begin{prop}

We have

\bes
\label{cpctKint}
\int_K a(h k \exp(iH), \mu) e^{i t \psi(k, H, \mu)} dk \ll \prod_{\alpha \in \tD^+} ( 1 + t |\alpha(H)| )^{-1/2}
\ees
for all $H \in B$ and $a \in C^\infty(V, B^*)$.

\end{prop}

\subsection{The critical set of $\psi$}
\label{psicrit}

We shall prove Proposition \ref{cpctKint} by uniformising $\psi$ as in $\mathsection$\ref{phibounds}.  We begin by establishing the following analogues of Propositions \ref{critical} and \ref{hessian} in the compact case.  We recall that $K_H$ is the stabiliser of $H$ in $K$.

\begin{prop}
\label{cpctphicrit}

The function $\psi(k, H, \mu)$ is right invariant under $K_H$, and its critical point set is equal to $M' K_H$.

\end{prop}

\begin{proof}

The invariance of $\psi$ under $K_H$ is immediate.  To determine the critical point set, we shall first assume that $k$ is a critical point of $\psi$ and show that $k = w k_H$ for some $w \in M'$ and $k_H \in K_H$.  Choose a vector $X \in \gk_0$ and use the diffeomorphism $\Phi : V \simeq B_1 M/M \times Q_0$ to write $h k \exp( t X) \exp(iH)$ as

\bes
h k \exp( t X) \exp(iH) = \overline{k}_1(t) \exp(V(t)) \in V
\ees
for $t$ near 0, where $\overline{k}_1(t)$ and $V(t)$ are smooth functions that take values in $K/M$ and $P_0$ respectively.  We have

\bes
(\partial / \partial t) \psi(k \exp(tX), H, \mu) \Big|_{t = 0} = i \mu( V'(0) ).
\ees
Let $s : K/M \rightarrow K$ be a section of the quotient map that is defined in a neighbourhood of $\overline{k}_1(0)$.  If we define $k_1(t) = s( \overline{k}_1(t))$ for $t$ near 0, this gives a smooth Cartan decomposition

\begin{equation}
\label{cartdecomp}
h k \exp( t X) \exp(iH) = k_1(t) \exp(V(t)) k_2(t) \in U.
\end{equation}
If we set $a = \exp(V(0))$ and $k_i = k_i(0)$ for $i = 1,2$, and define $X_i \in \gk_0$ so that

\bes
k_i(t) = k_i(0) \exp( t X_i + O(t^2) ),
\ees
then differentiating (\ref{cartdecomp}) at $t = 0$ gives

\begin{align*}
\Ad_{\exp(iH)}^{-1} X & = \Ad_{a k_2}^{-1} X_1 + \Ad_{k_2}^{-1} V'(0) + X_2 \\
\Ad_{k_2} \Ad_{\exp(iH)}^{-1} X & = \Ad_a^{-1} X_1 + V'(0) + \Ad_{k_2} X_2.
\end{align*}
As $\Ad_a^{-1} X_1$ and $\Ad_{k_2} X_2$ both lie in $\gk + \q$, we see that $V'(0)$ is equal to the projection of $\Ad_{k_2} \Ad_{\exp(iH)}^{-1} X$ to $\ga$.  Our assumption that $k$ is a critical point of $\psi$ then implies that

\begin{equation}
\label{critcondition}
\langle H_\mu,  \Ad_{k_2} \Ad_{\exp(iH)}^{-1} X \rangle = 0 \quad \text{for all } X \in \gk.
\end{equation}

Let $\gl_0$ be the centraliser of $H$ in $\g_0$, and let $L$ be the Levi subgroup of $G$ with Lie algebra $\gl_0$.  $\gl_0$ is stable under $\theta$, and we write its Cartan decomposition as $\gl_0 = \gk_{L,0} + \p_{L,0}$.  Let $K_L = L \cap K$, which is a maximal compact subgroup of $L$ because it is compact and has Lie algebra $\gk_{L,0}$.  We note that $K_L \subseteq K_H$.  After shrinking $B$ if necessary, our assumption that $H \in B$ implies that the projection of $\Ad_{\exp(iH)}^{-1} \gk$ to $\p$ is equal to $\p_L^\perp \subset \p$, and so condition (\ref{critcondition}) holds iff $\Ad_{k_2}^{-1} H_{\mu} \in \p_L$.  The inclusion $\Ad_{k_2}^{-1} H_{\mu} \in \p_L$ implies that there is an element $k_L \in K_L$ such that $\Ad_{k_L} \Ad_{k_2}^{-1} H_{\mu} \in \ga$, and as $H_\mu$ is regular this implies that $k_2 = w k_L$ for some $w \in M'$.  Substituting this into (\ref{cartdecomp}) at $t = 0$ gives

\begin{align*}
h k & = k_1 a k_2 \exp(-iH) \\
& = k_1 a w k_L \exp(-iH) \\
& = k_1 a \exp( -i wH) w k_L.
\end{align*}
We have $h = \exp(H_0)$ and $a \exp( -i wH) = \exp( V(0) - iwH)$, and if $B$ is chosen small enough we will have both $H_0 \in P_0$ and $V(0) - iwH \in P_0$.  Proposition \ref{cartan} then gives $k \in M' k_L$ as required.  This shows that the critical point set of $\psi$ is contained in $M' K_H$, and the reverse inclusion follows from the right $K_H$-invariance of $\psi$ and the easily observed fact that $\psi$ is critical on $M'$.

\end{proof}

Choose $l \in K$, and let $\ga_L$ and $K_L$ be as in $\mathsection$\ref{phinotation}.  We again write $l = w l_0$ with $l_0 \in K_L$ and $w \in M'$ fixed, and let $X_\alpha = \Ad_l^{-1} Y_\alpha$.  It follows from Proposition \ref{cpctphicrit}, just as in the proof of Lemma \ref{al}, that $l$ is a critical point of $\phi$ exactly when $H \in \ga_{0,L} \cap B$.

\begin{prop}
\label{cpctphihess}

There are positive analytic functions 

\begin{equation*}
F_\alpha : B \cap \ga_{0,L} \rightarrow \R, \quad \alpha \in \tD^+,
\end{equation*}
such that when $H \in B \cap \ga_{0,L}$, the Hessian of $\psi$ at $l$ with respect to the vector fields $\{ X_\alpha | \alpha \in \tD_0^+ \}$ is diagonal, and satisfies

\begin{align*}
(D \psi)_{\alpha \alpha} & = \langle \mu, \alpha \rangle F_\alpha(H) \sin(\alpha(wH)), \quad \alpha \in \tD^+, \\
(D \psi)_{\alpha \alpha} & = 0, \quad \alpha \in \tD^+_0 \setminus \tD^+.
\end{align*}

\end{prop}

\begin{proof}

Showing that $(D\psi)_{\alpha \beta} = 0$ when $\alpha$ or $\beta$ lie over the zero root is simple, and left to the reader.  Let $\alpha, \beta \in \tD^+$.  We wish to calculate

\bes
\frac{\partial^2}{\partial s \partial t} \mu( \rho( h l \exp(sX_\alpha) \exp(tX_\beta) \exp(iH) ) ) \Big|_{s=t=0}.
\ees
As $l = w l_0$ and $H \in \ga_L$, we may rewrite the argument of $\rho$ above as

\bes
h l \exp(sX_\alpha) \exp(tX_\beta) \exp(iH) = h \exp(s Y_\alpha)\exp(t Y_\beta) \exp(iwH) l,
\ees
and so we may instead calculate

\bes
\frac{\partial^2}{\partial s \partial t} \mu( \rho( h \exp(s Y_\alpha)\exp(t Y_\beta) \exp(iwH) ) ) \Big|_{s=t=0}.
\ees
As in the proof of Proposition \ref{cpctphicrit}, we choose a smooth Cartan decomposition

\be
\label{cartan1}
h \exp(s Y_\alpha)\exp(t Y_\beta) \exp(iwH) = k_1(s,t) \exp( V(s,t)) k_2(s,t)
\ee
for $s$ and $t$ near 0.  Moreover, because $h \exp( i wH) = \exp(H_0 + iwH)$ and $H_0 + iwH \in P_0$ for $B$ small, we may assume that $k_1(0,0) = k_2(0,0) = e$ so that we can write

\bes
k_i(s,t) = \exp( X_{i,s}s + X_{i,t}t + X_{i,st} st + O(s^2) + O(t^2) ).
\ees
We define $V = V(0,0) = H_0 + iwH$ and $V_{st} = V_{st}(0,0)$, and let $a = \exp(V)$.  Writing the approximation to (\ref{cartan1}) involving terms $s$, $t$ and $st$ gives

\begin{align*}
h \exp(s Y_\alpha)\exp(t Y_\beta) \exp(iwH) & = \exp( X_{1,s}s + X_{1,t}t + X_{1,st} st) a \exp( V_{st} st) \\
& \qquad \exp( X_{2,s}s + X_{2,t}t + X_{2,st} st) + O(s^2) + O(t^2)\\
a \exp(s \Ad_{\exp(iwH)}^{-1} Y_\alpha) \exp(t \Ad_{\exp(iwH)}^{-1} Y_\beta) & = a \exp( \Ad_a^{-1} [X_{1,s}s + X_{1,t}t + X_{1,st} st ]) \exp( V_{st} st) \\
& \qquad \exp( X_{2,s}s + X_{2,t}t + X_{2,st} st) + O(s^2) + O(t^2)\\
\end{align*}
Combining exponentials using Baker-Campbell-Hausdorff gives

\begin{multline}
\label{secondorder}
\exp( \Ad_{\exp(iwH)}^{-1} (sY_\alpha + t Y_\beta + st/2 [Y_\alpha, Y_\beta] ) ) = \exp ((\Ad_a^{-1} X_{1,s} + X_{2,s})s + (\Ad_a^{-1} X_{1,t} + X_{2,t})t \\
+ (V_{st} + \Ad_a^{-1} X_{1,st} + X_{2,st} +1/2[ \Ad_a^{-1} X_{1,s}, X_{2,t}] + 1/2[ \Ad_a^{-1} X_{1,t}, X_{2,s}] ) st + O(s^2) + O(t^2) ),
\end{multline}
and equating first order terms, we have

\be
\label{Yequation}
\Ad_{\exp(iwH)}^{-1} Y_\alpha = \Ad_a^{-1} X_{1,s} + X_{2,s}, \qquad \Ad_{\exp(iwH)}^{-1} Y_\beta = \Ad_a^{-1} X_{1,t} + X_{2,t}.
\ee
If we let $Y_\alpha = V_\alpha + V_{-\alpha}$ where $V_{\pm \alpha} \in \g_{0,\pm \alpha}$, and likewise for $\beta$, then (\ref{Yequation}) becomes

\bes
e^{-i \alpha(wH)} V_\alpha + e^{i \alpha(wH)} V_{-\alpha} = \Ad_a^{-1} X_{1,s} + X_{2,s}, \qquad e^{-i \beta(wH)} V_\beta + e^{i \beta(wH)} V_{-w\beta} = \Ad_a^{-1} X_{1,t} + X_{2,t}.
\ees
Because $V = H_0 + iwH$ was generic we may solve this to obtain

\begin{equation}
\label{kdiff}
X_{1,s} \in \frac{ \sin(\alpha(wH)) }{ \sin( \alpha(V)/i) } (V_\alpha + V_{-\alpha}) + \gm_0, \qquad X_{2,s} \in \frac{ \sin( \alpha(V)/i - \alpha(wH) ) }{ \sin( \alpha(V)/i) } (V_\alpha + V_{-\alpha}) + \gm_0,
\end{equation}
and likewise for $X_{i,t}$ and $\beta$.

Equating the $st$ terms in (\ref{secondorder}) gives

\bes
1/2 \Ad_{\exp(iwH)}^{-1} [Y_\alpha, Y_\beta] = V_{st} + \Ad_a^{-1} X_{1,st} + X_{2,st} +1/2[ \Ad_a^{-1} X_{1,s}, X_{2,t}] + 1/2[ \Ad_a^{-1} X_{1,t}, X_{2,s}].
\ees
We have $ [Y_\alpha, Y_\beta] \in \gk$, so that $\Ad_{\exp(iwH)}^{-1} [Y_\alpha, Y_\beta]$ and $\Ad_a^{-1} X_{1,st} + X_{2,st}$ both lie in $\gk + \q$.  This implies that

\be
\label{vstproj}
V_{st} = - \text{proj}_\ga \left( 1/2[ \Ad_a^{-1} X_{1,s}, X_{2,t}] + 1/2[ \Ad_a^{-1} X_{1,t}, X_{2,s}] \right),
\ee
where $\text{proj}_\ga$ is the orthogonal projection onto $\ga$.

We first consider the case where $\alpha \neq \beta$.  If $\alpha \neq \beta$ as elements of $\Delta^+$ rather than just $\tD^+$, we see that (\ref{vstproj}) must vanish because the commutators of the form $[V_\alpha, V_{-\beta}]$ must lie in root spaces corresponding to nonzero roots.  If $\alpha = \beta$ in $\Delta^+$, the vanishing of (\ref{vstproj}) follows from our assumption that the vectors $V_{\pm \alpha}$, $V_{\pm \beta}$ were orthogonal, and the identity $[I,J] = H_\gamma \langle I,J \rangle$ for $I \in \g_\gamma$ and $J \in \g_{-\gamma}$.

We now assume that $\alpha = \beta$.  In this case, $X_{i,s} = X_{i,t}$ so that (\ref{vstproj}) becomes

\bes
V_{st} = - \text{proj}_\ga \left( [ \Ad_a^{-1} X_{1,s}, X_{2,s}] \right).
\ees
Substituting the values of $X_{1,s}$ and $X_{2,s}$ from (\ref{kdiff}) and noting that $[\gm, \gk] \perp \ga$, we have

\begin{align*}
-[\Ad_a^{-1} X_{1,s}, X_{2,s}] & \in - \frac{ \sin(\alpha(wH)) \sin( \alpha(V)/i - \alpha(wH) ) }{ \sin( \alpha(V)/i)^2 } \\
& \qquad \qquad \times [ e^{-\alpha(V)} V_\alpha + e^{\alpha(V)} V_{-\alpha}, V_\alpha + V_{-\alpha} ] + \ga^\perp \\
& \in \frac{ 2i \sin( \alpha(wH)) \sin( \alpha(V)/i - \alpha(wH) ) }{ \sin( \alpha(V)/i) } [V_\alpha, V_{-\alpha} ] + \ga^\perp \\
& \in \frac{ 2i \sin( \alpha(wH)) \sin( \alpha(V)/i - \alpha(wH) ) }{ \sin( \alpha(V)/i) } H_\alpha \langle V_\alpha, V_{-\alpha} \rangle + \ga^\perp.
\end{align*}
As $\langle V_{\alpha}, V_{-\alpha} \rangle = -1/2$, we therefore have

\bes
V_{st} = \frac{ -i \sin( \alpha(wH)) \sin( \alpha(V)/i - \alpha(wH) ) }{ \sin( \alpha(V)/i) } H_\alpha.
\ees
If we define

\bes
F_\alpha(H) = \frac{\sin( \alpha(V)/i - \alpha(wH) ) }{ \sin( \alpha(V)/i) },
\ees
then if $B$ is sufficiently small, $F_\alpha$ will be a positive real analytic function on $B$.  We then have

\begin{align*}
(d / dt)^2 \psi( l \exp(tX_\alpha), H, \mu) \Big|_{t=0} & = i (d / dt)^2 \mu( \rho( h l \exp(tX_\alpha) \exp(iH) ) ) \Big|_{t=0} \\
& = i \mu(V_{st}) \\
& = \sin( \alpha(wH)) F_\alpha(H) \langle \mu, \alpha \rangle,
\end{align*}
which completes the proof.

\end{proof}

\subsection{Uniformisation of $\psi$}
\label{compactsph}

Proposition \ref{cpctKint}, and hence Theorem \ref{cpctphibd}, follows as in $\mathsection$\ref{phiphase} after proving analogues of Theorem \ref{uniform1} and Proposition \ref{uniform3} for $\psi$.  The analogue of Proposition \ref{uniform3} follows in a straightforward way from Propositions \ref{cpctphicrit} and \ref{cpctphihess}, but adapting Theorem \ref{uniform1} requires some comments.  Choose $l \in K$ with $l \notin M'$, and a flag $F \in \cF$, and retain all the notation of $\mathsection$\ref{phinotation} and $\mathsection$\ref{complexification}, including a choice of $s = (u,\mu) \in (\pi^{-1}_\cA(B) \cap \cJ \cap \cA_L) \times \ga_0^*$ with  $\mu$ regular.  We now denote points in $(S,s)$ by $s' = (u', \mu')$.  If $x' \in X$, we let $\mu(x')$ denote its projection to $\ga^*$.  We may apply Propositions \ref{cpctphicrit} and \ref{cpctphihess} exactly as in $\mathsection$\ref{phiuniform} to prove the following analogue of Corollary \ref{uniformcor}.

\begin{prop}
\label{cpctuniform2}

There exists a subspace $(Y, x) \subset (X, x)$ and an isomorphism

\bes
\xymatrix{ (X, x) \ar[rr]^f \ar[dr]_{\pi_S}&& (Y,x) \times (\C^d, 0)  \ar[dl]^{\pi_S \times 0}\\& (S, s)}
\ees
with the following properties:

\begin{enumerate}[(a)]

\item

$f|_Y$ is the identity.

\item
$(Y, x)$ is invariant under $c$, and $f$ commutes with $c$.

\item
The projection $(Y,x) \rightarrow (S, s)$ is regular.

\item
We have $Y_{s'} = l K_{q,\C}$ when $s' \in S_q$, and $Y_{s'} \subseteq l K_{i,\C}$ when $s' \in S_i$ with $0 \le i < q$.

\item

We have

\be
f_* \psi(y,z) = \psi(y) - \sum_{\alpha \in \tE_q} \langle \mu(y), \alpha \rangle w^{-1} \alpha_X(y) z_\alpha^2.
\ee

\end{enumerate}

\end{prop}

Next, we derive the analogue of Theorem \ref{uniform1} from Proposition \ref{cpctuniform2}, which completes the proof of Theorem \ref{cpctphibd}.

\begin{thm}
\label{cpctuniform1}

There is an isomorphism $f$

\bes
\xymatrix{ (X, x) \ar[rr]^f \ar[dr]_{\pi_S}&&  (\C^d \times \C^{d'}, 0) \times (S, s) \ar[dl]^{0 \times \textup{id}}\\& (S, s)},
\ees
a function $\psi_S \in \cO( S, s)$, and a non-constant affine-linear map $L : \C^{d'} \longrightarrow \C$, such that $f$, $\psi_S$ and $L$ all commute with $c$, and such that

\bes
f_* \psi(z, z', s') = \psi_S(s') - \sum_{\alpha \in \tE_q} \langle \mu', \alpha \rangle w^{-1}\alpha(u') z_\alpha^2 + Q(u')L(z').
\ees

\end{thm}

\begin{proof}

The deduction of Theorem \ref{cpctuniform1} from Proposition \ref{cpctuniform2} follows much as in the case of noncompact type.  After proving the analogue of Lemma \ref{Qdiff}, we are given $H_L \in \ga_{0,L}$ and $H^L \in \ga_0 \setminus \ga_{0,L}$, and have to show that $(\partial / \partial t) Y_\alpha \psi(l, H_L + t H^L, \mu) \neq 0$ for some $\alpha \in \tD_L^+$.  Let $\alpha \in \tD_L^+$, and for small $s, t \in \R$, choose a smooth Cartan decomposition

\bes
h l \exp( s Y_\alpha) \exp( H_L + tH^L) = k_1(s,t) a(s,t) k_2(s,t)
\ees
with the property that $k_2(0,0) = l$.  Reasoning as in the proof of Proposition \ref{cpctphicrit} with $X = Y_\alpha = V_\alpha + V_{-\alpha}$, we have

\bes
Y_\alpha \psi( l, H_L + tH^L, \mu) = i \langle \Ad_{k_2(0,t)}^{-1} H_\mu, e^{-i t \alpha(H^L)} V_\alpha + e^{i t \alpha(H^L)} V_{-\alpha} \rangle.
\ees
As in the proof of Lemma \ref{Wpsi1}, this gives

\bes
(\partial / \partial t) Y_\alpha \psi(l, H_L + tH^L, \mu)\big|_{t=0} = \alpha(H^L) \langle \Ad_{l}^{-1} H_\mu, V_{\alpha} - V_{-\alpha} \rangle.
\ees
Lemma \ref{Wpsi2} implies that this quantity is nonzero for some $\alpha \in \tD_L^+$, which completes the proof.

\end{proof}

\subsection{Bounds for $L^p$ norms in compact type}
\label{cpctLp}

Let $B_1 \subset i \p_0$ be a round ball around the origin with respect to the Killing form such that $2B_1 \cap i \ga_0 \subset B$ where $B$ is as in Theorem \ref{cpctphibd}.  Let $b \in C^\infty(S)$ be a nonnegative real valued $K$-biinvariant function with support in $\exp(B_1)$, and that satisfies $b(e) = 1$ and $b(u) = b(u^{-1})$.  Let $k_{t}^0 = t^{n-r} b \varphi_{s^* t\mu}$, let $K_t^0$ be the point pair invariant kernel on $S$ associated to $k_t^0$, and let $T_t$ be the operator with integral kernel $K_t^0$.

\begin{prop}
\label{kmutransform}

The spherical transform of $k_t^0$ is real, and satisfies $\widehat{k}_t^0(s^*t\mu) \gg 1$ and $\widehat{k}_t^0(t\nu) \ll_{A, \delta} t^{-A}$ if $\| s^* \mu - \nu \| > \delta$.

\end{prop}

\begin{proof}

To prove that $\widehat{k_t^0}(\lambda)$ is real, the identity $\varphi_\lambda(u^{-1}) = \overline{\varphi}_{\lambda}(u)$ implies that

\begin{align*}
\overline{\widehat{k_t^0}}(\lambda) & = t^{n-r} \int_U b(u) \overline{\varphi}_{s^* t\mu}(u) \overline{\varphi}_{s^* \lambda}(u) du \\
& = t^{n-r} \int_U b(u) \varphi_{s^* t\mu}(u^{-1}) \varphi_{s^* \lambda}(u^{-1}) du \\
& = t^{n-r} \int_U b(u) \varphi_{s^* t\mu}(u) \varphi_{s^* \lambda}(u) du \\
& = \widehat{k}_t^0(\lambda)
\end{align*}
as required.  The assertion that $\widehat{k}_t^0(s^* t\mu) \gg 1$ follows in a similar way from $\varphi_{s^* t\mu} = \overline{\varphi}_{t\mu}$ and Lemma \ref{Hregular}, which implies that $t^{n-r} |\varphi_{s^* t\mu}|^2$ has mass $\gg 1$ in any ball about the origin.

We prove the last assertion in the same way as Proposition \ref{phiphi}, using the integral representation (\ref{planewavephi}).  We have

\begin{align*}
\widehat{k}_t^0( t\nu) = \int_S t^{n-r} b(s) \varphi_{s^* t\mu}(s) \varphi_{s^* t\nu}(s) ds,
\end{align*}
and after substituting the representation (\ref{planewavephi}), this becomes

\bes
\widehat{k}_t^0( \nu) = \int_K \int_S t^{n-r} b(s) \varphi_{ts^* \mu}^0(h k s) \varphi_{ts^* \nu}^0(h s) ds dk.
\ees
Our assumption on $B_1$ implies that $h k s \in V$ for $s \in \text{supp}(b)$, and so we may apply the asymptotic expansion of Lemma \ref{planewavelemma} which reduces us to proving that

\bes
\int_K \int_S t^{n-r} b(s) a_1( hks, s^* \mu) a_2(hs, s^* \nu) \exp(- ts^* \mu( \rho(h k s)) - ts^* \nu( \rho(h s)) ) ds dk \ll_{A,\delta} t^{-A}
\ees
Under the identification of $T^* S$ with $U \times_K i\p^*$, the differentials of $- s^* \mu(\rho(s))$ and $-s^* \nu(\rho(s))$ lie in $U \times \Ad_K i\mu$ and $U \times \Ad_K i\nu$ respectively.  The assumption $\| s^* \mu - \nu \| > \delta$ implies that $U \times \Ad_K i\mu$ and $U \times \Ad_K (-i\nu)$ are separated, and result now follows from integration by parts as in Proposition \ref{phiphi}.

\end{proof}

It follows from Proposition \ref{kmutransform} that $T_t$ is a selfadjoint approximate spectral projector onto the parameter $t\mu$.  It follows that if we define $k_t = k_t^0 * k_t^0$ and let $K_t$ be the point pair invariant associated to $k_t$, then $K_t$ is the integral kernel of $T_t T^*_t$.  We may prove Theorem \ref{main} as in the noncompact case, by performing a radial decomposition of $K_t$ and estimating the $L^1 \rightarrow L^\infty$ and $L^2 \rightarrow L^2$ norms of the truncated pieces.  This works in exactly the same way once we have a pointwise bound for $k_t$ analogous to that of Lemma \ref{kpoint}, and a bound for the Harish-Chandra transform of the truncated pieces of $k_t$.  The pointwise bound is given by the following lemma.

\begin{lem}
\label{kbound}

We have

\bes
k_t(\exp(H)) \ll t^{n-r} \prod_{\alpha \in \tD^+} (1 + t|\alpha(H)|)^{-1/2},
\ees
uniformly for $H \in B$ and $\mu \in B^*$.

\end{lem}

\begin{proof}

Inverting the spherical transform of $k_t$ and substituting $s = e$ gives

\begin{align}
\notag
\sum_{\nu \in \Lambda} d(\nu) \widehat{k}_t(\nu) & = k_t(e) \\
\notag
& = k_t^0 * k_t^0 (e) \\
\notag
& \le \| k_t^0 \|_2^2 \\
\label{transformsum}
& \ll d(t\mu) \ll t^{n-r}.
\end{align}
If we choose $\delta > 0$ and let $B_1^* \subset \ga_0^*$ be the ball of radius $\delta$ about $s^* \mu$, we may also apply Proposition \ref{kmutransform} to obtain

\begin{align*}
k_t(s) & = \sum_{\nu \in \Lambda \cap t B_1^*} d(\nu) \widehat{k}_t (\nu) \varphi_{\nu}(s) + O_A(t^{-A}).
\end{align*}
Combining this with (\ref{transformsum}) and the positivity of $\widehat{k}_t$ implies that

\begin{equation*}
|k_t(s)| \ll t^{n-r} \underset{\nu \in \Lambda \cap t B_1^*}{\sup} | \varphi_\nu (s) | + O_A(t^{-A}),
\end{equation*}
and as we may assume that $s \in \exp(B)$, the result now follows from Theorem \ref{cpctphibd}.

\end{proof}

The $L^2 \rightarrow L^2$ bound for the truncated pieces is proven by combining Lemma \ref{kbound} with the method of Proposition \ref{kmutransform}.  Theorem \ref{main} now follows as in $\mathsection$\ref{Lp}.

\bigskip
\footnotesize
\noindent\textit{Acknowledgments.}
We would like to thank Sigurdur Helgason, Peter Sarnak, Christopher Sogge, Melissa Tacy, Veeravalli Varadarajan, and Steve Zelditch for many helpful discussions, and the referee for a careful reading of the manuscript.
This research was partly supported by NSF grants DMS-1201321 and DMS-1509331.

\end{document}